\documentclass[letterpaper,11pt,oneside,reqno]{amsart}
\usepackage{bm}
\usepackage{mathrsfs}
\usepackage{amsfonts,amsmath, amssymb,amsthm,amscd,stmaryrd}
\usepackage[height=9.6in,width=5.95in]{geometry}
\usepackage[mpexclude,DIV13]{typearea}
\usepackage{verbatim}
\usepackage{hyperref}
\usepackage{graphicx}
\usepackage[latin1]{inputenc}
\usepackage{latexsym}
\usepackage{lscape}
\usepackage{epsfig}
\include{bibtex}

\usepackage{setspace}

\usepackage{parskip}
\begin{document}
\bibliographystyle{alpha}
\newcommand{\cn}[1]{\overline{#1}}
\newcommand{\e}[0]{\epsilon}
\newcommand{\bbf}[0]{\mathbf}

\newcommand{\Pfree}[5]{\ensuremath{\mathbb{P}^{#1,#2,#3,#4,#5}}}
\newcommand{\PfreeShort}{\ensuremath{\mathbb{P}^{BB}}}

\newcommand{\WH}[8]{\ensuremath{\mathbb{W}^{#1,#2,#3,#4,#5,#6,#7}_{#8}}}
\newcommand{\Wfree}[5]{\ensuremath{\mathbb{W}^{#1,#2,#3,#4,#5}}}
\newcommand{\WHShort}[3]{\ensuremath{\mathbb{W}^{#1,#2}_{#3}}}
\newcommand{\WHShortCouple}[2]{\ensuremath{\mathbb{W}^{#1}_{#2}}}

\newcommand{\walk}[3]{\ensuremath{X^{#1,#2}_{#3}}}
\newcommand{\walkupdated}[3]{\ensuremath{\tilde{X}^{#1,#2}_{#3}}}
\newcommand{\walkfull}[2]{\ensuremath{X^{#1,#2}}}
\newcommand{\walkfullupdated}[2]{\ensuremath{\tilde{X}^{#1,#2}}}

\newcommand{\PH}[8]{\ensuremath{\mathbb{Q}^{#1,#2,#3,#4,#5,#6,#7}_{#8}}}
\newcommand{\PHShort}[1]{\ensuremath{\mathbb{Q}_{#1}}}
\newcommand{\PHExp}[8]{\ensuremath{\mathbb{F}^{#1,#2,#3,#4,#5,#6,#7}_{#8}}}

\newcommand{\D}[8]{\ensuremath{D^{#1,#2,#3,#4,#5,#6,#7}_{#8}}}
\newcommand{\DShort}[1]{\ensuremath{D_{#1}}}
\newcommand{\partfunc}[8]{\ensuremath{Z^{#1,#2,#3,#4,#5,#6,#7}_{#8}}}
\newcommand{\partfuncShort}[1]{\ensuremath{Z_{#1}}}
\newcommand{\bolt}[8]{\ensuremath{W^{#1,#2,#3,#4,#5,#6,#7}_{#8}}}
\newcommand{\boltShort}[1]{\ensuremath{W_{#1}}}
\newcommand{\boltNew}{\ensuremath{W}}
\newcommand{\QTLH}{\ensuremath{\mathfrak{H}}}
\newcommand{\QTLHgen}{\ensuremath{\mathfrak{L}}}

\newcommand{\whitenoise}{\ensuremath{\mathscr{\dot{W}}}}
\newcommand{\mf}{\mathfrak}

\newcommand{\EE}{\ensuremath{\mathbb{E}}}
\newcommand{\PP}{\ensuremath{\mathbb{P}}}
\newcommand{\var}{\textrm{var}}
\newcommand{\N}{\ensuremath{\mathbb{N}}}
\newcommand{\R}{\ensuremath{\mathbb{R}}}
\newcommand{\C}{\ensuremath{\mathbb{C}}}
\newcommand{\Z}{\ensuremath{\mathbb{Z}}}
\newcommand{\Q}{\ensuremath{\mathbb{Q}}}
\newcommand{\T}{\ensuremath{\mathbb{T}}}
\newcommand{\E}[0]{\mathbb{E}}
\newcommand{\OO}[0]{\Omega}
\newcommand{\F}[0]{\mathfrak{F}}
\def \Ai {{\rm Ai}}
\newcommand{\G}[0]{\mathfrak{G}}
\newcommand{\ta}[0]{\theta}
\newcommand{\w}[0]{\omega}
\newcommand{\ra}[0]{\rightarrow}
\newcommand{\vectoro}{\overline}
\newcommand{\crairy}{\mathcal{CA}}
\newcommand{\nc}{\mathsf{NoTouch}}
\newcommand{\ncf}{\mathsf{NoTouch}^f}
\newcommand{\wxy}{\mathcal{W}_{k;\bar{x},\bar{y}}}
\newcommand{\AP}{\mathfrak{a}}
\newcommand{\cm}{\mathfrak{c}}
\newtheorem{theorem}{Theorem}[section]
\newtheorem{partialtheorem}{Partial Theorem}[section]
\newtheorem{conj}[theorem]{Conjecture}
\newtheorem{lemma}[theorem]{Lemma}
\newtheorem{proposition}[theorem]{Proposition}
\newtheorem{corollary}[theorem]{Corollary}
\newtheorem{claim}[theorem]{Claim}
\newtheorem{experiment}[theorem]{Experimental Result}

\def\todo#1{\marginpar{\raggedright\footnotesize #1}}
\def\change#1{{\color{green}\todo{change}#1}}
\def\note#1{\textup{\textsf{\color{blue}(#1)}}}

\theoremstyle{definition}
\newtheorem{rem}[theorem]{Remark}

\theoremstyle{definition}
\newtheorem{com}[theorem]{Comment}

\theoremstyle{definition}
\newtheorem{definition}[theorem]{Definition}

\theoremstyle{definition}
\newtheorem{definitions}[theorem]{Definitions}

\theoremstyle{definition}
\newtheorem{conjecture}[theorem]{Conjecture}

\newcommand{\airysh}{\mathcal{A}}
\newcommand{\hfixed}{\mathcal{H}}
\newcommand{\afixed}{\mathcal{A}}
\newcommand{\canopynoarg}{\mathsf{C}}
\newcommand{\canopy}[3]{\ensuremath{\mathsf{C}_{#1,#2}^{#3}}}
\newcommand{\argmax}{x_{{\rm max}}}
\newcommand{\zmax}{z_{{\rm max}}}

\newcommand{\Rkle}{\ensuremath{\mathbb{R}^k_{>}}}
\newcommand{\Ronele}{\ensuremath{\mathbb{R}^k_{>}}}
\newcommand{\ewxy}{\mathcal{E}_{k;\bar{x},\bar{y}}}

\newcommand{\bxyf}{\mathcal{B}_{\bar{x},\bar{y},f}}
\newcommand{\bxyflr}{\mathcal{B}_{\bar{x},\bar{y},f}^{\ell,r}}

\newcommand{\bxyfone}{\mathcal{B}_{x_1,y_1,f}}

\newcommand{\Cpop}{G}
\newcommand{\cpop}{g}

\newcommand{\ptac}{p}
\newcommand{\ptact}{v}

\newcommand{\fext}{\mathfrak{F}_{{\rm ext}}}
\newcommand{\gext}{\mathfrak{G}_{{\rm ext}}}
\newcommand{\xext}{{\rm xExt}(\mathfrak{c}_+)}

\newcommand{\dd}{\, {\rm d}}
\newcommand{\signc}{\Sigma}
\newcommand{\wxylr}{\mathcal{W}_{k;\bar{x},\bar{y}}^{\ell,r}}
\newcommand{\wxylrprime}{\mathcal{W}_{k;\bar{x}',\bar{y}'}^{\ell,r}}
\newcommand{\Rklezero}{\ensuremath{\mathbb{R}^k_{>0}}}
\newcommand{\XYfM}{\textrm{XY}^{f}_M}

\newcommand{\upright}{SC}
\newcommand{\staircase}{SC}
\newcommand{\energy}{E}
\newcommand{\xmax}{{\rm max}_1}
\newcommand{\ymax}{{\rm max}_2}
\newcommand{\lppls}{\mathcal{L}}
\newcommand{\lpplsre}{\mathcal{L}^{{\rm re}}}
\newcommand{\lpplsarg}[1]{\mathcal{L}_{n}^{\fa \to #1}}
\newcommand{\larg}[3]{\mathcal{L}_{n}^{#1,#2;#3}}
\newcommand{\BP}{M}
\newcommand{\weight}{\mathsf{Wgt}}
\newcommand{\pairweight}{\mathsf{PairWgt}}
\newcommand{\sumweight}{\mathsf{SumWgt}}
\newcommand{\mpgood}{\mathcal{G}}
\newcommand{\mpg}{\mathsf{Fav}}
\newcommand{\mcgone}{\mathsf{Fav}_1}
\newcommand{\radnik}[2]{\mathsf{RN}_{#1,#2}}
\newcommand{\size}[2]{\mathsf{S}_{#1,#2}}
\newcommand{\pdr}{\mathsf{PolyDevReg}}
\newcommand{\pwr}{\mathsf{PolyWgtReg}}
\newcommand{\lwr}{\mathsf{LocWgtReg}}
\newcommand{\hwp}{\mathsf{HighWgtPoly}}
\newcommand{\fbr}{\mathsf{ForBouqReg}}
\newcommand{\bbr}{\mathsf{BackBouqReg}}
\newcommand{\fsc}{\mathsf{FavSurCon}}
\newcommand{\maxpoly}{\mathrm{MaxDisjtPoly}}
\newcommand{\maxswf}{\mathsf{MaxScSumWgtFl}}
\newcommand{\emaxswf}{\e \! - \! \maxswf}
\newcommand{\minswf}{\mathsf{MinScSumWgtFl}}
\newcommand{\eminswf}{\e \! - \! \minswf}
\newcommand{\surreg}{\mathcal{R}}
\newcommand{\scf}{\mathsf{FavSurgCond}}
\newcommand{\disjtpoly}{\mathsf{DisjtPoly}}
\newcommand{\intint}[1]{\llbracket 1,#1 \rrbracket}
\newcommand{\maxsym}{*}
\newcommand{\polynum}{\#\mathsf{Poly}}
\newcommand{\dlp}{\mathsf{DisjtLinePoly}}
\newcommand{\lowb}{\underline{B}}
\newcommand{\highb}{\overline{B}}
\newcommand{\tottt}{t_{1,2}^{2/3}}
\newcommand{\tot}{t_{1,2}}
\newcommand{\btone}{{\bf{t}}_1}
\newcommand{\bttwo}{{\bf{t}}_2}
\newcommand{\formerE}{C}
\newcommand{\rcon}{r_0}
\newcommand{\para}{Q}
\newcommand{\Cstrong}{E}

\newcommand{\mc}{\mathcal}
\newcommand{\vect}{\mathbf}
\newcommand{\bt}{\mathbf{t}}
\newcommand{\scB}{\mathscr{B}}
\newcommand{\scBres}{\mathscr{B}^{\mathrm{re}}}
\newcommand{\rightshadow}{\mathrm{RS}Z}
\newcommand{\dbm}{D}
\newcommand{\edgedbm}{D^{\rm edge}}
\newcommand{\gue}{\mathrm{GUE}}
\newcommand{\edgegue}{\mathrm{GUE}^{\mathrm{edge}}}
\newcommand{\eqdist}{\stackrel{(d)}{=}}
\newcommand{\geqdist}{\stackrel{(d)}{\succeq}}
\newcommand{\leqdist}{\stackrel{(d)}{\preceq}}
\newcommand{\scal}{{\rm sc}}
\newcommand{\fa}{x_0}
\newcommand{\hit}{H}
\newcommand{\scaledle}{\mathsf{Sc}\mc{L}}
\newcommand{\cenleup}{\mathscr{L}^{\uparrow}}
\newcommand{\cenledown}{\mathscr{L}^{\downarrow}}
\newcommand{\eln}{T}
\newcommand{\xmin}{{\rm Corner}^{\mfl,\mc{F}}}
\newcommand{\ymin}{{\rm Corner}^{\mfr,\mc{F}}}
\newcommand{\barxmin}{\overline{\rm Corner}^{\mfl,\mc{F}}}
\newcommand{\barymin}{\overline{\rm Corner}^{\mfr,\mc{F}}}
\newcommand{\qmin}{Q^{\mc{F}^1}}
\newcommand{\barqmin}{\bar{Q}^{\mc{F}^1}}
\newcommand{\test}{T}
\newcommand{\mfl}{\mf{l}}
\newcommand{\mfr}{\mf{r}}
\newcommand{\gfl}{\ell}
\newcommand{\gfr}{r}
\newcommand{\jre}{J}
\newcommand{\highfl}{{\rm HFL}}
\newcommand{\flyleap}{\mathsf{FlyLeap}}
\newcommand{\touch}{\mathsf{Touch}}
\newcommand{\notouch}{\mathsf{NoTouch}}
\newcommand{\close}{\mathsf{Close}}
\newcommand{\abovepar}{\mathsf{High}}
\newcommand{\vecint}{\bar{\iota}}
\newcommand{\cornthree}{{\rm Corner}^\mc{G}_{k,\mfl}}
\newcommand{\cornfour}{{\rm Corner}^\mc{H}_{k,\fa}}
\newcommand{\mpgg}{\mathsf{Fav}_{\mc{G}}}

\newcommand{\lefta}{M_{1,k+1}^{[-2\eln,\gfl]}}
\newcommand{\mida}{M_{1,k+1}^{[\gfl,\gfr]}}
\newcommand{\righta}{M_{1,k+1}^{[\gfr,2\eln]}}

\newcommand{\ipdval}{d}
\newcommand{\ctemp}{d_0}

\newcommand{\wien}{W}
\newcommand{\pole}{P}
\newcommand{\pp}{p}

\newcommand{\const}{D_k}
\newcommand{\numcone}{14}
\newcommand{\numctwo}{13}
\newcommand{\numcthree}{6}
\newcommand{\rsC}{C}
\newcommand{\rsc}{c}
\newcommand{\cone}{c_1}
\newcommand{\Cone}{C_1}
\newcommand{\Ctwo}{C_2}
\newcommand{\smallc}{c_0}
\newcommand{\smallcprime}{c_1}
\newcommand{\smallcanother}{c_2}
\newcommand{\smallcnew}{c_3}
\newcommand{\Cda}{D}
\newcommand{\Kzero}{K_0}
\newcommand{\Rmac}{R}
\newcommand{\rmac}{r}
\newcommand{\conseqmac}{H}
\newcommand{\constn}{C'}
\newcommand{\coninit}{\Psi}
\newcommand{\condee}{\hat{H}}
\newcommand{\conbrac}{\hat{C}}
\newcommand{\Cnew}{\tilde{C}}
\newcommand{\Cbig}{C^*}
\newcommand{\Ctbd}{C_+}
\newcommand{\Ctbs}{C_-}

\newcommand{\stitch}{s}

\newcommand{\imax}{i_{{\rm max}}}

\newcommand{\wlp}{{\rm WLP}}

\newcommand{\canopynumber}{\mathsf{Canopy}{\#}}

\newcommand{\cannum}{{\#}\mathsf{SC}}

\newcommand{\boundgood}{\mathsf{G}}
\newcommand{\lshift}{\mc{L}^{\rm shift}}
\newcommand{\deltapi}{\theta}
\newcommand{\rootneigh}{\mathrm{RNI}}
\newcommand{\rootneighuse}{\mathrm{RNI}}
\newcommand{\manycan}{\mathsf{ManyCanopy}}
\newcommand{\specialpt}{\mathrm{spec}}

\newcommand{\dist}{\vert\vert}
\newcommand{\fik}{\mc{F}_i^{[K,K+1]^c}}
\newcommand{\mcfa}{\mc{H}[\fa]}
\newcommand{\tent}{{\rm Tent}}
\newcommand{\goodk}{\mc{G}_{K,K+1}}
\newcommand{\pairsep}{{\rm PS}}
\newcommand{\mbf}{\mathsf{MBF}}
\newcommand{\nbd}{\mathsf{NoBigDrop}}
\newcommand{\bd}{\mathsf{BigDrop}}
\newcommand{\jleft}{j_{{\rm left}}}
\newcommand{\jright}{j_{{\rm right}}}
\newcommand{\smalljfluc}{\mathsf{SmallJFluc}}
\newcommand{\mfone}{M_{\mc{F}^1}}
\newcommand{\mfthree}{M_{\mc{G}}}
\newcommand{\rhomac}{P}
\newcommand{\phimac}{\varphi}
\newcommand{\chimac}{\chi}
\newcommand{\xnmac}{z_{\mathcal{L}}}
\newcommand{\Cwb}{E_0}
\newcommand{\initcond}{\mathcal{I}}
\newcommand{\neargeod}{\mathsf{NearGeod}}
\newcommand{\polyunique}{\mathrm{PolyUnique}}
\newcommand{\latecoal}{\mathsf{LateCoal}}
\newcommand{\nolatecoal}{\mathsf{NoLateCoal}}
\newcommand{\normalcoal}{\mathsf{NormalCoal}}
\newcommand{\regfluc}{\mathsf{RegFluc}}
\newcommand{\mdeltaweight}{\mathsf{Max}\Delta\mathsf{Wgt}}
\newcommand{\ovbar}[1]{\mkern 1.5mu\overline{\mkern-1.5mu#1\mkern-1.5mu}\mkern 1.5mu}

\newcommand{\maxmin}{\pwr}
\newcommand{\nmac}{N}

\newcommand{\indexset}{\mathcal{I}}

\def\fff#1{&{{\pageref{#1}}}\cr}
\def\hfff#1{\label{#1}}

\title[Regularity of polymer weight profiles]{A patchwork quilt sewn from Brownian fabric: \\
regularity of polymer weight profiles \\ in Brownian last passage percolation}

\author[A. Hammond]{Alan Hammond}
\address{A. Hammond\\
  Departments of Mathematics and Statistics\\
 U.C. Berkeley \\
  Evans Hall \\
  Berkeley, CA, 94720-3840 \\
  U.S.A.}
  \email{alanmh@berkeley.edu}
  \thanks{The author is supported by NSF grants DMS $1512908$ and $1855550$.}
  \subjclass{$82D30$ (primary); $82C22$, $82B23$, $60H15$ (secondary).}

\begin{abstract} 
In last passage percolation models lying in the KPZ universality class, the energy of long energy-maximizing paths may be studied as a function of the paths' pair of endpoint locations.
Scaled coordinates may be introduced, so that these maximizing paths, or polymers, now cross unit distances with unit-order fluctuations, and have scaled energy, or weight, of unit order. In this article, we consider Brownian last passage percolation in these scaled coordinates.
In the narrow wedge case, when one endpoint of such polymers is fixed, say at $(0,0) \in \R^2$, and the other is varied horizontally, over $(z,1)$, $z \in \R$, the polymer weight profile as a function of $z \in \R$
is locally Brownian; indeed, by \cite[Theorem~$2.11$ and Proposition~$2.5$]{BrownianReg},
 the law of the profile is known to enjoy a very strong comparison to Brownian bridge on a given compact interval, with a Radon-Nikodym derivative in every $L^p$ space for $p \in (1,\infty)$, uniformly in the scaling parameter,
provided that an affine adjustment is made to the weight profile before the comparison is made.
In this article, we generalize this narrow wedge case and study polymer weight profiles begun from a very general initial condition. We prove that the profiles on a compact interval resemble Brownian bridge in a uniform sense: splitting the compact interval into a random but controlled number of patches, the profile in each patch after affine adjustment has a Radon-Nikodym derivative that lies in every $L^p$ space for $p \in (1,3)$. This result is proved by harnessing an understanding of the uniform coalescence structure in the field of polymers developed in~\cite{NonIntPoly} using techniques from~\cite{BrownianReg} and~\cite{ModCon}.

\end{abstract}

\maketitle

\tableofcontents

\section{Introduction}

\subsection{KPZ universality}
Consider a discrete model of random growth in which the initially healthy integer lattice sites in the upper half-plane become infected. At time zero, a certain subset of such sites on the $x$-axis are infected. At any given positive integer time, one uninfected site in the upper half-plane that is a nearest neighbour of a presently infected site is selected, uniformly at random. The site then becomes infected, with a `transmission' edge being added from the newly infected site into the already infected set, this edge selected uniformly from the available possibilities.

\begin{figure}[ht]
\begin{center}
\vspace{0.6cm}
\includegraphics[height=3.3cm]{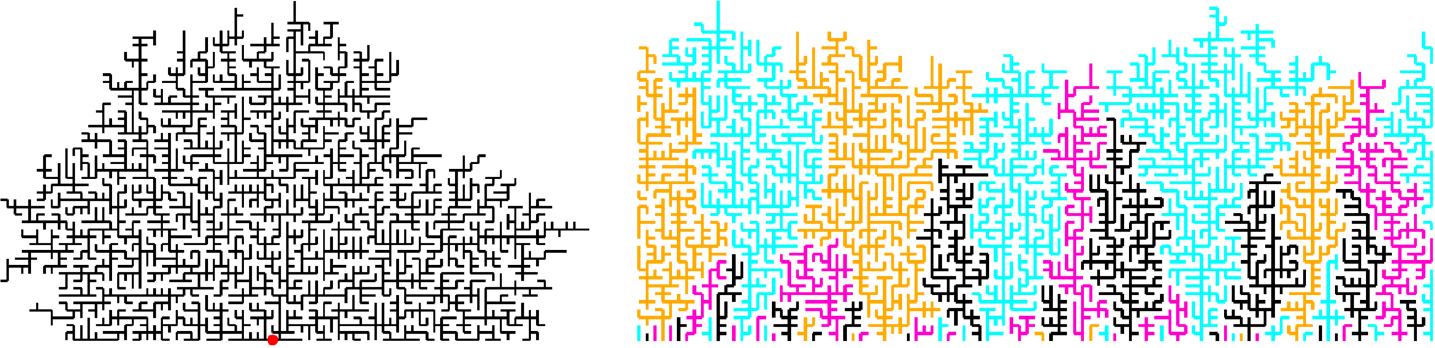}
\vspace{0.5cm}
\caption{
The Eden-like infection model growing from two initial conditions. On the left, only the origin is initially infected. The model may also be defined so that an infinite set of initially infected vertices is permitted. On the right, every integer lattice site on the $x$-axis is so infected.}
\label{f.eden}
\end{center}
\vspace{-0.1cm}
\end{figure}

In this way, the infected region grows, one site at a time. 
As Figure~\ref{f.eden} illustrates,  this region  is at any given moment the collection of vertices abutting the present set of transmission edges. The transmission edge-set is partitioned into a collection of trees, each rooted at one of the initially infected sites.

\vspace{0.4cm}

The $1 + 1$ dimensional Kardar-Parisi-Zhang (KPZ) universality class
includes a wide range of interface models suspended over a one-dimensional domain, in which growth in a direction normal to the surface competes with a smoothening surface tension in the presence of a local randomizing force that roughens the surface.
The infection model, which is a variant of the Eden model introduced in~\cite{Eden}, is expected to lie in the KPZ class. The boundary of the infected region is an interface.
When the infected region has diameter of order $n$, the interface height above a given point has, according to KPZ prediction, typical deviation from the mean of order $n^{1/3}$, while non-trivial correlations in this height as the horizontal coordinate is varied are encountered on scale~$n^{2/3}$. Moreover, an exponent of one-half dictates the interface's regularity, with the interface height being expected to vary between a pair of locations at distance of order at most $n^{2/3}$
on the order of the square root of the distance between these locations.

The broad range of interface models that are rigorously known or expected to lie in the KPZ universality class includes many last passage percolation models, 
in which the interface models the maximum obtainable value of paths, where a path is assigned a value by integrating over its course weights specified by a product measure random environment. 
 The scaling assertions associated to these values of $(1/3,1/2,2/3)$ have been rigorously demonstrated for only a few random growth models, each of which enjoys an integrable structure: for example, the seminal work of Baik, Deift and Johansson~~\cite{BDJ1999}
 rigorously established the one-third exponent, and moreover obtained the GUE Tracy-Widom distributional limit, for the case of Poissonian last passage percolation, while the two-thirds power law for transversal fluctuation was derived for this model by Johansson~~\cite{Johansson2000}.

 \subsection{Probabilistic problems and proof techniques for KPZ}

The theory of KPZ universality has advanced through physical insights, numerical analysis, and several techniques of integrable or algebraic origin. We will not hazard a summary of literature to support this one-sentence history, but refer to the reader to~\cite{IvanSurvey} for a KPZ survey from 2012; in fact, integrable and analytic approaches to KPZ have attracted great interest around and since that time. 
Now, it is hardly  deep or controversial to say that  many problems and models in KPZ are intrinsically random: the random interface of the infected region boundary in the two simulations in Figure~\ref{f.eden}, and the geometry of the randomly evolving forest composed of transmission edges, are two important characters. It would thus seem  valuable to approach the problems of KPZ universality from a predominately probabilistic perspective.

The present article is the culmination of a four-paper study of KPZ from such a viewpoint. The companion papers are~\cite{BrownianReg},~\cite{ModCon} and~\cite{NonIntPoly}.
Our main result, Theorem~\ref{t.unifpatchcompare}, draws on important technical ingredients from the other works in order to make a significant probabilistic inference about KPZ universality.

In order to explain the problem  in question, we first develop our discussion of KPZ growth.
A growth model in the KPZ universality class may be started at time zero from many different initial profiles. One very special case, in which growth is initiated from a unique point, goes by the name `narrow wedge'. This corresponds to the growth from a single infected site, the origin, seen in the left sketch of Figure~\ref{f.eden}.
In this case, the limiting description of the late-time interface, suitably scaled in light of the one-third and two-thirds powers and up to the subtraction of a parabola, is offered by the {\rm Airy}$_2$ process,
 which is a random  function $\mc{A}: \R \to \R$, whose finite dimensional distributions are specified by Fredholm determinants, that was introduced by~\cite{PrahoferSpohn}. The one-half power law for interface regularity is expressed by the H{\"o}lder-$1/2-$-continuity of~$\mc{A}$. It has been anticipated that $\mc{A}$ has a locally Brownian structure.
Such a statement as this may be interpreted either by taking a {\em local} limit, in which, for given $x \in \R$, the Gaussianity of $\e^{-1/2} \big( \mc{A}(x + \e) - \mc{A}(x) \big)$ is investigated after the low $\e$ limit is taken. An in essence stronger, and much more useful, assertion of Brownian structure would concern a {\em unit-order} scale: examples would be the absolute continuity of the random function $[x,x+K] \to \R$ given by $y \to \mc{A}(y) - \mc{A}(x)$ with respect to Brownian motion (whose rate is two in view of the convention of definition for $\mc{A}$) on any compact interval~$[x,x+K]$; or, stronger still, that the resulting Radon-Nikodym derivative lies in $L^p$-spaces for values of $p$ exceeding one.

Both types of   Brownian comparison have been made for the Airy$_2$ process. Finite dimensional distributional convergence  in the local limit has been proved in~\cite{Hagg}; 
this conclusion is strengthened to a functional convergence in~\cite[Theorem 3]{CatorPimentel}  by means of stochastically sandwiching the process between variants of the process at equilibrium. 
Regarding unit-order comparison, 
the process, after the subtraction of a suitable parabola, may be embedded as the top curve of an $\N$-indexed ordered ensemble of curves that are in effect mutually avoiding Brownian bridges.
Such ensembles of random curves satisfy a simple `Brownian Gibbs' resampling property which is a very valuable tool for their analysis. Indeed, the tool of  Brownian Gibbs resampling was used in~\cite{AiryLE} to prove the absolute continuity statement  mentioned in the preceding paragraph. These resampling ideas have been refined in the first article~\cite{BrownianReg} in our four-paper study to prove a result,~\cite[Theorem~$1.9$]{BrownianReg},
which asserts that the Radon-Nikodym derivative (of the Airy$_2$ process with respect to Brownian motion on a compact interval) lies in all $L^p$ spaces, for $p \in (1,\infty)$. Actually, there is  one caveat: the comparison of the two processes is made only after an affine shift, so that comparison is made with Brownian bridge. 

Growth may be initiated from a much more general initial condition than in this narrow wedge case. For example, in Figure~\ref{f.eden}(right), each site on the $x$-axis is initially infected. This corresponds to growth from a flat (or zero) initial condition. In this case, it is the {\rm Airy}$_1$ process which putatively describes the interface at late time. The H{\"o}lder-$1/2-$-continuity of this process has been proved by~\cite{QR13} by means of a description of its finite dimensional distributions via Fredholm determinants with kernels in $L^2(\R)$; the case of the Airy$_2$ process is also treated.  

For initial conditions that grow at most linearly, it has been anticipated that a limiting description of the suitably scaled late-time interface should also exist. Indeed, in a recent preprint~\cite{MQR17}, Matetski, Quastel and Remenik
have utilized a biorthogonal ensemble representation found by~\cite{Sas05,BFPS07} associated to the totally asymmetric exclusion process in order to  find Fredholm determinant formulas for the multi-point distribution of the height function of this growth process begun from an arbitrary initial condition. Using these formulas to take the KPZ scaling limit, the authors construct a scale invariant Markov process that lies at the heart of the KPZ universality class. The time-one evolution of this Markov process may be applied to such general initial data as we have mentioned, and the result is the scaled profile begun from such data, which generalizes the Airy$_2$ process seen in the narrow wedge case. It is natural to ask what form of local regularity this profile enjoys: is it H\"older-$1/2-$ continuous, and what comparison to Brownian motion can be made? Theorem~$4.4$ in~\cite{MQR17}
 asserts such H\"older continuity, and also makes a local limit Brownian comparison, proving that there is convergence to Brownian motion  (with rate two) in finite dimensional distributions in this limit.
Such local limit results for general initial condition profiles have also been derived by~\cite{Pimentel17}, for geometric last passage percolation models. 

The distinction between asserting Brownian structure in a local limit, and doing so on a unit-order scale, is an important one.
The aim of proving~\cite[Conjecture~$1.5$]{Johansson2003} illustrates the difference. Johansson's conjecture states that  the Airy$_2$ process after subtraction of the parabola~$x^2$ has a unique maximizer; it is important because the maximizer describes the scaled location of the endpoint of a point-to-line maximizing path in last passage percolation models.
Local limit Gaussianity does not rule out the presence of two or more maximizers, but the result is a direct consequence~\cite[Theorem~$4.3$]{AiryLE} of unit-order Brownian comparison for the  Airy$_2$ process.  
(The conjecture in fact has several proofs:  Moreno Flores, Quastel and Remenik~\cite{MFQR13} via an explicit formula for the maximizer, and an argument of  Pimentel~\cite{Pimentel14} showing that any stationary process minus a parabola has a unique maximizer.) Moreover, stronger assertions of unit-order Brownian structure, involving finiteness for higher $L^p$-norms of the Radon-Nikodym derivative,
imply that Brownian motion (or bridge) characteristics obtain on the unit scale in a very strong sense -- see \cite[Theorem~$1.10$]{BrownianReg}.

The problem of unit-order scale, rather than local limit, Brownian comparison for scaled height profiles begun from general initial data, is an important example of an intrinsically probabilistic question in the theory of KPZ universality.
In this article, we study it, for
almost arbitrary initial data. 

An idea at the heart of our approach is represented by the system of trees seen in the right sketch of Figure~\ref{f.eden}. 
(This is only a conceptual connection: the infection model illustrates ideas, and is not being investigated here.)
As the trees grow, they compete,
with the interface at any given time being partitioned into canopies of presently surviving trees. One important aspect of forest geometry is that it may be expected to respect KPZ scaling, so that a tree surviving when growth has advanced on order~$n$ has a canopy of length of order $n^{2/3}$. In the companion papers of our four-paper study, such a scaling for forest geometry has been proved.
Here we exploit this understanding to study scaled interfaces.

Although our approach is predominately probabilistic, it does rely crucially on certain limited integrable inputs (which are lacking for example in the infection model).  The model that we choose for study is Brownian last passage percolation, a model in the KPZ universality class that enjoys attractive probabilistic (and integrable) features.
Our principal conclusion, Theorem~\ref{t.unifpatchcompare}, makes  a strong Brownian comparison for scaled interface profiles  uniformly over all high choices of the length scale parameter for these microscopic models. 

We now prepare to state this principal conclusion, first specifying the model under study.



\subsection{Brownian last passage percolation}\label{s.lpp}
We will call this model Brownian LPP.
On a probability space carrying a law labelled~$\PP$,
let $B:\Z \times \R \to \R$ denote an ensemble of independent  two-sided standard Brownian motions $B(k,\cdot):\R\to \R$, $k \in \Z$.

Let $i,j \in \Z$ with $i \leq j$.
We denote the integer interval $\{i,\cdots,j\}$ by $\llbracket i,j \rrbracket$.
Further let $x,y \in \R$ with $x \leq y$.
Consider the collection of  non-decreasing lists 
 $\big\{ z_k: k \in \llbracket i+1,j \rrbracket \big\}$ of values $z_k \in [x,y]$. 
With the convention that $z_i = x$ and $z_{j+1} = y$,
we associate an energy $\sum_{k=i}^j \big( B ( k,z_{k+1} ) - B( k,z_k ) \big)$ to any such list.
We then define  the maximum energy
$$
M^1_{(x,i) \to (y,j)} \, = \, \sup \, \bigg\{ \, \sum_{k=i}^j \Big( B ( k,z_{k+1} ) - B( k,z_k ) \Big) \, \bigg\} \, , 
$$
where the supremum is taken over all such lists. The random process $M^1_{(0,1) \to (\cdot,n)}: [0,\infty) \to \R$ was introduced by~\cite{GlynnWhitt} and further studied in~\cite{Baryshnikov},~\cite{GTW} and~\cite{O'ConnellYor}.


The one-third and two-thirds KPZ scaling considerations that we have outlined are manifest in Brownian LPP. 
When the ending height $j$ exceeds the starting height $i$ by a large positive integer $n$, and the location $y$ exceeds $x$ also by $n$, then the maximum energy grows linearly, at rate $2n$,
and has a fluctuation about this mean of order $n^{1/3}$. Moreover, if $y$ is permitted to vary from this location, then it is changes of $n^{2/3}$ in its value that result in a non-trivial correlation of the maximum energy from its original value.

These facts prompt us to introduce scaled coordinates to describe the two endpoint locations, and a notion of scaled maximum energy, which we will refer to as weight. 
Let  $n$ be an element in the {\em positive} integers $\N$, and 
suppose that $x,y \in \R$ satisfy 
 $y \geq x - 2^{-1} n^{1/3}$.
Define
\begin{equation}\label{e.weightmzeroone} 
  \weight_{n;(x,0)}^{(y,1)} \,     =  \,   2^{-1/2} n^{-1/3} \Big(  M^1_{(2n^{2/3}x,0) \to (n  + 2n^{2/3}y,n)} - 2n  -  2n^{2/3}(y-x) \Big) \, .
\end{equation}
Consistently with the facts just mentioned, the quantity  $\weight_{n;(x,0)}^{(y,1)}$ may be expected to be, for given real choices of $x$ and $y$, a unit-order random quantity, whose law is tight in the scaling parameter $n \in \N$. The quantity describes, in units chosen to achieve this tightness, the maximum possible energy associated to journeys which in the original coordinates occur between  $(2n^{2/3}x,0)$ and $(n  + 2n^{2/3}y,n)$.
In scaled coordinates, this is a journey between $(x,0)$ and $(y,1)$.
We view the first coordinate as space and the second as time, so this journey is between $x$ and $y$ over the unit time interval $[0,1]$.

 Underlying this definition is a geometric picture of scaled maximizing paths, or polymers, that achieve these weight values. This picture will be central to our study, and we will explain it shortly, in a way that may serve to further explain the above definition.

\subsection{Polymer weight profiles from general initial data}

For now, we continue on a rather direct route to stating our principal conclusion. The random function $y \to   \weight_{n;(0,0)}^{(y,1)}$
may be viewed as the weight profile obtained by scaled maximizing paths that travel from the origin at time zero to the variable location $y$ at time one. This insistence that the paths must begin at the origin  (the case that is called the narrow wedge) is of course rather special. 
We now make a more general definition, of the  $f$-rewarded line-to-point polymer weight  $\weight_{n;(*:f,0)}^{(y,1)}$. Here, $f$ is an initial condition, defined on the real line. 
Paths may begin anywhere on the real line at time zero; they travel to $y \in \R$ at time one. (Because they are free at the beginning and fixed at the end, we refer to these paths as `line-to-point'.) They begin with a reward given by evaluating $f$ at the starting location, and then gain the weight associated to the journey they make.  The value $\weight_{n;(*:f,0)}^{(y,1)}$, which we will define momentarily, denotes the maximum $f$-rewarded weight of all such paths.   In the notation $\weight_{n;(*:f,0)}^{(y,1)}$, we again use subscript and superscript expressions to refer to space-time pairs of starting and ending locations. The starting spatial location is being denoted $*:f$. The star is intended to refer to the free time-zero endpoint, which may be varied, and the $:f$ to the reward offered according to where this endpoint is placed.

A very broad class of initial conditions $f$
is next specified by an upper bound of uniform linear growth and a lower bound of non-degeneracy in a compact interval. This function class is
  suitable for a study of the weight profiles $y \to \weight_{n;(*:f,0)}^{(y,1)}$ for all sufficiently high~$n \in \N$.

\begin{definition}\label{d.if}
Writing  $\ovbar\coninit = \big( \coninit_1,\coninit_2,\coninit_3 \big) \in (0,\infty)^3$ for a triple of positive reals, we let $\initcond_{\ovbar\coninit}$~\hfff{initcond}
denote the set of measurable functions $f:\R \to \R \cup \{ - \infty \}$ such that
$f(x) \leq \coninit_1 \big( 1 + \vert x \vert \big)$
and $\sup_{x \in [-\coninit_2,\coninit_2]} f(x) > - \coninit_3$.
\end{definition}

For $f$ lying in one of the function spaces $\initcond_{\ovbar\coninit}$, we now formally define the $f$-rewarded line-to-point polymer weight  $\weight_{n;(*:f,0)}^{(y,1)}$~\hfff{fweight} according to  
$$
 \weight_{n;(*:f,0)}^{(y,1)} \,  = \,  \sup \, \Big\{ \,  \weight_{n;(x,0)}^{(y,1)}    + f(x) : x \in (-\infty,2^{-1}n^{1/3} + y] \, \Big\}  \, .
$$

Our principal conclusion, Theorem~\ref{t.unifpatchcompare}, roughly asserts that the weight profiles  $y \to  \weight_{n;(*:f,0)}^{(y,1)}$, viewed as functions of $y$ in a given compact real interval, 
enjoy a uniformly strong similarity with Brownian motion, even as the parameters $n$ and $f$ are permitted to vary over all sufficiently high integer values and over the function space $\initcond_{\ovbar\coninit}$. The phrase `uniformly strong similarity' is a simplification, however, and, in the next section, we explain what form of comparison we will make.

\subsection{Patchwork quilts sewn from Brownian pieces of fabric}\label{s.patchwork}

In essence, our theorem will assert that each weight profile is a Brownian patchwork quilt, with a uniform efficiency in manufacture.
Now we define what we mean by this new concept.

\subsubsection{Continuous paths, bridges, and the projection between them}\label{s.bridgeproject}

For $a,b \in \R$ with $a \leq b$, let $\mc{C}_{*,*} \big( [a,b], \R \big)$ denote the space of continuous real-valued functions of the interval~$[a,b]$.
The use of the pair of subscript stars is intended to indicate that there is no restriction placed on the endpoint values of the member functions.
An element $f$ of  $\mc{C}_{*,*} \big( [a,b], \R \big)$ that vanishes at the endpoints of $[a,b]$ is here called a bridge. Denote by $\mc{C}_{0,0} \big( [a,b], \R \big)$  the collection of bridges.
A natural projection  maps the first space onto the second, sending a continuous function $f:[a,b] \to \R$  to the bridge
$f(x) - (b-a)^{-1} \big(  (b-x)f(a) + (x-a) f(b) \big)$. The latter bridge will be denoted by $f^{[a,b]}:[a,b] \to \R$. The notation extends to stochastic processes. If $X$ is a $\mc{C}_{*,*} \big( [a,b], \R \big)$-valued random process defined under the law $\PP$, then its bridge projection $X^{[a,b]}$ is a  $\mc{C}_{0,0} \big( [a,b], \R \big)$-valued process under the same measure.
Any law $\nu$ on  $\mc{C}_{*,*} \big( [a,b], \R \big)$ naturally induces a push-forward law on  $\mc{C}_{0,0} \big( [a,b], \R \big)$.
When~$\nu$ is the law of standard Brownian motion $X:[a,b] \to \R$ with $X(a) = 0$, the bridge-valued push forward measure, which is the law of standard Brownian 
bridge, will be denoted by $\mc{B}_{0,0}^{[a,b]}$.~\hfff{bridge}

\subsubsection{Patchwork quilts}

How to sew a patchwork quilt? By cutting several pieces of fabric and stitching them together.

\begin{figure}[ht]
\begin{center}
\includegraphics[height=7cm]{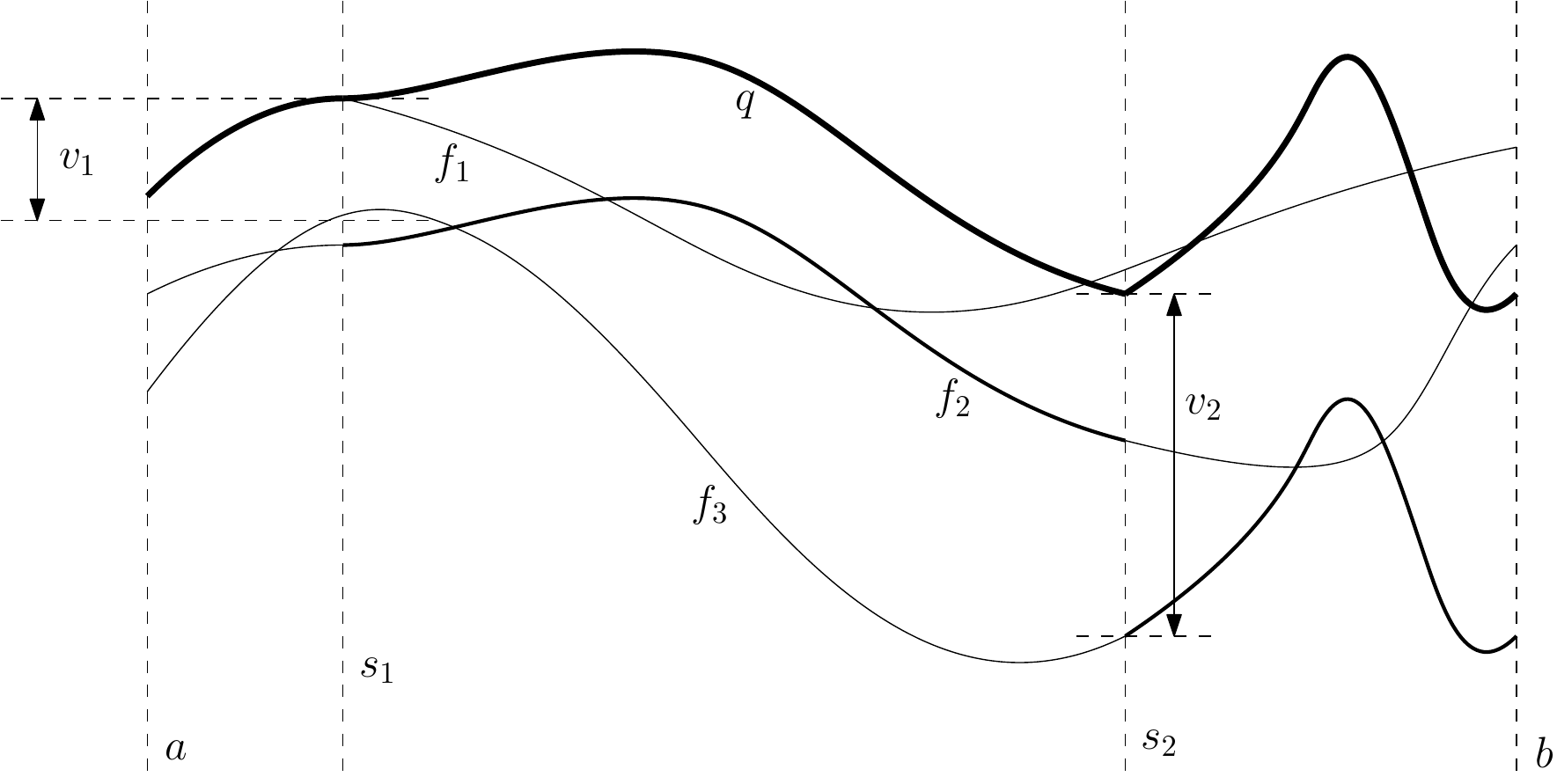}
\caption{A quilt $q:[a,b] \to \R$ is formed from three pieces of fabric.}
\label{f.quilt}
\end{center}
\end{figure}

Let $a,b \in \R$ with $a < b$ and 
let $k \in \N$. We aim to form a patchwork quilt, which will be a continuous function $q:[a,b] \to \R$. Our raw material consists of $k$ pieces of fabric, each a continuous real-valued function on $[a,b]$. This data is listed as the {\em fabric} sequence:
$f_i: [a,b] \to \R$ for $1 \leq i \leq k$.
It is our intention to sew together the consecutive pieces of fabric at a given set of $k - 1$ locations in $[a,b]$. We call these locations the {\em stitch} points, and list them in increasing order $a \leq \stitch_1 < \stitch_2 < \cdots < \stitch_{k-1} \leq b$. In this way, the interval $[a,b]$ is divided into $k$ {\em patches}. With the convention that $a = \stitch_0$ and $b = \stitch_k$, the $i\textsuperscript{th}$ patch is $[\stitch_{i-1},\stitch_i]$.
The $i\textsuperscript{th}$ piece of fabric is cut at the two ends of the $i\textsuperscript{th}$ patch, and the resulting pieces are displaced vertically, in order that they meet at their endpoints, and  stitched together.
The resulting {\em patchwork quilt} is the continuous function $q:[a,b] \to \R$ given by
$q(x) = f_i(x) + v_i$ for $x \in [\stitch_{i-1},\stitch_i]$ and $1 \leq i \leq k$.
Here, $v_i$, $1 \leq i \leq k$, are the vertical shifts: we demand that $v_0 = 0$ (so that $q(a) = f_1(a)$), and determine the later $v_i$ values by insisting on the continuity of $q$. See Figure~\ref{f.quilt}.

\subsubsection{Random continuous curves that uniformly withstand comparison to Brownian bridge}
Let $X:[a,b] \to \R$ denote a random continuous function defined on a probability space carrying the law~$\PP$; thus, $X$ is a $\mc{C}_{*,*} \big( [a,b], \R \big)$-valued random variable. The process 
$X^{[a,b]}$ is thus a random bridge on $[a,b]$ whose law may be compared to that of Brownian bridge.
For $\beta \in (0,\infty)$, we say that $X$ {\em withstands $L^{\beta-}$-comparison to Brownian bridge} if the Radon-Nikodym derivative of the law of $X^{[a,b]}$
with respect to $\mc{B}_{0,0}^{[a,b]}$
has a finite $L^\eta$-norm, for every $\eta \in (0,\beta)$. We may reexpress this condition in terms of the {\em deformation in probability of unlikely events}:
to this end, abbreviate $\mc{B} = \mc{B}_{0,0}^{[a,b]}$, and write $\mu$ for the law of $X$.  
The rephrased assertion is that, for any $\eta \in (0,1 - \beta^{-1})$, there exists an $\eta$-dependent constant $C_0$ such that, for any measurable subset $A$ of  $\mc{C}_{*,*} \big( [a,b], \R \big)$,
$\mu(A) \leq C_0 \big( \mc{B}(A) \big)^{1 - \beta^{-1} - \eta}$.

 Suppose that we instead consider a collection of random continuous functions, each defined under the law~$\PP$.
 The collection is said to {\em uniformly} withstand  $L^{\beta-}$-comparison to Brownian bridge if, 
 for any given  $\eta \in (0,1 - \beta^{-1})$,
the bound 
$\mu(A) \leq C_0 \big( \mc{B}(A) \big)^{1 - \beta^{-1} - \eta}$ holds  for any measurable subset $A$ of  $\mc{C}_{*,*} \big( [a,b], \R \big)$ and for every law $\mu$ of these random functions, with the $\eta$-dependent constant $C_0$
chosen independently of $\mu$. 
 
Suppose further that the index set of the random continuous functions takes the form $\N \times \mc{K}$,
so that a collection of sequences of such functions, indexed by $\mc{K}$, is being considered. Let $h:\N \to (0,1]$ be given. The collection is said to uniformly withstand   $L^{\beta-}$-comparison to Brownian bridge   {\em above scale~$h$} if, writing $\mu_{n,\kappa}$
for the law of the random function indexed by $(n,\kappa) \in \N \times \mc{K}$, we have that, for any given  $\eta \in (0,1 - \beta^{-1})$,
the bound 
$\mu_{n,\kappa}(A) \leq C_0 \big( \mc{B}(A) \big)^{1 - \beta^{-1} - \eta}$ holds  for any measurable subset $A$ of  $\mc{C}_{*,*} \big( [a,b], \R \big)$ that satisfies $\mc{B}(A) \geq h(n)$; here, the $\eta$-dependent constant $C_0$
may be chosen independently of $(n,\kappa)$. When $h$ decays rapidly to zero -- as it will in our application -- this definition varies the preceding one merely by eliminating from consideration sets whose $\mc{B}$-probability decays rapidly in the high~$n$ limit.

\subsubsection{Sewing a random patchwork quilt} 
On a probability space with measure $\PP$,
suppose given a {\em fabric} sequence $\big\{ F_i: i \in \N \big\}$ of random continuous functions defined on $[a,b]$.
Suppose further that an almost surely finite random {\em stitch points}
set $S$ satisfying  $S \subset [a,b]$ is defined on this probability space.

The patchwork quilt formed by this fabric sequence and stitch point set is itself a random continuous function, defined under the law~$\PP$, that we will denote by ${\rm Quilt}[\overline{F},S]$, (where the bar notation indicates a vector, in this case, indexed by $\N$). At almost every sample point in the probability space, $S$ is finite; at such a point, the quilt is formed from the first $\vert S \vert + 1$ fabric sequence elements, and the later elements play no role.

\subsubsection{When random functions may be formed as Brownian patchwork quilts with uniform efficiency} 
Let $\indexset$ be an arbitrary index set. 
Suppose given a collection of random continuous functions $X_{n,\alpha}:[a,b] \to \R$, indexed by $(n,\alpha) \in \N \times \indexset$, defined under the law~$\PP$. 

 Let $\beta_1 > 0$, $\beta_2 \geq 1$, $\beta_3 > 0$ and $\beta_4 > 0$. This collection is said to be uniformly Brownian patchwork $(\beta_1,\beta_2,\beta_3,\beta_4)$-quiltable if
 \begin{itemize}
 \item
 there exist  sequences $p,q :\N \to[0,1]$ 
verifying $p_j \leq j^{-\beta_1 + \e}$ and $q_j \leq j^{-\beta_3 + \e}$ for each $\e > 0$ and all $j$ sufficiently high
\item and a constant $g > 0$  
 \end{itemize}
such that, for each $(n,\alpha) \in \N \times \indexset$, we may construct under the law $\PP$,
\begin{enumerate}
\item an error event $E_{n,\alpha}$ that satisfies $\PP \big( E_{n,\alpha} \big) \leq q_n$;
\item a {\em fabric} sequence $\overline{F}_{n,\alpha} = \big\{ F_{n,\alpha;i} : i \in \N \big\}$ (consisting of continuous random functions on $[a,b]$), where the 
collection $\big\{ F_{n,\alpha;i} : \big(n,(\alpha,i) \big) \in \N \times \mc{K} \big\}$, with $\mc{K} = \indexset \times \N$,   uniformly withstands  $L^{\beta_2-}$-comparison to Brownian bridge above scale $\exp \big\{ - g n^{\beta_4} \big\}$; 
\item a {\em stitch points} set $S_{n,\alpha} \subset [a,b]$ whose cardinality verifies
$\PP \big( \vert S_{n,\alpha} \vert \geq \ell \big) \leq p_\ell$ for each $\ell \in \N$;
\item and all this in such a way that, for every $(n,\alpha) \in \N \times \indexset$, the random function  $X_{n,\alpha}$ is equal to the patchwork quilt ${\rm Quilt}[\overline{F}_{n,\alpha},S_{n,\alpha}]$ throughout the interval $[a,b]$, whenever the error event $E_{n,\alpha}$ does {\em not} occur.
\end{enumerate}
(The reason for the order of labelling of the four $\beta$-values is not apparent from this definition. High values of $\beta_1$ correspond to strong control on the size of the stitch points set; high values of $\beta_2$, to strong similarity of fabric functions to Brownian bridge. We regard the degree to which these properties are controlled as the more significant aspect in the assertion of the upcoming Theorem~\ref{t.unifpatchcompare}; while $\beta_3$ and $\beta_4$ are concerned with rare, and conceptually rather insignificant, error events.)

Uniformly sewn Brownian patchwork quilts manifest in a reasonably strong sense the notion of unit-order scale Brownian comparison that we discussed in the article's opening paragraphs.
A simple example shows, however, that more than a merely abstract tool must be used to remove a given stitch from a quilt.
Indeed, consider standard Brownian motion $B:[0,1] \to \R$, and let $t:[0,1] \to [0,1]$ denote the tent map that affinely interpolates the points $t(0) = 0$, $t(1/4) = 1$, $t(1/2) = 0$ and $t(1) = 0$.
Let $H$ denote the random value such that the modified process $B + tH:[0,1] \to \R$
achieves the same maximum value on $[0,1/2]$ as it does on $[1/2,1]$. 
The modified process is absolutely continuous with respect to a suitable vertical shift of Brownian motion on each of the intervals $[0,1/2]$ and $[1/2,1]$.
However, because $B$ has an almost surely unique maximizer on $[0,1]$, whereas $B + tU$ almost surely has two maximizers, the modified process is singular with respect to Brownian motion on $[0,1]$.
That is, there is a quilt description of the modified process with a stitch sewn at one-half, and there is no means of undoing that stitch. 

\subsection{The main result: a generic weight profile is a Brownian patchwork quilt}

\begin{theorem}\label{t.unifpatchcompare}
Let $\ovbar\coninit \in (0,\infty)^3$ satisfy $\coninit_2 \geq 1$. 
The collection of random continuous functions
$$
[-1,1] \to \R : y \to \weight_{n;(*:f,0)}^{(y,1)}
$$
 indexed by 
$(n,f) \in \N \times  \initcond_{\ovbar\coninit}$
is  uniformly Brownian patchwork $\big(2,3,1/252,1/12\big)$-quiltable. 
\end{theorem}

In order to summarise what has been achieved in proving this theorem, and in what sense further progress may be possible, it is useful to formulate a conjecture about an even stronger Brownian regularity of the weight profiles.
For $a \leq b$, we may write $\mc{C}_{0,*} \big( [a,b], \R \big)$ for the space of continuous $f:[a,b] \to \R$ with $f(a) = 0$. 
We write  $\mc{B}_{0,*}^{[a,b]}$ for the law of standard Brownian motion $B:[a,b] \to \R$, with $B(a) = 0$.
Let $X:[a,b] \to \R$ be a $\mc{C}_{0,*} \big( [a,b], \R \big)$-valued random variable. 
For $\beta \in (0,\infty)$, we say that $X$ {\em withstands $L^{\beta-}$-comparison to Brownian motion} if the Radon-Nikodym derivative of the law of $X$
with respect to $\mc{B}_{0,*}^{[a,b]}$
has a finite $L^\eta$-norm, for every $\eta \in (0,\beta)$.
When a collection of such processes~$X$ is instead considered,
 the collection is said to {\em uniformly} withstand  $L^{\beta-}$-comparison to Brownian motion if the above $L^\eta$-norm, for any given  $\eta \in (0,\beta)$, is bounded above by a finite quantity that is independent of the choice of $X$ in the collection.

\begin{conjecture}\label{c.bm}
Let $\ovbar\coninit \in (0,\infty)^3$ satisfy $\coninit_2 \geq 1$.
There exists $n_0 \in \N$ such that
the collection
$$
[-1,1] \to \R : y \to \weight_{n_0 + n;(*:f,0)}^{(y,1)} -  \weight_{n_0 + n;(*:f,0)}^{(-1,1)} \, ,
$$
 indexed by 
$(n,f) \in \N \times  \initcond_{\ovbar\coninit}$, uniformly withstands $L^{\infty-}$-comparison to Brownian motion. 
\end{conjecture}

If the conjecture is to be believed, then the use of patches, and affine-shifting to obtain bridges, could ultimately be dispensed with.
 Theorem~\ref{t.unifpatchcompare}
is a significant advance in understanding because 
\begin{itemize}
\item
it captures Brownian regularity across a very wide range of initial conditions;
\item
it does so without taking a local limit, under which Gaussianity would arise from a random function $W$ in a limit of low $\e$ for the scaled process $\e^{-1/2} \big( W(x+\e) - W(x) \big)$;
\item and it provides a stronger comparison to Brownian behaviour than absolute continuity statements, in which only the finiteness of the $L^1$-norm is proved.
\end{itemize}
The concepts that drive the proof of Theorem~\ref{t.unifpatchcompare}, including the rigorous tools that have been developed to derive the result, are a critical aspect of the theorem's significance. 
A guiding theme of our approach is to rely on algebraic inputs only in a very limited way -- for example,  in order to gain control of narrow wedge profiles at given points --
and to harness  probabilistic techniques in order to reach far stronger conclusions about profiles.  (This theme has something in common with the approach recently used    in~\cite{SlowBondSol} and~\cite{SlowBond} to resolve the slow bond conjecture for the totally asymmetric exclusion process.)
In~\cite{BrownianReg}, it has been understood, using probabilistic resampling techniques, that narrow wedge profiles closely resemble Brownian motion. Our task here is to relate much more general profiles to these special ones. Theorem~\ref{t.unifpatchcompare} will be derived by studying, and proving natural properties of, an important polymer forest structure that is associated to the problem of $f$-rewarded line-to-point polymers. (The trees in the forest roughly correspond to narrow wedge profiles, similarly to the way that the surviving trees in the right sketch of Figure~\ref{f.eden} are rooted at certain points on the $x$-axis, so that each one may be  compared to the single, narrow wedge case, tree depicted on the left.) As such, we view Theorem~\ref{t.unifpatchcompare} as an important practical and conceptual step towards Conjecture~\ref{c.bm}.
How close is the conjecture, given the theorem? 
All of the stitch points should be removed from the quilts, and the affine adjustment that is made, for the purpose of Brownian comparison, to the (putatively unique) fabric sequence element   should also be eliminated. These are genuine technical challenges, in whose resolution  the structure of the polymer forest may again have a role to play. Whatever the degree of difficulty these challenges may pose, the techniques leading to Theorem~\ref{t.unifpatchcompare} have built a probabilistic road from an algebraic point of depature into a far wider realm and as such they achieve for an interesting example the aim of broadening KPZ horizons by probabilistic means.
 
 \subsubsection{Brownian comparison, the KPZ fixed point and the Airy sheet}
 We review our principal conclusion in light of the recent construction of the KPZ fixed point in~\cite{MQR17} that we mentioned at the outset of the article.
 The authors of~\cite{MQR17} construct a limiting Markov evolution on profiles begun from a class of initial data $f$ that is very similar to the function space $\initcond_{\ovbar\coninit}$ for given 
 $\ovbar\coninit \in (0,\infty)^3$. It is tempting to say that, in our language, this object is the limiting profile $y \to \weight_{n;(*:f,0)}^{(y,1)}$ with $n = \infty$.
 However, the construction undertaken in~\cite{MQR17} at present works for totally asymmetric simple exclusion, not for Brownian LPP, so that the uniqueness of the $n=\infty$ limiting process is not yet a consequence of the theory developed by~\cite{MQR17}.
 
In fact, a recent advance, made in Dauvergne, Ortmann and Vir{\'a}g~\cite{DOV18} and assisted by~\cite{DV18}, does prove the uniqueness of this limiting process. The Airy sheet is a rich scaled structure in the KPZ universality class whose existence was mooted in~\cite{CQR2015}. In the present notation, it is the putatively unique limiting weight system $\R^2 \to \R: (x,y) \to  \weight_{\infty;(x,0)}^{(y,1)}$ -- a system that  encodes the weights of all polymers that cross the unit strip $\R \times [0,1]$.
 In~\cite{DOV18}, the Airy sheet -- and in fact a richer object in which the temporal coordinates of the polymer endpoints are permitted to be general -- is constructed via the Brownian LPP prelimit. An extension of the Robinson-Schensted-Knuth correspondence is used to express the construction in terms of a last passage percolation problem whose underlying environment is a copy of the distributional limit of an ensemble of curves whose uppermost is the narrow wedge profile in Brownian LPP. 

Thus, we may speak rigorously and unambiguously of $y \to \weight_{\infty;(*:f,0)}^{(y,1)}$.
It is natural to think, then, that a form of 
 Theorem~\ref{t.unifpatchcompare} could be asserted in the limiting case. Taking $n=\infty$ in the result, we would expect to assert a comparable theorem uniformly in $f \in  \initcond_{\ovbar\coninit}$; the fourth parameter~$1/12$ would play no role (and formally could be an arbitrary positive value), in the sense that the $n$-dependent non-smallness condition involved in the comparison with Brownian bridge would be absent, comparison being made {\em above scale zero}. All this is very plausible, but a task concerning weak convergence of measures should be carried out in order to pass the quilt description to the limit.
 
Given that the existence of  $y \to \weight_{\infty;(*:f,0)}^{(y,1)}$ follows from~\cite{DOV18}, it is natural to extend Johansson's conjecture, and thus to suggest that this random function has an almost surely unique maximizer for any $f$ lying in one of the function spaces  $\initcond_{\ovbar\coninit}$. Note that the example at the end of Section~\ref{s.patchwork} shows that the passage of the Brownian patchwork quilt description to $n = \infty$
 would not in itself be adequate for proving this extended conjecture. 
 
 \subsubsection{The broader study of Brownian LPP in scaled coordinates}
 This paper has written so that it may be read on its own. However, the paper forms part of a four-paper study of scaled Brownian LPP, alongside~\cite{BrownianReg},~\cite{ModCon} and~\cite{NonIntPoly}, and here we comment on this article's relation to the other papers.
 The first of the companion papers develops a theory of Brownian Gibbs line ensembles that applies to scaled Brownian LPP. 
 This theory is hidden for the reader of the present article since it is not needed directly; it is used during the proofs of the central results in~\cite{ModCon} and~\cite{NonIntPoly}, which will be applied here and recalled presently. The basic idea of the proof in this article, which we will explain in Section~\ref{s.roughguide}, after presenting some notation and tools, is crucial to the four-paper study, since the task of implementing it rigorously has required the development of tools throughout the study. The reader may consult \cite[Section~$1.2$]{BrownianReg} for an overview of the investigation of scaled Brownian LPP at large.  It is worth bearing in mind, however, that the material in \cite[Subsection~$1.2.3$]{BrownianReg}, which presents a conceptual overview, follows closely the upcoming Section~\ref{s.roughguide}.
 This choice of presentation in~\cite{BrownianReg} reflects the conceptual significance for the overall study of  the basic idea of the proof of Theorem~\ref{t.unifpatchcompare}.
 
\subsubsection{Acknowledgments.}
The author thanks Riddhipratim Basu, Ivan Corwin, Shirshendu Ganguly and Jeremy Quastel for valuable conversations at many stages of this project; Judit Z{\'a}dor for the simulations shown in Figure~\ref{f.eden}; and two referees for their insightful and comprehensive commentary.

 \section{A geometric view: staircases, zigzags and polymers}\label{s.geom}

We intend to explain more of the geometric meaning of Theorem~\ref{t.unifpatchcompare}, and give a rough idea of how this result will be proved, 
before we embark on the proof itself. Indeed, Section~\ref{s.roughguide} will offer such a rough guide to the principal concerned concepts. Before we can offer this outline, we need to develop some preliminaries. In the present section,
we revisit the basic setup of Brownian LPP, and discuss our use of scaled coordinates in a more geometric light.

\subsubsection{Staircases.}\label{s.staircases}
Recall from Section~\ref{s.lpp} that energy is ascribed to any given non-decreasing list $\llbracket i+1,j \rrbracket \to [x,y]$, and that $M_{(x,i) \to (y,j)}^1$ denotes  the supremum of the energy of all such lists.
In order to make a study of those lists that attain this maximum energy, 
we begin by noting that the lists are in bijection with certain subsets of $[x,y] \times [i,j] \subset \R^2$ that we call {\em staircases}.
Staircases offer a geometric perspective on Brownian LPP and perhaps help in visualizing the problems in question.


 The staircase~\hfff{staircase} associated to the non-decreasing list $\big\{ z_k: k \in \llbracket i+1,j \rrbracket \big\}$ is specified as the union of certain horizontal planar line segments, and certain vertical ones.
The horizontal segments take the form $[ z_k,z_{k+1} ] \times \{ k \}$ for $k \in \llbracket i , j \rrbracket$.
Here, the convention that $z_i = x$ and  $z_{j+1} = y$ is again adopted. 
The right and left endpoints of each consecutive pair of horizontal segments are interpolated by a vertical planar line segment of unit length. It is this collection of vertical line segments that form
the vertical segments of the staircase.

The resulting staircase may be depicted as the range of an alternately rightward and upward moving path from starting point $(x,i)$ to ending point $(y,j)$. 
The set of staircases with these starting and ending points will be denoted by $\staircase_{(x,i) \to (y,j)}$.
Such staircases are in bijection with the collection of non-decreasing lists considered earlier. Thus, any staircase $\phi \in \staircase_{(x,i) \to (y,j)}$
is assigned an energy~\hfff{energy} $E(\phi) = \sum_{k=i}^j \big( B ( k,z_{k+1} ) - B( k,z_k ) \big)$ via the associated $z$-list.

\subsubsection{Energy maximizing staircases are called geodesics.}
A staircase  $\phi \in \staircase_{(x,i) \to (y,j)}$ whose energy  attains the maximum value $M^1_{(x,i) \to (y,j)}$ is called a geodesic~\hfff{geodesic} from $(x,i)$ to~$(y,j)$.
It is a simple consequence of the continuity of the constituent Brownian paths $B(k,\cdot)$
that such a geodesic exists for all choices of $(x,y) \in \R^2$ with $x \leq y$.
 The geodesic with given endpoints is known to be almost surely unique: take $\ell = 1$ in Lemma~\ref{l.severalpolyunique}, which appears in Appendix~\ref{s.polyunique}.

\subsubsection{The scaling map.}
For $n \in \N$, consider the $n$-indexed {\em scaling} map\hfff{scalingmap} $R_n:\R^2 \to \R^2$ given by
$$
 R_n \big(v_1,v_2 \big) = \Big( 2^{-1} n^{-2/3}( v_1 - v_2) \, , \,   v_2/n \Big) \, .
$$ 
 
The scaling map acts on subsets $C$ of $\R^2$ by setting
$R_n(C) = \big\{ R_n(x): x \in C \big\}$.

\subsubsection{Scaling transforms staircases to zigzags.}
The image of any staircase under $R_n$
will be called an $n$-zigzag.~\hfff{zigzag} The starting and ending points of an $n$-zigzag $Z$ are defined to be the image under $R_n$
of such points for the staircase $S$ for which $Z = R_n(S)$.
 
Note that the set of horizontal lines is invariant under $R_n$, while vertical lines are mapped to lines of gradient  $- 2 n^{-1/3}$.
As such, an $n$-zigzag is the range of a piecewise affine path from the starting point to the ending point which alternately moves rightwards  along horizontal line segments  and northwesterly along sloping line segments, where each sloping line segment has gradient  $- 2 n^{-1/3}$; the first and last segment in this journey may be either horizontal or sloping.

 \subsubsection{Scaled geodesics are called polymers.}
 For $n \in \N$, the image of any geodesic under the scaling map $R_n$ will be called an $n$-polymer, or often simply a polymer.\hfff{polymer} This usage of the term `polymer' for `scaled geodesic' is apt for our study, due to the central role played by scaled geodesics. The usage is not, however, standard: the term `polymer' is often used to refer to typical realizations on the path measure in LPP models at positive temperature.

\subsubsection{Zigzags have weight}
Any $n$-zigzag $Z$ from $(x,i/n)$ to $(y,j/n)$  is ascribed a scaled weight~\hfff{weight} 
$\weight(Z) = \weight_n(Z)$ given by 
$$
 \weight(Z) =  2^{-1/2} n^{-1/3} \Big( E(S) - 2(j - i)  - 2n^{2/3}(y-x) \Big) 
$$
where $Z$ is the image under $R_n$ of the staircase $S$. Thus, a polymer maximizes weight among the zigzags that share its endpoints, just as a geodesic maximizes energy over staircases.

\subsubsection{Some basic notation}
For $k \geq 1$, we write $\R^k_\leq$
for the subset of $\R^k$ whose elements $(z_1,\cdots,z_k)$
are non-decreasing sequences. When the sequences are increasing, we instead write $\R^k_<$. We also use the notation $A^k_\leq$ and $A^k_<$.
Here, $A \subset \R$ and the sequence elements are supposed to belong to $A$.
We will typically use this notation when $k=2$.
 
 \subsubsection{Compatible triples.}
 Let $(n,t_1,t_2) \in \N \times \R^2_<$, which is to say that $n \in \N$ and $t_1,t_2 \in \R$ with $t_1 < t_2$.
 Taking $x,y \in \R$, does there exist an $n$-zigzag from $(x,t_1)$ to $(y,t_2)$?
 As far as the data $(n,t_1,t_2)$ is concerned, such an $n$-zigzag may exist only if 
 \begin{equation}\label{e.ctprop}
     \textrm{$t_1$ and $t_2$ are integer multiplies of $n^{-1}$} \, .
\end{equation}
We say that data $(n,t_1,t_2)  \in \N \times \R^2_<$ is a {\em compatible triple}~\hfff{comptriple} if it verifies the last condition.

An important piece of notation associated to a compatible triple is $\tot$,~\hfff{tot} which we will use to denote the difference $t_2 - t_1$. The law of the underlying Brownian ensemble $B: \Z \times \R \to \R$ is invariant under integer shifts in the first, curve indexing, coordinate. This translates to a distributional invariance of scaled objects under vertical shifts by multiples of $n^{-1}$, something that makes the parameter $\tot$
of far greater significance than $t_1$ or $t_2$.

\begin{figure}[ht]
\begin{center}
\includegraphics[height=7cm]{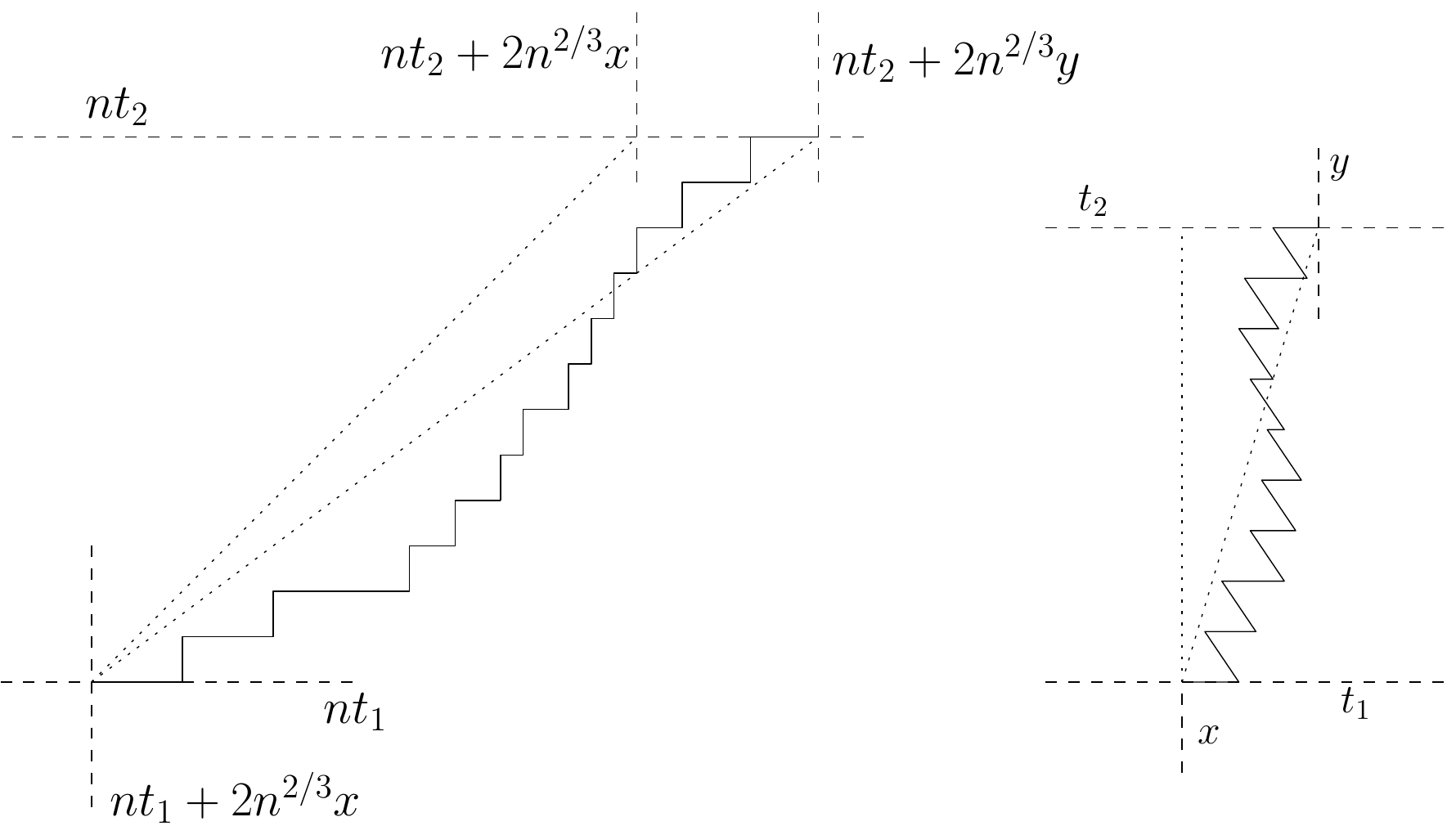}
\caption{Let $(n,t_1,t_2)$ be a compatible triple and let $x, y \in \R$. The endpoints of the geodesic in the left sketch have been selected so that, when the scaling map~$R_n$ is applied to produce the right sketch,
the $n$-polymer $\rho_{n;(x,t_1)}^{(y,t_2)}$ results.}
\label{f.scaling}
\end{center}
\end{figure}

Supposing now that $(n,t_1,t_2)$ is indeed a compatible triple, 
the condition 
 $y \geq x - 2^{-1} n^{1/3} \tot$ ensures that the preimage of $(y,t_2)$ under the scaling map $R_n$ lies northeasterly of the preimage of $(x,t_1)$.
 Thus, an  $n$-zigzag from $(x,t_1)$ to $(y,t_2)$ exists in this circumstance. We have mentioned that geodesics exist uniquely for given endpoints almost surely.
 Taking the image under scaling, this translates to the almost sure existence and uniqueness of the $n$-polymer from  $(x,t_1)$ to $(y,t_2)$. This polymer will be denoted $\rho_{n;(x,t_1)}^{(y,t_2)}$:~\hfff{polynot} see Figure~\ref{f.scaling}. 
 In notation that encompasses the special case of $t_1 =0$ and $t_2=1$ seen in~(\ref{e.weightmzeroone}), the weight of this polymer will be recorded as~$\weight_{n;(x,t_1)}^{(y,t_2)}$;~\hfff{maxweight} (in fact, the latter weight can be ascribed meaning as the weight maximum for the given endpoints even in the probability zero event that the polymer in question is not well defined).
 These two pieces of notation share something in common with a number of later examples, in which objects are described in scaled coordinates.
Round bracketed expressions in the subscript or superscript will refer to  
a space-time pair, with the more advanced time in the superscript.
Typically some aspect of the $n$-polymer from $(x,t_1)$ to $(y,t_2)$ is being described when this ${\mathsf{TBA}}_{n;(x,t_1)}^{(y,t_2)}$ notation is used.  
A glossary of notation in 
Appendix~\ref{s.glossary} includes several examples that employ or develop this usage.

 \section{Polymer forests: a rough guide to the proof of the main result}\label{s.roughguide}
 
Some of the central players in our study have been introduced: the use of scaled coordinates, polymers with given endpoints such as $\rho_{n;(x,t_1)}^{(y,t_2)}$, and polymer weights such as $\weight_{n;(x,t_1)}^{(y,t_2)}$.

Our principal conclusion, Theorem~\ref{t.unifpatchcompare}, asserts that the polymer weight profile $y \to \weight_{n;(*:f,0)}^{(y,1)}$ is a Brownian patchwork quilt, with the quality of the sewing (and the resulting comparison to Brownian motion) not deteriorating even as the initial condition $f$ is permitted to vary over the very broad function class $\initcond_{\ovbar\coninit}$, and $n$ over all sufficiently high integers. 
 
A valuable point of departure for beginning to understand how we will obtain this theorem is to consider the special case where $f:\R \to \R \cup \{-\infty\}$ equals zero at the origin and is otherwise minus infinity. This is the narrow wedge initial condition, in which all competing polymers must begin at zero. A crucial first observation for our proof is that, in this special case, the resulting weight profile is known 
enjoy a very strong comparison to Brownian motion.
Indeed, a result in~\cite{BrownianReg} that we will later recall as Theorem~\ref{t.bridge} implies that, for this $f$,   
 the random function $y \to \weight_{n;(*:f,0)}^{(y,1)}$, 
 defined on any given compact interval, withstands $L^{\infty-}$-comparison to Brownian bridge. 
 Weight profiles arising from such narrow wedge initial conditions, with fixed starting point, would thus seem to represent excellent candidates for the role of fabric sequence elements in the quilt description that we seek to offer of a much more general polymer weight profile. 
 
\begin{figure}[ht]
\begin{center}
\includegraphics[height=12cm]{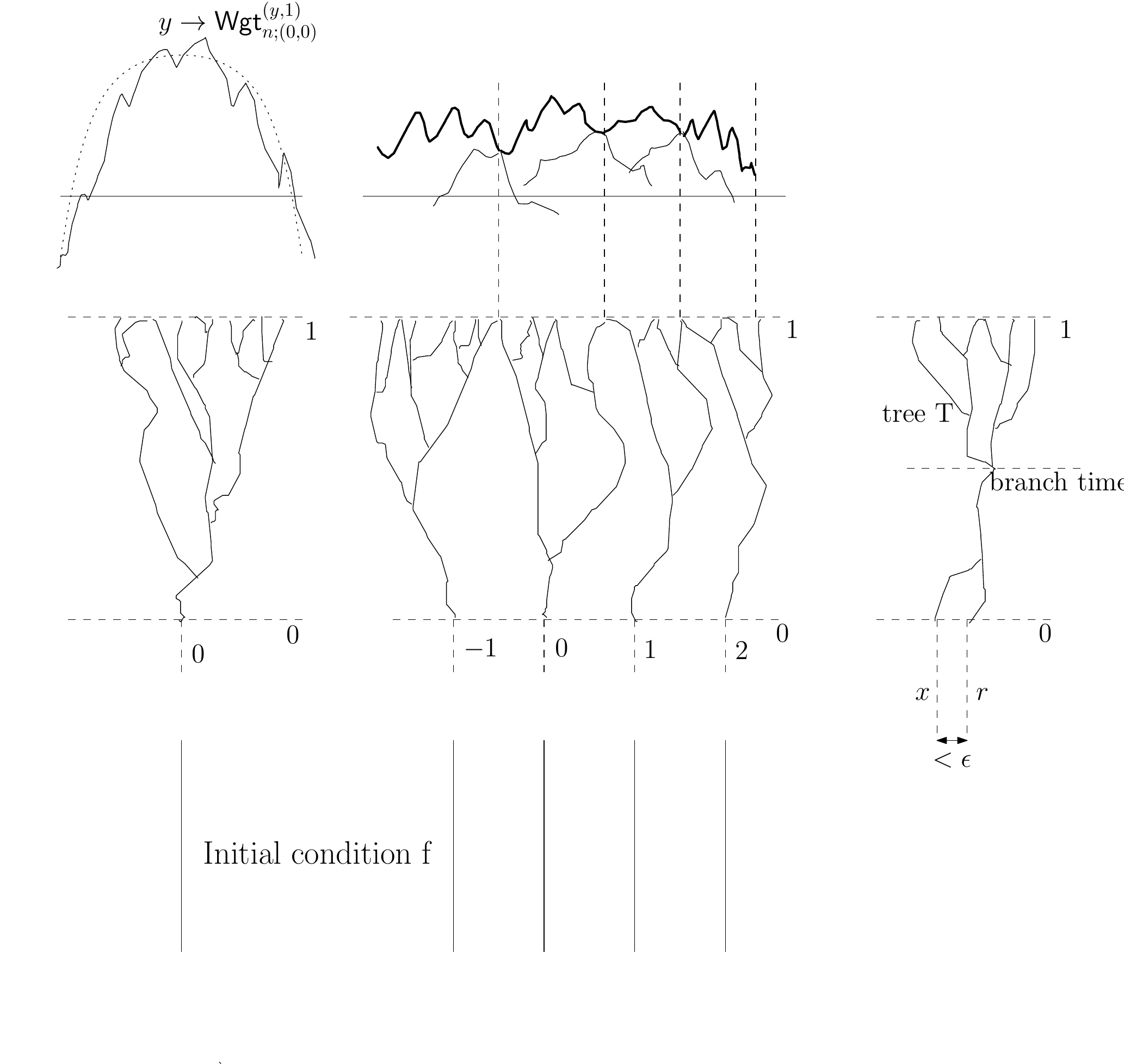}
\caption{In the left and middle sketches, the weight profile $y \to \weight_{n;(*:f,0)}^{(y,1)}$, the $f$-rewarded polymer forest, and the function $f$ are depicted.
 On the left, $f$ is zero at zero and otherwise minus infinity, so that the profile is the narrow wedge case $y \to \weight_{n;(0,0)}^{(y,1)}$.
In the middle, $f$ is instead zero at integer points. The  profile $y \to \weight_{n;(*:f,0)}^{(y,1)}$ is depicted in bold in the upper middle sketch; it is the maximum of the also indicated profiles $y \to \weight_{n;(k,0)}^{(y,1)}$, where $k \in \Z$ indexes the canopy roots in this case. The dashed vertical lines that contact the line at height one indicate canopy boundary points.
The right sketch depicts one tree, rooted at $r$, in an $f$-rewarded polymer forest. The proposed rerooting of this tree to a nearby discrete mesh element $x \in \e \Z$ is illustrated.}
\label{f.manytrees}
\end{center}
\end{figure}

 In the left sketch of Figure~\ref{f.manytrees} is illustrated a geometric view of the narrow wedge weight profile. Note that $\weight_{n;(0,0)}^{(y,1)}$
 is the weight of the polymer $\rho_{n;(0,0)}^{(y,1)}$. These polymers, indexed by $y \in \R$, all stream out of the origin at time zero, to arrive at their various ending locations $y$ at time one. The almost sure uniqueness of polymers with given endpoints suggests that, once separated, polymers will not meet again. Thus, the system of polymers should be viewed as a tree, with a root at $(0,0)$, and a canopy $\R \times \{ 1 \}$ of ending locations.

 Return now to a general initial condition $f: \R \to \R \cup \{ - \infty \}$.
 The $f$-rewarded line-to-point polymers may be traced backwards in time from locations $(y,1)$ with $y \in \R$. They arrive at time zero at a variety of locations, in contrast to the narrow wedge case. 
 They share something with that case, however: as time decreases from one to zero, two polymers that meet will stay together; this may again be expected on the basis of the uniqueness of polymers with given endpoints. 
 
 This fact has the implication that we may view the collection of polymers  -- indexed by $y \in \R$, and labelled $\rho_{n;(*:f,0)}^{(y,1)}$ in a natural extension of notation -- as a forest, which we may call the $f$-rewarded polymer forest.
 Each constituent tree has a root lying on the $x$-axis, and a canopy that consists of an interval lying in the line at height one. 
 Indeed, we may partition this copy of the real line into this set of canopies. 
 The polymer weight profile $y \to \weight_{n;(*:f,0)}^{(y,1)}$, when restricted to any given canopy, 
 would seem to have much in common with the narrow wedge weight profile. After all, all the concerned polymers, ending at locations $(y,1)$, for points $y$ in the given canopy, share their starting location, namely $(r,0)$, where $r \in \R$ is the root of the canopy in question: see the middle sketch of  Figure~\ref{f.manytrees}.

 This then is the basis of our Brownian patchwork quilt description of the general weight profile 
 $y \to \weight_{n;(*:f,0)}^{(y,1)}$. The patches of the quilt should be the canopies in the $f$-rewarded polymer forest, so that the stitch point set consists of boundary points between consecutive canopies. The fabric sequence elements take the form  $y \to \weight_{n;(r,0)}^{(y,1)}$ where $r$ is the root associated to a tree in the polymer forest. Since $r$ is fixed, the fabric functions should withstand a strong comparison to Brownian bridge, since the narrow wedge weight profiles do. 
 
 This plan is a useful guide to the approach by which we will prove Theorem~\ref{t.unifpatchcompare}. There are however significant challenges to be overcome in implementing it. The structure of the problems is elucidated by voicing two objections to the plan:
  
 \begin{itemize}
 \item Under the plan, the stitch point set is said to be the collection of canopy boundary points. In order to have a useful description, we want this set to be not too big. Certainly, we want to argue that, for all high $n$, the stitch point set on a unit interval has cardinality that is tight. This raises a geometric question: can we show that the forest of $f$-rewarded polymers has a tight number of trees per unit length?\footnote{If such an assertion is proved, it might seem desirable to demonstrate that every canopy is a non-trivial interval, rather than a singleton set. This fact is not needed for the proposed analysis, however, because a singleton set canopy would contribute a vacuous element to the fabric sequence, and thus would play no role in specifying a patchwork quilt.} 
 \item The plan proposes that the $f$-rewarded weight profile on a given canopy be viewed as   $y \to \weight_{n;(r,0)}^{(y,1)}$ where $r$ is the root in question. We may hope that because $r$ is fixed, this object is statistically similar to the narrow wedge case, where say $r$ is zero. There is a significant objection, however. The location $r$ is random, and it may be exceptional. As such, properties of the narrow profile    $y \to \weight_{n;(x,0)}^{(y,1)}$, such as locally Brownian structure, that obtain for typical values of $x \in \R$, may disappear for the possibly exceptional choices of root $r$ that we are considering.
 \end{itemize}
 
 We end this rough guide by giving the briefest indication of how we will resolve these problems. 
 
 Regarding the first, consider again the $f$-rewarded polymer forest. Pick a location in the interior of each canopy, so that a point process on $\R \times \{ 1 \}$ results.
 From each point, draw the $f$-rewarded line-to-point polymer back to time zero. This system of polymers is pairwise disjoint, because the time-one endpoints lie in different canopies.
 If the tightness alluded to in the first problem were to fail, we would be drawing, for high $n$, many disjoint polymers, all ending on a given unit-order interval. 
This behaviour has been proved to be a rarity for Brownian LPP. As we will recall in Theorem~\ref{t.maxpoly}, it has been proved in~\cite[Theorem~$1.4$]{NonIntPoly}  
that the probability that there exist $k$ disjoint polymers, beginning and ending in a unit interval, has a superpolynomial tail in $k$, uniformly in high choices of the scaling parameter $n \in \N$.

Regarding the second difficulty, that root locations $r$ may be exceptional, we may try to solve the problem by a rerooting procedure: see the right sketch in Figure~\ref{f.manytrees}. Close to the location $r$ is an element $x$ in the discrete $\e$-mesh $\e \Z$, where $\e > 0$ is small but fixed. The polymer weight profile that we are trying to describe,  $y \to \weight_{n;(r,0)}^{(y,1)}$, may be  compared to the narrow wedge profile  $y \to \weight_{n;(x,0)}^{(y,1)}$. Because $x$ lies in a discrete mesh, it may be viewed as a typical location, so that the objection raised in the second bullet point does not arise.
But why should the nearby profile  $y \to \weight_{n;(x,0)}^{(y,1)}$ be expected to offer an accurate description of  $y \to \weight_{n;(r,0)}^{(y,1)}$? 
We will now present an heuristic argument that, since $r$ and $x$ are close, these two functions should, with a high probability determined by this closeness, differ by a random constant as $y$ varies over the canopy of the tree of which $r$ is the root. The tree in question has a {\em branch} time, with the constituent polymers following a shared course (the tree trunk, if you like) until that time, after which they may go their separate ways to the various locations in the canopy. 
Now, when the tree is rerooted the short distance from $(r,0)$ to $(x,0)$, we may expect that its structure changes by a modification in the form of the polymers only close to time zero: see again the right sketch of   Figure~\ref{f.manytrees}.  Provided this modification has finished by the branch time, it will affect only the form of the tree trunk, and will be shared by all concerned polymers, no matter at  which canopy point $(y,1)$ they end. 
For this reason, the difference $y \to \weight_{n;(x,0)}^{(y,1)} - \weight_{n;(r,0)}^{(y,1)}$ may be expected to equal the discrepancy in weight between the new and the old tree trunks with high probability,
and thus typically be independent of the canopy location~$y$. 
  

The actual solution we undertake will modify these suggestions. In particular, the resolution of the second difficulty will proceed by a means that does involve rerooting, but not quite in the way we have just suggested. In any case, this rough guide may help to set the reader's bearings, as we now turn to set out precisely some of the tools that will be needed when it is implemented rigorously.
 
At roughly this moment in a paper, it is conventional to offer an overview of the structure of the remainder. The story of the proof will unfold in several steps: for example, after the next section of tools, we will return to the rough guide and explain more carefully how rerooting will resolve the second objection. It is probably not helpful to attempt to indicate the global structure of the article for now, beyond pointing that  the actual proof of Theorem~\ref{t.unifpatchcompare} will appear in the final Section~\ref{s.quilt}, and that the road to that section will  alternate between heuristic elaboration of the rough guide and rigorous development of the machinery necessary for the rough guide's implementation.
We also mention that, beyond the glossary of notation in Appendix~\ref{s.glossary}, the paper has two further
 appendices, whose respective roles will be described in Section~\ref{s.polymerbasic} and Subsection~\ref{s.rolehyp}. 
 
\section{Important tools}
 
 If the rough guide is to be implemented, several significant inputs will be needed. It is the task of this section to quote the necessary results.

\subsection{A comment about explicit constants}\label{s.explicit}

The results that we are quoting state hypotheses on parameters in fairly explicit terms, and this means that we must define some explicit constants before beginning to state the results.
Certain results have been proved using a framework of Brownian Gibbs line ensembles, which are in essence systems of mutually avoiding Brownian bridges that arise in Brownian LPP via the Robinson-Schensted-Knuth correspondence. Several of the results we will be quoting are proved using a regularity property of such ensembles, with this property specified in terms of  two positive parameters $(c,C)$. We mention these things now merely in order to explain that these two parameters $(c,C)$ are fixed throughout the paper.  It is~\cite[Proposition~$4.2$]{NonIntPoly}  
 that determines their values.  
 The reader who wishes to understand more about the Brownian Gibbs framework may consult~\cite{BrownianReg}, in particular, this article between the start of Section~$1.3$ and the end of Chapter~$2$.

The results to be quoted will make reference to two fixed sequences of constants that are expressed in terms of the given positive values  $C$ and $c$.
These  two sequences are  $\big\{ C_k : k \geq 2 \big\}$ and $\big\{ c_k: k \in \N \big\}$. They are specified by setting, for each $k \geq 2$,   
$$
 \formerE_k = \max \Big\{  10 \cdot 20^{k-1} 5^{k/2} \Big( \tfrac{10}{3 - 2^{3/2}} \Big)^{k(k-1)/2} C \, , \, e^{c/2} \Big\} 
$$
as well as $C_1 = 140 C$; and
$$
 c_k =   \big( (3 - 2^{3/2})^{3/2} 2^{-1} 5^{-3/2} \big)^{k-1} c_1 \, ,
$$ 
with  $c_1 = 2^{-5/2} c \wedge 1/8$. 

\subsection{The rarity of many disjoint polymers}

In the rough guide, we saw that this rarity will be needed to resolve the first objection. In fact, we will use two related results, Theorems~\ref{t.maxpoly} and~\ref{t.disjtpoly}.

Let $(n,t_1,t_2) \in \N \times \R^2_<$ be a compatible triple, and let $I,J \subset \R$ be intervals. Set $\maxpoly_{n;(I,t_1)}^{(J,t_2)}$~\hfff{maxcard}\\
 equal to  the maximum cardinality of a disjoint set of   $n$-polymers each of whose start and end points have the respective forms $(x,t_1)$ and $(y,t_2)$ 
where $x$ is some element of $I$ and $y$ is some element of~$J$.


(This $\maxpoly$ notation fits within the framework discussed in the last paragraph of Section~\ref{s.geom}. In this case, the first element in the space-time pairs $(I,t_1)$ and $(J,t_2)$ is an interval, rather than a point.)

The next result is~\cite[Theorem~$6.2$]{NonIntPoly}:
 the maximum disjoint polymer cardinality has a tail that decays super-polynomially.
\begin{theorem}\label{t.maxpoly}
There exist constants $K_0 \geq 1$, $a_0 \in (0,1)$,  $k_0 > 0$ and  a sequence of positive constants $\big\{ \conseqmac_i: i \in \N \big\}$ for which $\sup_{i \in \N} \conseqmac_i \exp \big\{ - 2 (\log i)^{11/12} \big\}$ is finite, such that the following holds.
Let
$(n,t_1,t_2) \in \N \times \R^2_<$ be a compatible triple. Let $x,y \in \R$, $a,b \in \N$ and $k \in \N$. Write $h = a \vee b$. Suppose that 
 $$
 k \geq   k_0 \vee \big( \vert x - y \vert \tot^{-2/3} + 2h \big)^3  
 $$ 
 and
\begin{eqnarray}
 n \tot  & \geq  & \max \bigg\{ \,    2(K_0)^{(12)^{-2} (\log \log k)^2} \big( \log k \big)^{K_0}
    \, , \,  a_0^{-9} \big( \vert y - x \vert  \tot^{-2/3} + 2h \big)^9 \, ,
    \label{e.nmaxpoly} \\
  & & \qquad \qquad    
     10^{325}       \rsc^{-36} k^{465}  \max \big\{   1  \, , \,   \big(\vert y - x \vert \tot^{-2/3} + 2h \big)^{36}   \big\}  
   \,  \bigg\} \, . \nonumber
\end{eqnarray}

Then
$$
 \PP \bigg(  \maxpoly_{n;\big([x,x+a\tot^{2/3}],t_1\big)}^{\big([y,y+b \tot^{2/3}],t_2\big)} \geq k \bigg)  \, \leq \,  k^{-   (145)^{-1} ( \log \beta)^{-2} (0 \vee \log \log k)^2}  \cdot h^{(\log \beta)^{-2} (\log \log k)^2/{288}  + 3/2} \conseqmac_k  \, .
$$
Here, $\beta$ is specified to be $e \vee \limsup_{i \in \N} \beta_i^{1/i}$, where the sequence of constants $\big\{ \beta_i: i \in \N \big\}$, which verifies $\limsup \beta_i^{1/i} < \infty$, is supplied by~\cite[Corollary~$5.2$]{NonIntPoly}.
\end{theorem}
We make an expository comment about result statements which is illustrated by the last theorem: we have chosen to be fairly explicit in stating hypothesised parameter bounds.
One effect of this is to make some statements look complicated, and we now suggest to the reader some ways of reading the results that may serve to focus attention on their basic meaning. 
First, a basic scaling principle, explained in \cite[Section~$2.3$]{ModCon}  means that there is no loss of generality in taking $t_1 = 0$ and $t_2=1$, and thus $\tot=1$, in results stated using this parameter pair.
Choosing $x$, $y$ bounded above in absolute value, by one say, may also be useful. Also taking $a=b=1$ in the present case, we see that the requirements in the theorem that $k$ is bounded below by a universal constant, and that $n$ is at least a polynomial in $k$, are quite weak, since our interest is in statements that hold uniformly in high $n$. In this light, the basic conclusion of the last theorem is seen to be that the probability that $k$ disjoint polymers cross a unit interval in unit time is at most $k^{- \Theta(1) (\log \log k)^2}$, uniformly in high $n$.

The next result,~\cite[Theorem~$6.1$]{NonIntPoly}, 
 asserts a more local form of disjoint polymer rarity:
the probability that a given number $m \in \N$
of disjoint polymers begin and end in given intervals of a small length $\e > 0$ 
is at most $\e^{m^2 - 1 +o(1)}$.
It shares the parameters $K_0$ and $a_0$,  and  the sequence  $\big\{ \beta_i: i \in \N \big\}$, with the preceding result. It uses a further {\em positive} constant $\eta_0$.
\begin{theorem}\label{t.disjtpoly}
Let
$(n,t_1,t_2) \in \N \times \R^2_<$ be a compatible triple.
Let $m \in \N$,  $\e > 0$ and $x,y \in \R$ satisfy the conditions that $m \geq 2$,
\begin{equation}\label{e.epsilonbound}
 \e \leq \min \Big\{ \, (\eta_0)^{4m^2} \, , \, 10^{-616}  c_m^{22}  m^{-115} \,  , \, \exp \big\{ - C^{3/8} \big\} \, \Big\} \, ,
\end{equation}
\begin{equation}\label{e.nlowerbound}
    n \tot \geq \max \bigg\{ \, 2(K_0)^{m^2} \big( \log \e^{-1} \big)^{K_0}
    \, , \, 10^{606}   c_m^{-48}  m^{240}  \rsc^{-36} \e^{-222}  \max \big\{   1  \, , \,   \vert x - y \vert^{36} \tot^{-24}  \big\} \, , \, 
  a_0^{-9} \vert y - x \vert^9 \tot^{-6}  \,  \Bigg\} \, ,
\end{equation}
as well as $\vert y - x \vert \tot^{-2/3} \leq \e^{-1/2} \big( \log \e^{-1} \big)^{-2/3} \cdot 10^{-8} c_m^{2/3} m^{-10/3}$.
Then
\begin{eqnarray*}
 & & \PP \bigg( \maxpoly_{n;\big([x-\tot^{2/3}\e,x+\tot^{2/3}\e],t_1\big)}^{\big([y-\tot^{2/3}\e,y+\tot^{2/3}\e],t_2\big)}  \geq m \bigg) \\
 & \leq & 
  \e^{(m^2 - 1)/2}  \cdot   10^{32m^2}  m^{15m^2}  c_m^{-3m^2} C_m   \big( \log \e^{-1} \big)^{4m^2}   \exp \big\{ \beta_m \big( \log \e^{-1} \big)^{5/6} \big\} \, .
\end{eqnarray*}
\end{theorem}
For this result, it is helpful to bear in mind that $m$ is fixed independently of $n$, so that the probability that $m$ polymers disjointly cross between two similarly placed intervals of a short width~$\e$ in unit time is found to be at most $\e^{m^2 - 1 + o(1)}$, uniformly in high $n$.

\subsection{Narrow-wedge polymer weight profiles bear strong comparison to Brownian bridge}

This is a fundamental idea in the rough guide and we now recall the result in question, which is in essence~\cite[Theorem~$2.11$]{BrownianReg}.
For the statement,
recall from Subsection~\ref{s.bridgeproject} our notation for the space of continuous functions, and that of bridges, defined on a compact real interval, and for the bridge projection between them, as well as the standard Brownian bridge law.

\begin{theorem}\label{t.bridge}
 Let $K \in \R$ and $\ipdval \geq 1$. 

Let $x \in \R$. For $n \in \N$ satisfying $K \geq x - 2^{-1} n^{1/3}$, we define a $\mc{C}_{*,*} \big( [K,K + \ipdval]  , \R \big)$-valued stochastic process $\mc{L} = \mc{L}_{n,x}:[K,K + \ipdval] \to \R$
by setting  $\mc{L}(y) = \weight_{n;(x,0)}^{(y,1)}$ for each $y \in [K,K + \ipdval]$.

 Endowing the bridge space $\mc{C}_{0,0} \big( [K,K + \ipdval]  , \R \big)$ with the topology of uniform convergence, let  
 $A$ be any measurable subset of this space. Set $a = \mc{B}_{0,0}^{[K,K + \ipdval]}(A)$. Then there exists a  constant $G > 0$ such that, provided that $n \geq G$ and 
$\min \big\{ G^{-1} , \exp \{ - d^3 \} \big\} \geq a \geq \exp \big\{   - G^{-1} n^{1/12} \big\}$, 
we have that 
  $$
\PP \big( \mc{L}^{[K,K+\ipdval]}  \in A \big) \leq a \cdot \ipdval^{1/2}  G  \exp \big\{ G \ipdval^2 (\log a^{-1})^{5/6}  \big\} \, .
 $$ 
 The constant $G > 0$ may be selected independently of  the choice of $(n,K,d,x)$ subject to $\vert K - x \vert + \ipdval \leq 2^{-1} c (n+1)^{1/9}$. 
 \end{theorem}
 \noindent{\bf Proof.} 
 The result~\cite[Theorem~$2.11$]{BrownianReg}
 depends on the apparatus of the lengthy paper~\cite{BrownianReg} in a substantial way.
 The reader who wants to understand the proof should consult that work. Here, we merely explain formally how to derive the theorem from this result.
 We mentioned in Section~\ref{s.explicit} that certain regular Brownian Gibbs line ensembles play an important role, and this is certainly true of this theorem.
 Indeed,~\cite[Proposition~$4.2$]{NonIntPoly} 
  implies that the stochastic process $y \to \mc{L}(x+y)$ is the lowest indexed curve in a line ensemble with $n+1$ curves
 that, in the language of~\cite{BrownianReg}, is a $(\bar\phimac,c,C)$-regular Brownian Gibbs ensemble, where $\bar\phimac$ is the vector $(1/3,1/9,1/3)$ and the positive constants $c$ and $C$ are furnished by~\cite[Proposition~$4.2$]{NonIntPoly}. (We mention that the basic result that underlies~\cite[Proposition~$4.2$]{NonIntPoly} is~\cite[Proposition~$2.5$]{BrownianReg}.) 
  For this reason,~\cite[Theorem~$2.11$]{BrownianReg} 
  implies the result; note that the index $n$ there assumes the value $n+1$ given by the present context.
  The uniform choice of $G$ subject to $\vert K - x \vert + \ipdval \leq 2^{-1} c (n+1)^{1/9}$
  is a consequence of the hypothesis $[K+x,K+x+d] \subset c/2 \cdot [-n^{1/9},n^{1/9}]$ that is needed to apply~\cite[Theorem~$2.11$]{BrownianReg} and the explicit form of the relevant constants when this result is applicable. \qed

\subsection{Basic results on polymers and weights}\label{s.polymerbasic}

Some such basic results  will be needed, for example, to make sense of polymer forests. Here we state the relevant tools.

Let $(n,t_1,t_2)$ be a compatible triple.
We introduce an ordering relation $\preceq$ on the space of $n$-polymers with lifetime $[t_1,t_2]$. To define it, let $(x_1,x_2),(y_1,y_2) \in \R^2$ and  consider
a polymer  $\rho_1$ from $(x_1,t_1)$ to $(y_1,t_2)$ and another
 $\rho_2$ from $(x_2,t_1)$ to $(y_2,t_2)$. We declare that $\rho_1 \preceq \rho_2$ if `$\rho_2$ lies on or to the right of $\rho_1$': formally, if $\rho_2$ is contained in the union of the closed horizontal planar line segments whose left endpoints lie in $\rho_1$.

First is a simple sandwiching result,~\cite[Lemma~$5.7$]{NonIntPoly}.
\begin{lemma}\label{l.sandwich}
Let $(n,t_1,t_2)$ be a compatible triple, and let  $(x_1,x_2),(y_1,y_2) \in \R_\leq^2$. Suppose that there is a unique $n$-polymer from $(x_i,t_1)$ to $(y_i,t_2)$, both when $i=1$ and $i=2$. (The shortly upcoming Lemma~\ref{l.densepolyunique}(1) asserts that this circumstance occurs almost surely; the resulting polymers have been labelled $\rho_{n;(x_1,t_1)}^{(y_1,t_2)}$ and  $\rho_{n;(x_2,t_1)}^{(y_2,t_2)}$.)
Now let $\rho$ denote any $n$-polymer that begins in $[x_1,x_2] \times \{ t_1\}$
and ends in  $[y_1,y_2] \times \{ t_2 \}$. Then $\rho_{n;(x_1,t_1)}^{(y_1,t_2)} \preceq \rho \preceq \rho_{n;(x_2,t_1)}^{(y_2,t_2)}$.
\end{lemma}

Next we discuss $f$-rewarded polymers.

\begin{definition}\label{d.polyunique}
\begin{enumerate}
\item
Let $\polyunique_{n;0}^1$~\hfff{polyunique} denote the set of pairs $(x,y) \in \R^2$
for which there exists an $n$-zigzag from $(x,0)$ to $(y,1)$ whose weight uniquely attains the value $\mathsf{W}_{n;(x,0)}^{(y,1)}$.
In other words, $(x,y) \in \polyunique_{n;0}^1$ precisely when the polymer $\rho_{n;(x,0)}^{(y,1)}$ is well defined.
\item Let $f: \R \to \R \cup \{ - \infty \}$ be measurable. 
Let $\polyunique_{n;(*:f,0)}^1$ denote the set of $y \in \R$
such that there exists a unique choice of $(x,\phi)$, with $x \in \R$ and $\phi$ an $n$-zigzag from $(x,0)$ to $(y,1)$, such that 
the weight $\weight (\phi)$ is equal to
  $\mathsf{W}_{n;(*:f,0)}^{(y,1)}$. That is, $y \in \R$ is an element of $\polyunique_{n;(*:f,0)}^1$~\hfff{fpolyunique} precisely when the $f$-rewarded line-to-point polymer that ends at  $(y,1)$
  is well defined. It is this polymer that we call $\rho_{n,(*:f,0)}^{(y,1)}$.~\hfff{fpolymer}
\end{enumerate}  
\end{definition}

For $\ovbar\coninit \in (0,\infty)^3$,
recall the function space $\initcond_{\ovbar\coninit}$ from Definition~\ref{d.if}. 
The next result will be proved in Appendix~\ref{s.polyunique}.

\begin{lemma}\label{l.densepolyunique}
 \begin{enumerate}
 \item Let $x,y \in \R$. 
 There exists an $n$-zigzag from $(x,0)$ to $(y,1)$
 if and only if $y \geq x - n^{1/3}/2$.
 When the last condition is satisfied, there is almost surely a unique $n$-polymer from $(x,0)$ to $(y,1)$, which is to say, $(x,y) \in \polyunique_{n;0}^1$ almost surely. 
 \item Suppose that $n \in \N$ satisfies
$n > 2^{-3/2} \coninit_1^3 \vee 8 (\coninit_2  + 1)^3$, and that $f \in \initcond_{\ovbar\coninit}$ for some $\ovbar\coninit \in (0,\infty)^3$.
Then $[-1,1] \setminus \polyunique_{n;(*:f,0)}^1$ is almost surely of Lebesgue measure zero (and contains any given value $y \in [-1,1]$ with zero probability).
\end{enumerate}
\end{lemma} 
We mention that the point-to-point polymer uniqueness Lemma~\ref{l.densepolyunique}(1) has been expressed in scaled coordinates. Indeed, we have been eager to move promptly to the use of scaled coordinates throughout this article, since our whole perspective makes them essential. However, Brownian LPP, in its unscaled coordinates, underlies everything we do.  Lemma~\ref{l.densepolyunique}(1)  is an example: the basic result on which it depends is Lemma~\ref{l.severalpolyunique}, which is its unscaled counterpart; in fact, the result is more general than this, treating uniqueness of certain multi-geodesics.
Although, in keeping with our emphasis on scaled coordinates,  this lemma appears only at the end of the paper, in Appendix~\ref{s.polyunique},  it is an interesting basic result concerning Brownian LPP. 

We also need~\cite[Lemma~$2.2(1)$]{ModCon}.
\begin{lemma}\label{l.basic}
Let $(n,t_1,t_2) \in \N \times \R^2_<$
be a compatible triple. 
 The random function $(x,y) \to \weight_{n;(x,t_1)}^{(y,t_2)}$,
which is defined on the set of $(x,y) \in \R^2$ 
satisfying $y \geq x - 2^{-1} n^{1/3} \tot$, is continuous almost surely.
\end{lemma}
\subsection{Polymer fluctuation}

The rerooting procedure to resolve the second difficulty in the rough guide will make use of an important aspect of polymer geometry: polymers have H\"older-$2/3-$ regularity. We now quote two results from~\cite{NonIntPoly} to this effect.

Some notational preparation is needed.
Let $(n,t_1,t_2) \in \N \times \R^2_<$ be a compatible triple, and let $x,y \in \R$.
The polymer $\rho_{n;(x,t_1)}^{(y,t_2)}$ has been defined to be a subset of $\R \times [t_1,t_2]$ containing $(x,t_1)$ and $(y,t_2)$, but really as $n$ rises towards infinity, it becomes more natural to seek to view it as a random function that maps its lifetime $[t_1,t_2]$ to the real line.
In choosing to adopt this perspective, we will abuse notation: taking $t \in [t_1,t_2]$, we will speak of the value  $\rho_{n;(x,t_1)}^{(y,t_2)}(t) \in \R$, as if the polymer were in fact a function of $[t_1,t_2]$.
Some convention must be adopted to resolve certain microscopic ambiguities as we make use of this new notation, however.
First, we will refer to  $\rho_{n;(x,t_1)}^{(y,t_2)}(t)$ only when $t \in [t_1,t_2]$ satisfies $nt \in \Z$, a condition that ensures that the intersection of the set $\rho_{n;(x,t_1)}^{(y,t_2)}$  with the line at height $t$ takes place along a horizontal planar interval.

Second, we have to explain which among the points in this interval $\rho_{n;(x,t_1)}^{(y,t_2)} \cap \{  (\cdot,t): \cdot \in \R  \}$
 we wish to denote by $\rho_{n;(x,t_1)}^{(y,t_2)}(t)$.
 To present and explain our convention in this regard, we let
   $\ell_{(x,t_1)}^{(y,t_2)}$ denote the planar line segment whose endpoints are $(x,t_1)$ and $(y,t_2)$. Adopting the same perspective as for the polymer, we abuse notation to view  $\ell_{(x,t_1)}^{(y,t_2)}$ as a function from $[t_1,t_2]$ to $\R$, so that  $\ell_{(x,t_1)}^{(y,t_2)}(t) = \tot^{-1} \big( (t_2 - t) x + (t - t_1)y \big)$.
   
Our convention will be to set  $\rho_{n;(x,t_1)}^{(y,t_2)}(t)$ equal to~$z$ where $(z,t)$ is that point in the horizontal segment  $\rho_{n;(x,t_1)}^{(y,t_2)} \cap \{  (\cdot,t): \cdot \in \R  \}$ whose distance from   $\ell_{(x,t_1)}^{(y,t_2)}(t)$ is maximal. (An arbitrary tie-breaking rule, say $\rho_{n;(x,t_1)}^{(y,t_2)}(t) \geq  \ell_{(x,t_1)}^{(y,t_2)}(t)$, resolves the dispute if there are two such points.)
The reason for this very particular convention is that our purpose in using it is to explore, in the soon-to-be-stated Theorem~\ref{t.polyfluc},
upper bounds on the probability of large fluctuations between the polymer  $\rho_{n;(x,t_1)}^{(y,t_2)}$ and the line segment  $\ell_{(x,t_1)}^{(y,t_2)}$ that interpolates the polymer's endpoints.
Our convention ensures that the form of the theorem would remain valid were any other convention instead adopted.     

In order to study the intermediate time $(1-a)t_1 + at_2$ (in the role of $t$ in the preceding), we now  let $a \in (0,1)$ and impose that $a \tot n \in \Z$: doing so ensures that, as desired, $t \in n^{-1} \Z$, where $t = (1-a)t_1 + at_2$.

Consider also $r > 0$. Define the {\em polymer deviation regularity} event~\hfff{polydevreg} 
$$
\pdr_{n;(x,t_1)}^{(y,t_2)}\big(a,r\big) \, = \, \bigg\{  \,  \Big\vert \, \rho_{n;(x,t_1)}^{(y,t_2)} \big(  (1-a) t_1 +  a t_2 \big) - \ell_{(x,t_1)}^{(y,t_2)}  \big(  (1-a) t_1 +  a t_2 \big) \, \Big\vert \, \leq \, r t_{1,2}^{2/3} \big(  a \wedge (1-a) \big)^{2/3} \bigg\} \, ,
$$
where $\wedge$ denotes minimum. For example, if $a \in (0,1/2)$,
the polymer's deviation from the interpolating line segment, at height $(1-a)t_1 + at_2$ (when the polymer's journey has run for time $a\tot$), is measured in the natural time-to-the-two-thirds  scaled units obtained by division by $(a \tot)^{2/3}$, and compared to the given value $r > 0$.

For intervals $I,J \subset \R$, we extend this definition by setting
$$
 \pdr_{n;(I,t_1)}^{(J,t_2)}\big(a,r\big)   = \bigcap_{x \in I, y \in J}  \pdr_{n;(x,t_1)}^{(y,t_2)}\big(a,r\big) \, .
$$ 
The perceptive reader may notice a problem with the last definition. 
The polymer $\rho_{n;(x,t_1)}^{(y,t_2)}$
is well defined almost surely for given endpoints, but this property is no longer assured as the parameters vary over $x \in I$ and $y \in J$.
In the case of exceptional $(x,y)$ where several $n$-polymers
move from $(x,t_1)$ to $(y,t_2)$, we interpret $\rho_{n;(x,t_1)}^{(y,t_2)}$ as the union of all these polymers, for the purpose of defining $\rho_{n;(x,t_1)}^{(y,t_2)}(t)$. This convention permits us to identify worst case behaviour, so that the {\em complementary} event $\neg \, \pdr_{n;(I,t_1)}^{(J,t_2)}\big(a,r\big)$ is triggered by a suitably large fluctuation on the part of any concerned polymer.

The next two results are~\cite[Theorem~$1.5$]{NonIntPoly}  and~\cite[Proposition~$1.6$]{NonIntPoly}. In essence, they both assert that a unit lifetime polymer fluctuates, in a short initial or final period of duration $a$, by more than $a^{2/3}r$, (with $r$ capable of assuming a very broad range of positive values), only with probability
at most $\exp \big\{ - \Theta(1) r^{3/4} \big\}$, uniformly in high choices of $n \in \N$. In addition, Theorem~\ref{t.polyfluc} asserts this uniformly over polymers that begin and end in a unit-order interval.
\begin{theorem}\label{t.polyfluc}
Let $(n,t_1,t_2) \in \N \times \R^2_<$ be a compatible triple,  and let $x,y \in \R$. 
\begin{enumerate}
\item
Let  $a \in \big[1 - 10^{-11} c_1^2 , 1 \big)$ satisfy $a \tot \in n^{-1} \Z$.  
 Suppose that $n \in \N$ satisfies 
$$
 n \tot   \geq  \max \bigg\{ 
 10^{32} (1-a)^{-25} c^{-18} \, \, , \, \, 
 10^{24} c^{-18} (1-a)^{-25} \vert x - y \vert^{36} \tot^{-24} 
 \bigg\} \, .
$$
Let $r > 0$ be a parameter that satisfies
$$ 
 r  \geq 
 \max \bigg\{  10^9 c_1^{-4/5} \, \, , \, \,  15 C^{1/2} \, \, , \, \,  87(1-a)^{1/3} \tot^{-2/3} \vert x - y \vert   \bigg\}
 $$
and
$r \leq 3  (1-a)^{25/9} n^{1/36} \tot^{1/36}$. 

Writing $I = \big[ x, x +  t_{1,2}^{2/3} (1-a)^{2/3} r \big]$
and  $J = \big[ y, y +  t_{1,2}^{2/3} (1-a)^{2/3} r \big]$, we have that
$$
\PP \Big( \neg \, \pdr_{n;(I,t_1)}^{(J,t_2)}\big(a,2r\big) \, \Big) \leq  
 44  C r    \exp \big\{ - 10^{-11} c_1  r^{3/4} \big\} 
 \, .
$$
\item The same statement holds verbatim when appearances of $a$ are replaced by $1-a$.
\end{enumerate}
\end{theorem}

\begin{proposition}\label{p.polyfluc}
Let $(n,t_1,t_2) \in \R^2_<$ be a compatible triple.  Let $x,y \in \R$ and   let $a \in \big[1 - 10^{-11} c_1^2 , 1 \big)$ satisfy $a \tot \in n^{-1} \Z$. Suppose that 
$$
 n \tot   \geq  \max \bigg\{ 
 10^{32} (1-a)^{-25} c^{-18} \, \, , \, \, 
 10^{24} c^{-18} (1-a)^{-25} \vert x - y \vert^{36} \tot^{-24} 
 \bigg\} \, .
$$
Let $r > 0$ be a parameter that satisfies
$$ 
 r  \geq 
 \max \bigg\{  10^9 c_1^{-4/5} \, \, , \, \,  15 C^{1/2} \, \, , \, \,  87(1-a)^{1/3} \tot^{-2/3} \vert x - y \vert   \bigg\}
 $$
and
$r \leq 3  (1-a)^{25/9} n^{1/36} \tot^{1/36}$. 
Then
$$
 \PP \Big(  \neg \, \pdr_{n;(x,t_1)}^{(y,t_2)}\big(a,r\big) \, \Big) \leq 
22  C r    \exp \big\{ - 10^{-11} c_1  r^{3/4} \big\}    \,  .
$$
\end{proposition}

\subsection{Control on the fluctuation of line-to-point polymers}

To implement the rough guide, we will want to know that $f$-rewarded line-to-point polymers that end at time one in $[-1,1]$ begin at time zero at unit-order distance from this interval.
We now quote a result, \cite[Lemma~$6.2$]{ModCon}, that, when allied with the polymer sandwich Lemma~\ref{l.sandwich}, will make such an assertion.

Recalling Definition~\ref{d.if}, let  $\ovbar\coninit \in (0,\infty)^3$ and $f \in \initcond_{\ovbar\coninit}$.
For $R \geq 0$, define the event~\hfff{regfluc}
$$
\regfluc_{n;(*:f,0)}^{\big( \{-1,1 \}, 1 \big)}(R) = \Big\{ \rho_{n;(*:f,0)}^{(-1,1)}(0) \geq  -(R+1)  \, , \,   \rho_{n;(*:f,0)}^{(1,1)}(0) \leq   R + 1  \Big\} \, .
$$
\begin{lemma}\label{l.regfluc}
Let $n \in \N$, $R > 0$ and $\ovbar\coninit \in (0,\infty)^3$ satisfy 
$$
  n \geq  c^{-18} \max \Big\{  (\coninit_2 + 1)^9 \, , \,   10^{23} \coninit_1^9     \, , \, 3^{9}  \Big\} \, ,
$$
 $$
    R  \geq \max \Big\{ \,    39 \coninit_1  \,  , \, 5  \, , \,   3 c^{-3}  \, , \, 2 \big( (\coninit_2 + 1 )^2 +  \coninit_3 \big)^{1/2} \, \Big\} \, , 
  $$
  and
 $R \leq 6^{-1}c n^{1/9}$.
Then, for any $f \in \initcond_{\ovbar\coninit}$,
$$
 \PP \Big( \neg \, 
\regfluc_{n;(*:f,0)}^{\big( \{-1,1 \}, 1 \big)}(R)  \Big) \leq 38 R \rsC       \exp \big\{ -  2^{-6} \rsc   R^3 \big( 2^{ - 1/2} - 2^{-1} \big)^{3/2}  \big\}  \, .
$$
\end{lemma}


\subsection{Conventions governing the presentation of upcoming proofs}

We close this section of tools by explaining two such conventions, even though they will become operative only when we begin giving proofs in Section~\ref{s.rarelatecoal}.

\subsubsection{Boldface notation for quoted results}\label{s.boldface}

During the upcoming proofs, we will naturally be making use of the various tools that we have recalled in this section. The statements of such results involve several parameters, in several cases including $(n,t_1,t_2)$, spatial locations $x$ and $y$, and positive real parameters such as $r$.
We will employ a notational device that will permit us to disregard notational conflict between the use of such parameters in the context of the ongoing proof in question, and their role in the statements of quoted results. When specifying the parameter settings of a particular application, we will allude to the parameters of the result being applied in boldface, and thus permit the reuse of the concerned symbols.

\subsubsection{The role of hypotheses invoked during proofs}\label{s.rolehyp}

When we quote results in order to apply them, we will take care, in addition to specifying the parameters according to the just described convention, to indicate explicitly what the conditions on these parameters are that will permit the quoted result in question to be applied. Of course, it is necessary that the hypotheses of the result being proved imply all such conditions. The task of verifying that the hypotheses of a given result are adequate for the purpose of obtaining all conditions needed to invoke the various results used during its proof 
may be called the calculational derivation of the result in question. 
This derivation
is necessary, but also in some cases lengthy and unenlightening: a succession of trivial steps. 
In the case of two results proved in this article, Proposition~\ref{p.latecoal} and Lemma~\ref{l.normalcoal}, 
we have chosen to separate the calculational derivation from the rest of the proof. 
These derivations may be found in Appendix~\ref{s.calcder} at the end of the paper. 
After both of the two results' statements a remark is made intended for the reader who is disinclined to review or verify the associated calculational derivation. Each remark presents the preceding result's hypotheses in a simpler and less explicit form and offers a little guidance about how to review the ensuing proof.

\section{The rough guide elaborated: how rerooting will be carried out}\label{s.reroot}

We have assembled the key tools for the proof of Theorem~\ref{t.unifpatchcompare}. In order to present a clear heuristic for how the proof will proceed, we now revisit the rough guide from Section~\ref{s.roughguide}. Recall that we proposed that the second objection, that root locations may be exceptional, may be resolved by rerooting the trees in the $f$-rewarded polymer forest to nearby elements in an $\e$-mesh lying in the $x$-axis, where $\e$ is small.

There are some problems in implementing this idea. The weight profile from a given root and from an adjacent $\e$-mesh point will differ by a constant only if the concerned polymers that emanate from these two locations merge quickly enough, before the branch time at which the polymers running from the root to different locations in the canopy in question split apart. This consideration puts downward pressure on $\e$, in the sense that this assertion will be become easier to prove if $\e$ is small. However,
when we come to say that the polymer weight profile $y \to \weight_{n;(x,0)}^{(y,1)}$ for $x$ belonging to the $\e$-mesh withstands a strong comparsion to Brownian bridge by invoking Theorem~\ref{t.bridge},
our comparison will deteriorate as $\e$ decreases, because it will be necessary to sum estimates with a union bound over the order $\e^{-1}$ points in a bounded portion of the $\e$-mesh. This will put an opposing, upward, pressure on our choice of the parameter $\e$.

Naturally a tension exists between these downward and upward pressures on the value of $\e$. We resolve the tension by lessening the downward pressure, considering a formulation of the assertion of constancy for the difference of weight profiles associated to rerooting that can be proved to have a probability of working that converges rather rapidly to one in the limit of $\e \searrow 0$. 
Recall the right sketch in Figure~\ref{f.manytrees}, and that we have considered the root $r$  of a tree~$T$ in the $f$-rewarded polymer forest,
as well as the left-adjacent $\e$-mesh point $x$. In proposing to reroot the tree from a beginning at $(r,0)$ to one at $(x,0)$,
we want to ensure, in order to achieve the desired constant difference of the associated polymer weight profiles, that the coalescence time of the polymer pair running to a given canopy point from $(r,0)$ and from $(x,0)$ occurs {\em before} the branch time of the tree~$T$.

Since $x$ and $r$ differ by at most the small quantity $\e > 0$, it may be expected that this merging occurs quickly; in fact, the order of time at which it occurs is $\e^{3/2}$.
In order to work with such a merging event that has a provably high probability, we will modify the trees we work with, in order that their branch times are necessarily close to one.

\begin{figure}[ht]
\begin{center}
\includegraphics[height=10cm]{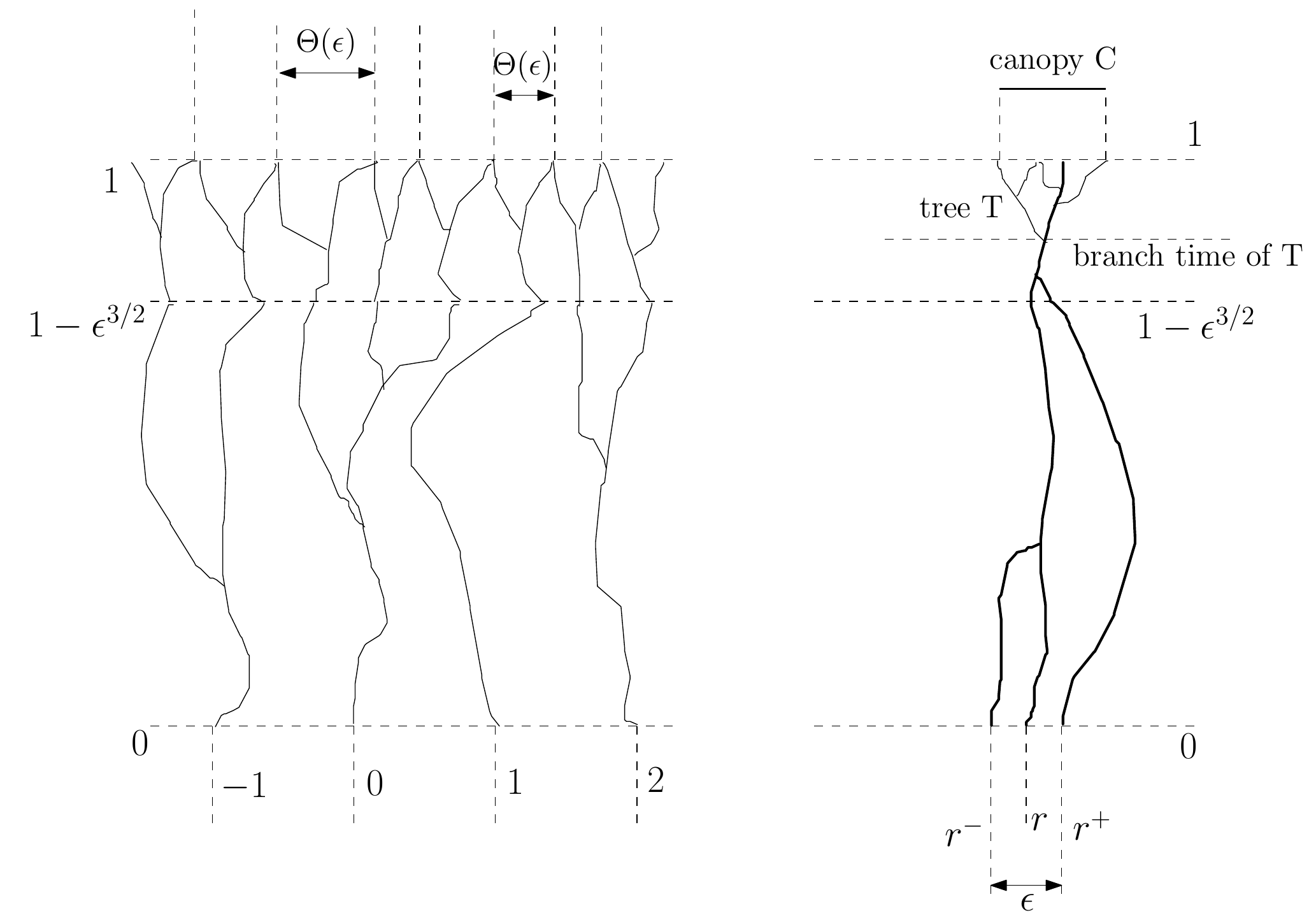}
\caption{{\em Left:} The $f$-rewarded polymer forest corresponding to the $f$ that is zero on $\Z$, and otherwise minus infinity, is depicted. The tree tops are cut along the line at height $1 - \e^{3/2}$, with the boundary points of the resulting new canopies indicated by the vertical dashed lines that touch the height one line. {\em Right:}
The $f$-rewarded polymer tree $T$ associated to one of these new canopies $C$. If the canopy does not experience late coalescence then, for any $y \in C$, one of the consecutive pairs in the emboldened triple of polyers must merge before time $1 - \e^{3/2}$. In this instance, the polymer pair from $r^-$ and $r$ so merge, raising the prospect of rerooting the tree $T$ from its root $r$ to the nearby mesh element $r^- \in \e\Z$.}
\label{f.newcanopy}
\end{center}
\end{figure}

To specify the new trees, picture the $f$-rewarded polymer forest sitting inside the strip $\R \times [0,1]$, and draw a horizontal line of height $1 - \e^{3/2}$ through the strip. The line cuts the trees into smaller tree tops sitting in the strip $\R \times [1-\e^{3/2},1]$, each with its own canopy. These new canopies, which we will be introduced formally as $(n,*:f,1-\e^{3/2})$-canopies later on, are intervals that subdivide the original collection of canopies into finer pieces. Each new canopy has a portion of an $f$-rewarded polymer tree associated to it, with a root at time zero, as before. These new trees have the merit of having a branch time that necessarily exceeds $1 - \e^{3/2}$. Consistently with the two-thirds power law that governs polymer fluctuation, the lengths of the new canopies are of order $\big(\e^{3/2}\big)^{2/3} = \e$, something illustrated in the left sketch of Figure~\ref{f.newcanopy}.

Consider now one of the new canopies~$C$, to which is associated a time-zero root $r$. Working again with the $\e$-mesh $\e\Z$ embedded in the $x$-axis, we will  
 consider not only the mesh point to the left of~$r$, but also the one to the right: call these points $r^-$ and $r^+$. For any canopy location $y \in C$, we may consider the triple of polymers $\rho_{n;(r^-,0)}^{(y,1)}$,  $\rho_{n;(r,0)}^{(y,1)}$ and  $\rho_{n;(r^+,0)}^{(y,1)}$.  
 As we have mentioned, we may expect the merging time between 
  any pair of these polymers to be rather small. We now have the advantage that, with the aim of proving such merging times are less than the branch time of the tree, we need merely prove that they are typically less than $1 - \e^{3/2}$: that is, not necessarily early, but not very late. Moreover, in order to work with a circumstance of suitably high probability, we will consider an event in which it is demanded that only one pair of the concerned polymers coalesces by this time.
  
  Thus, we want to investigate the probability of the 
   event that at least one of the coalescence times between the first and second of the polymers, and between the second and the third, is at most $1 - \e^{3/2}$.

If there exists $y \in C$ for which this presumably typical event fails to occur, we say that the $(n,*:f,1-\e^{3/2})$-canopy $C$  experiences {\em late} coalescence.
How may we bound above the probability of late coalescence? When this event takes place, there exists an element $y \in C$ such that the above triple of polymers remains disjoint throughout the strip $\R \times [0,1-\e^{3/2}]$. These three polymers are heading to a shared endpoint at time one, namely $y$. The two-thirds power law for polymer fluctuation, given rigorous expression by Proposition~\ref{p.polyfluc}, implies that, typically, the three polymers at time $1 - \e^{3/2}$ will lie within a distance of order $\e$ of $y$. 
Since the polymers begin at $r^-$, $r$ and $r^+$, which are locations lying in a given interval of length $\e$, we see that the polymers restricted to the strip in question begin and end within distance $\e$ of each other, and remain pairwise disjoint.

When the length-$\e$ intervals at which such a triple of polymers begin and end are given, we may use Theorem~\ref{t.disjtpoly} to find an upper bound on the probability of the existence of such disjoint polymers. Taking $m=3$ in that result, we find an upper bound of $\e^{(3^2 - 1)/2 + o(1)} = \e^{4 + o(1)}$.

Now if we use this approach to determine an upper bound on the probability that there exists an $(n,*:f,1-\e^{3/2})$-canopy  that intersects a given unit interval and for which the late coalescence event occurs, we must sum the $\e^{4 + o(1)}$ bound over an order of $\e^{-1}$ such canopies (each, after all,  has length of order $\e$) as well as an order of $\e^{-1}$ time-zero length-$\e$ intervals that capture the root location associated to a given canopy. Thus an entropy term of $\e^{-2}$ appears, leading to an upper bound of the form $\e^{2 + o(1)}$ on this event of a uniform absence of late coalescence.

The battle between $4$ and $-2$ has ended with a positive outcome, showing, as desired, the rarity of any instance of late coalescence. Note that the use of the advanced time $1-\e^{3/2}$ and of a triple, rather than merely a pair, of polymers is needed for this to happen.

\begin{figure}[ht]
\begin{center}
\includegraphics[height=7cm]{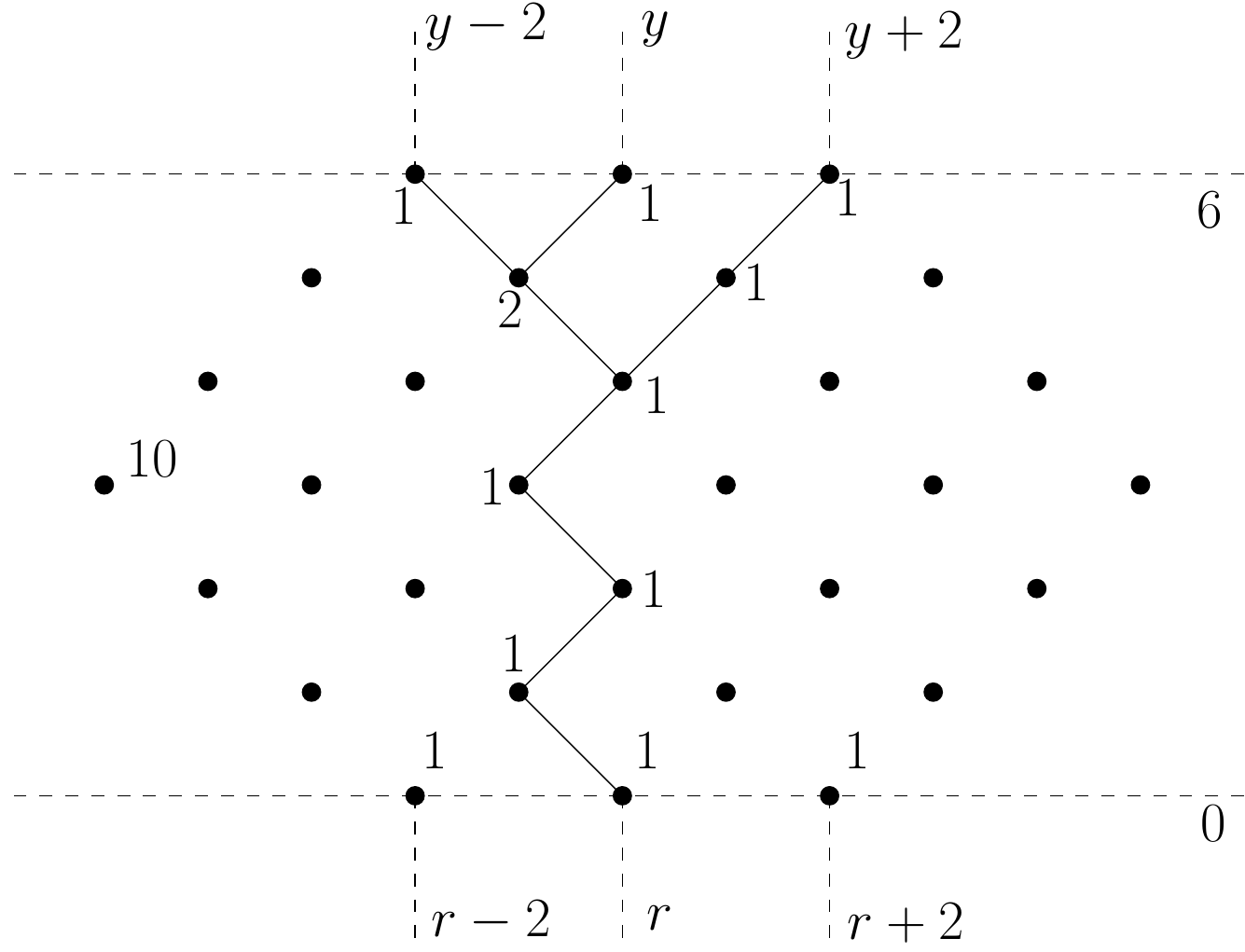}
\caption{In this discrete example of an LPP model, directed paths move from height zero to height six by means of diagonal adjacency. A path's value is naturally the sum of the lattice values that it finds in its way. All lattice values are zero except those indicated. Depicted is a geodesic tree~$T$, rooted at $(r,0)$, with a three element canopy $\{ y-2,y, y+2 \}$ at height six. (As it happens, $y=r$.) The tree has branch time equal to four. The canopy does not experience late coalescence, in the sense that, for each $z \in \{ y-2,y, y+2 \}$, either the pair of geodesics from $(r-2,0)$ and $(r,0)$ to $(z,6)$ merges by time four, or the pair from $(r,0)$  and $(r+2,0)$ does. Despite this, the difference between the geodesic value profiles that map $z \in \{ y-2,y, y+2 \}$ to the value of the geodesic from $(u,0)$ to $(z,6)$, for $u = r$ and $u = r-2$, is not constant. Indeed, the geodesic from $(r-2,0)$ to $(y-2,6)$ collects the value of ten located at $(r-5,3)$, while such an opportunity is not available to the geodesic from $(r,0)$ to $(y-2,6)$. The canopy can however be split into two subintervals so that this problem of non-constancy does not occur.}
\label{f.discrete}
\end{center}
\end{figure}

The rerooting procedure will then be applied in this typical scenario of absence of late coalescence. For any location $y$ belonging to one of the new canopies $C$, either the polymers $\rho_{n;(r^-,0)}^{(y,1)}$ and $\rho_{n;(r,0)}^{(y,1)}$ meet before time $1 - \e^{3/2}$, or the polymers $\rho_{n;(r,0)}^{(y,1)}$ and $\rho_{n;(r^+,0)}^{(y,1)}$ do. Suppose for definiteness that the former case applies. If we vary $y$ within the order-$\e$ length interval~$C$, the various polymers emanating from $(r,0)$ will share a common course until their branch time, which exceeds $1 - \e^{3/2}$.  Thus, the weight profiles $y \to \weight_{n;(r^-,0)}^{(y,1)}$ and $y \to \weight_{n;(r,0)}^{(y,1)}$ would appear to differ merely by a random constant as $y$ varies over~$C$: see the right sketch in Figure~\ref{f.newcanopy}. If this can be rigorously demonstrated, it will allow us to treat the  $(n,*:f,1-\e^{3/2})$-canopies   as the patches in our Brownian patchwork quilt, with the associated fabric sequence element being either  $y \to \weight_{n;(r^-,0)}^{(y,1)}$ or $y \to \weight_{n;(r^+,0)}^{(y,1)}$. In fact, there is a perhaps surprising counterexample to the assertion that this difference of weight profiles is necessarily constant. An example of a difficulty is given in Figure~\ref{f.discrete}. When we implement this approach rigorously, we will employ a trick whereby any given  $(n,*:f,1-\e^{3/2})$-canopy will be broken into two, in a way that ensures the desired constancy of the two weight profiles on the resulting pieces.

It is these rerooting ideas that we will develop to resolve the second objection in the rough guide. The quantity $\e > 0$ will be chosen to be random in fact: we will search for a choice on successively smaller dyadic scales until the uniform absence of late coalescence is achieved. Sometimes $\e$ will end up being very small, and then the $\e$-mesh locations, such as $r^-$ in the preceding discussion, will number a large quantity of order $\e^{-1}$. This will cause some deterioration in comparison to Brownian bridge when weight profiles such as  $y \to \weight_{n;(r^-,0)}^{(y,1)}$  are so compared by means of Theorem~\ref{t.bridge}. The result is that while this theorem implies that the weight profile  $y \to \weight_{n;(x,0)}^{(y,1)}$ for a given choice of $x$ withstands $L^{\infty - }$ comparison to Brownian bridge, the fabric sequence elements in our Brownian quilt representation of the $f$-rewarded weight profiles will merely be proved to withstand $L^{3-}$ comparison to Brownian bridge.

\section{Rarity of late coalescence}\label{s.rarelatecoal}

In this section, we specify the late coalescence event. We then state and prove Proposition~\ref{p.latecoal}; the proof will use Lemma~\ref{l.maxpolytriple}. The proposition asserts the discussed $\e^{2 + o(1)}$ 
upper bound on the probability of uniform absence of late coalescence. The lemma
asserts that it is typical in the case of this absence that there exists a triple of disjoint polymers crossing the strip $\R \times [0,1-\e^{3/2}]$ with $\e$-distant endpoints. Thus, the role of the lemma is in essence to reduce the proposition to an application of Theorem~\ref{t.disjtpoly}.

For $n \in \N$, $\e > 0$ and $K > 0$, define the {\em late coalescence} event~\hfff{latecoal}
$\latecoal_{n;([-K,K],0)}^{([-1,1],1)}(\e)$
to be the event that  there exist $x \in [-K,K]$ and  $y \in [-1,1]$ such that,
\begin{itemize}
\item  setting $i = \lfloor x \e^{-1} \rfloor$,
each of the pairs $(i\e,y)$, $(x,y)$ and $\big((i+1)\e,y \big)$ lies in $\polyunique_{n;0}^1$, which is to say, each of the polymers 
$\rho_{n;(i\e,0)}^{(y,1)}$, $\rho_{n;(x,0)}^{(y,1)}$ and $\rho_{n;\big((i+1)\e,0 \big)}^{(y,1)}$  is well defined;
\item
and, for any pair among this triple of polymers, the intersection of the pair contains no point whose vertical coordinate is less than or equal to $1 - \e^{3/2}$.
\end{itemize}
We write $\nolatecoal_{n;([-K,K],0)}^{([-1,1],1)}(\e)$ for the complementary event $\neg \, \latecoal_{n;([-K,K],0)}^{([-1,1],1)}(\e)$.
This last event is the presumably typical event under which there is no late coalescence in the concerned unit-order region that we discussed in the preceding heuristic.

\begin{proposition}\label{p.latecoal}
Let $K \in [2,\infty)$, $\e > 0$ and $n \in \N$ satisfy $n \e^{3/2} \in \N$, 
$$
  \e \leq \min \Big\{   \,      (\eta_0)^{72} \, , \, 10^{- 1342}  c_3^{44}  \big( K  + 2 \big)^{-8}  \,  , \, \exp \big\{ - 2 C^{3/8} \big\} \, \Big\} \, ,
$$ 
and
$$
     n  \geq 2 \max \bigg\{ \, 2(K_0)^9 \big( \log \e^{-1} \big)^{K_0}
    \, , \,   
 10^{728}     c_3^{-84}  \e^{-222}   (K+2)^{36}  \, , \, 
  a_0^{-9} (K+2)^9 2^{6}  \,  \Bigg\} \, .
$$
Then
$$
 \PP \Big( \latecoal_{n;([-K,K],0)}^{([-1,1],1)}(\e) \Big) \leq  \e^2 \cdot 10^{420}  K  c_3^{-33}      C_3   \big( \log \e^{-1} \big)^{42}     \exp \big\{ \beta_3 \big( \log \e^{-1} \big)^{5/6} \big\} \, . 
$$
(The positive constants $\eta_0$, $K_0$, $a_0$ and $\beta_3$ are each specified either in Theorem~\ref{t.maxpoly} or Theorem~\ref{t.disjtpoly}. See Section~\ref{s.explicit} for the values of $C_3$ and $c_3$ and in regard to the origin of the constant $C$.)
\end{proposition}
\noindent{\em Remark.}
This result is the first of the two whose calculational derivation has been separated, as we discussed in Subsection~~\ref{s.rolehyp}. 
The reader may choose to interpret the proposition's hypotheses as taking, for some large universal constant $\Theta$, the forms  $K \geq 2$, $\e \leq \Theta^{-1} K^{-\Theta}$ and $n \geq K^\Theta \e^{-\Theta}$. 
In the upcoming proof, each bound that is needed is recorded explicitly; each of these may be straightforwardly seen to be compatible with these conditions on $\e$ and $n$.

\noindent{\bf Proof of Proposition~\ref{p.latecoal}.} The next result is an important component.
\begin{lemma}\label{l.maxpolytriple}
Let $r > 0$. For $\e \in (0,2^{-2/3}]$, the occurrence of 
\begin{equation}\label{e.lcassume}
 \latecoal_{n;([-K,K],0)}^{([-1,1],1)}(\e) \cap \pdr_{n;([-K-1,K+1],0)}^{([-1,1],1)}\big(1 - \e^{3/2},r\big)
\end{equation} 
entails the existence of a pair $(i,j)$ of integers satisfying $\vert i \vert \leq K\e^{-1} + 1$ and $\vert j \vert \leq \e^{-1} + 1$ such that
\begin{equation}\label{e.maxpolytriple}
 \maxpoly_{n;\big([i\e,(i+1)\e],0\big)}^{\Big( \big[ j\e-(r+1)\e -  (K+2)\e^{3/2} \, , \, j\e+(r+1)\e +  (K+2)\e^{3/2} \big] ,1 - \e^{3/2} \Big)} \, \geq \, 3 \, .
\end{equation}
\end{lemma}
 This use of the $\pdr$ event requires that $1 - \e^{3/2} \in n^{-1}\Z$. This holds due to the   Proposition~\ref{p.latecoal} hypothesis that $n \e^{3/2} \in \N$.

\noindent{\bf Proof of Lemma~\ref{l.maxpolytriple}.}
Suppose that the event~(\ref{e.lcassume}) occurs.

The occurrence of $\latecoal_{n;([-K,K],0)}^{([-1,1],1)}(\e)$ furnishes  $x \in [-K,K]$ and $y \in [-1,1]$. By setting 
$i = \lfloor x \e^{-1} \rfloor$
and $j = \lfloor y \e^{-1} \rfloor$,
note that the pair $(i,j)$ satisfies the bounds stated in Lemma~\ref{l.maxpolytriple}.
Moreover, the intersections of the polymers 
$\rho_{n;(i\e,0)}^{(y,1)}$, $\rho_{n;(x,0)}^{(y,1)}$ and $\rho_{n;\big((i+1)\e,0 \big)}^{(y,1)}$ 
with the strip $\R \times [0,1-\e^{3/2}]$
form a pairwise disjoint triple of polymers; the three values $i\e$, $x$ and $(i+1)\e$ are thus necessarily distinct.

To confirm that this polymer triple verifies~(\ref{e.maxpolytriple}), the event $\pdr_{n;([-K-1,K+1],0)}^{([-1,1],1)}\big(1 - \e^{3/2},r\big)$ will be invoked to gain understanding of where the  endpoints of these three polymers lie.
Specifically, we will allow  $u$ to take any value in $\{ i\e,x,(i+1)\e \}$, and argue that the upper endpoints satisfy
\begin{equation}\label{e.uinterval}
\big\vert \rho_{n;(u,0)}^{(y,1)} (  1 - \e^{3/2} ) - j\e \big\vert \leq (r+1)\e + (K+2)\e^{3/2} \, .
 \end{equation}
 Since the lower endpoints clearly satisfy $u \in [i\e,(i+1)\e]$,
 (\ref{e.uinterval}) will be sufficient to verify~(\ref{e.maxpolytriple}).
 
To begin arguing for~(\ref{e.uinterval}), note that, since $\e \leq 1$, the event  $\pdr_{n;(u,0)}^{(y,1)}\big(1 - \e^{3/2},r\big)$ occurs whenever  $u \in \{ i\e,x,(i+1)\e \}$.  Since  $1 - \e^{3/2} \geq 1/2$, this entails that 
 $$
   \ \Big\vert \, \rho_{n;(u,0)}^{(y,1)} \big(  1 - \e^{3/2} \big) - \ell_{(u,0)}^{(y,1)}  \big(    1 - \e^{3/2} \big) \, \Big\vert \, \leq \, r  \e 
$$
for such $u$.
Since our choice of $u$ lies in $[-K-1,K+1]$, and $y \in [-1,1]$, the gradient of the planar line segment $\ell_{(u,0)}^{(y,1)}$ is a constant that is at least $(K+2)^{-1}$ in absolute value,
so that  $\big\vert \ell_{(u,0)}^{(y,1)}(1-\e^{3/2}) - y \big\vert \leq (K+2)\e^{3/2}$. We see then that 
 $\big\vert \rho_{n;(u,0)}^{(y,1)} (  1 - \e^{3/2} ) - y \big\vert \leq r\e + (K+2)\e^{3/2}$, which alongside $\vert y - j\e \vert \leq \e$ completes the derivation of~(\ref{e.uinterval}) and thus also of (\ref{e.maxpolytriple}). \qed

\medskip

We let $r > 0$ be a parameter whose value will be specified later in terms of $\e$; it may be useful to note for now that $r$ will take the form of a bounded multiple of a certain power of $\log \e^{-1}$.  For any given pair $(i,j)$ of integers such that $\vert i \vert \leq K\e^{-1} + 1$ and $\vert j \vert \leq \e^{-1} + 1$,
the event in~(\ref{e.maxpolytriple}) is a subset of 
$$
 \maxpoly_{n;\big(
\big[ i\e-(r+1)\e -  (K+2)\e^{3/2} \, , \, i\e+(r+1)\e +  (K+2)\e^{3/2} \big],0\big)}^{\Big( \big[ j\e-(r+1)\e -  (K+2)\e^{3/2} \, , \, j\e+(r+1)\e +  (K+2)\e^{3/2} \big] ,1 - \e^{3/2} \Big)} \, \geq \, 3 \, ,
$$
and the probability of the latter event
may be bounded above by
applying Theorem~\ref{t.disjtpoly}
with parameter settings ${\bf t_1} = 0$, ${\bf t_2} = 1 - \e^{3/2}$,  ${\bf x} = i \e$,
${\bf y} = j \e$, ${\bf m} = 3$, and with ${\bm \e}$ chosen so that
$2 \, ({\bf \tot})^{2/3} {\bm \e} = (r+1)\e +  (K+2)\e^{3/2}$. Since $\e \leq (1 - 2^{-3/2})^{2/3}$, $({\bf \tot})^{2/3} \geq 1/2$; also using ${\bf \tot} \leq 1$, we find that 
$$
\tfrac{1}{2}\, \Big( (r+1) \e +  (K+2)\e^{3/2}   \Big)  \leq  {\bm \e} \leq (r+1) \e +  (K+2)\e^{3/2}  \, .
$$
(This is the first use of boldface notation, explained in Subsection~\ref{s.boldface}.)
Since $r \geq 1$ and $K+2 \leq \e^{-1/2}$, these bounds imply
\begin{equation}\label{e.etwobounds}
\e   \leq  {\bm \e} \leq (r  + 2)\e   \, .
\end{equation}
Using these bounds, we learn from this application of Theorem~\ref{t.disjtpoly} that the event  in~(\ref{e.maxpolytriple}) has probability at most
$$
 (r + 2)^{4} \e^{4} \cdot 
 10^{288}  3^{135}  c_3^{-27} C_3   \big( \log \e^{-1} \big)^{36}     \exp \big\{ \beta_3 \big( \log \e^{-1} \big)^{5/6} \big\} \, .
$$
Since $x - {\bf x}$ and $y - {\bf y}$ belong to $[0,\e)$, we see that $\vert {\bf x} - {\bf y} \vert \leq \vert x - y \vert + \e$.

Regarding the hypotheses that are needed for this application of Theorem~\ref{t.disjtpoly}, note that in view of the latter bound in~(\ref{e.etwobounds}), the condition~(\ref{e.epsilonbound}) is met
provided that 
$$
 (r+2)\e \leq \min \Big\{ \, (\eta_0)^{36} \, , \, 10^{-616}  c_3^{22}  3^{-115} \,  , \, \exp \big\{ - C^{1/4} \big\} \, \Big\} \, ,
$$
while since ${\bf \tot} \geq 2^{-1}$ and  $\vert {\bf x} - {\bf y} \vert \leq \vert x - y \vert + \e \leq K+2$, (\ref{e.nlowerbound})
is met when
$$
    n  \geq 2 \max \bigg\{ \, 2(K_0)^9 \big( \log \e^{-1} \big)^{K_0}
    \, , \, 10^{606}   c_3^{-48}  3^{240}  \rsc^{-36} \e^{-222}  \max \big\{   1  \, , \,   (K+2)^{36} 2^{24}  \big\} \, , \, 
  a_0^{-9} (K+2)^9 2^{6}  \,  \Bigg\} \, .
$$
Finally, the hypothesis that  $({\bf \tot})^{-2/3} \vert {\bf y} - {\bf x} \vert \leq  {\bm \e}^{-1/2} \big( \log {\bm \e}^{-1} \big)^{-2/3} \cdot 10^{-8} c_3^{2/3} 3^{-10/3}$ is ensured by  
$$
  2 \big( K +1 + 2\e \big) \leq (r+2)^{-1/2} \e^{-1/2}  \big( \log (r+2)^{-1}\e^{-1} \big)^{-2/3} \cdot 10^{-8} c_3^{2/3} 3^{-10/3}
$$
since $\vert {\bf x} - {\bf y} \vert \leq K +1 + 2\e$ and  $({\bf \tot})^{2/3} \geq 1/2$, where we also used that the function $u \to u^{-1/2} \big( \log u^{-1} \big)^{-2/3}$ is decreasing on $(0,e^{-4/3})$ alongside $(r+2) \e \leq e^{-4/3}$ and~(\ref{e.etwobounds}).

 We now
 apply Lemma~\ref{l.maxpolytriple},
taking a union bound over the integer pair $(i,j)$ 
that the lemma provides, and using the just obtained upper bound on the probability of~(\ref{e.maxpolytriple}). We learn that  
\begin{eqnarray}
 & & \PP \Big( \latecoal_{n;([-K,K],0)}^{([-1,1],1)}(\e) \cap \pdr_{n;([-K-1,K+1],0)}^{([-1,1],1)}\big(1 - \e^{3/2},r\big) \Big) \nonumber \\
 & \leq &  \big( 2 K \e^{-1} + 3 \big) \big( 2\e^{-1} + 3 \big) 
 \cdot 
 (r + 2)^{4} \e^{4} \cdot 
 10^{288}  3^{135}  c_3^{-27} C_3   \big( \log \e^{-1} \big)^{36}     \exp \big\{ \beta_3 \big( \log \e^{-1} \big)^{5/6} \big\}  \nonumber \\
 & \leq & \e^2 \cdot 3K \cdot 3 \cdot 3^4 r^4  \cdot 
 10^{288}  3^{135}  c_3^{-27} C_3   \big( \log \e^{-1} \big)^{36}     \exp \big\{ \beta_3 \big( \log \e^{-1} \big)^{5/6} \big\} \nonumber \\
 & \leq & \e^2 \cdot 10^{356}  K r^4  c_3^{-27} C_3   \big( \log \e^{-1} \big)^{36}     \exp \big\{ \beta_3 \big( \log \e^{-1} \big)^{5/6} \big\} \, . \label{e.pupperbound}
\end{eqnarray}
In the penultimate inequality, the bounds $\e \leq 1/3$,  $K \geq 1$ and $r \geq 1$ were used in the guise $2K \e^{-1} + 3 \leq 3K \e^{-1}$, $2\e^{-1} + 3 \leq 3 \e^{-1}$ and $(r+2)^4 \leq 3^4 r^4$.

We now find an upper bound on the probability of  
$\neg \, \pdr_{n;([-K-1,K+1],0)}^{([-1,1],1)}\big( 1 - \e^{3/2} , r \big)$
by applying Theorem~\ref{t.polyfluc}.
In this application, we will take ${\bf t_1} = 0$, ${\bf t_2} = 1$,  ${\bf a} = 1 - \e^{3/2}$ and ${\bf r} = r/2$; note that ${\bf a} \geq 1/2$ and that $({\bf \tot})^{2/3} (1-{\bf a})^{2/3} = \e$.  
The product set $[-K-1,K+1] \times [-1,1]$  may be covered by 
$\lfloor \big( 4(K+1) \e^{-1} r^{-1} + 1 \big)  \big( 4  \e^{-1} r^{-1} + 1 \big) \rfloor$ products of the form $[x,x+\e r/2] \times [y,y+\e r/2]$; Theorem~\ref{t.polyfluc} will be applied with $({\bf x},{\bf y})$ ranging over such choices of $(x,y)$. We find then that
\begin{equation}\label{e.pbound}
 \PP \Big( \neg \, \pdr_{n;([-K-1,K+1],0)}^{([-1,1],1)}\big( 1 - \e^{3/2} , r \big) \Big) \leq 64 (K+1) \e^{-2}   
\cdot  22 Cr    \exp \big\{ - 10^{-11} 2^{-3/4} c_1  r^{3/4} \big\} 
\end{equation}
where we also used $r \geq 1$, $K \geq 0$ and $\e \leq 1/4$.

In each of these applications of Theorem~\ref{t.polyfluc},
$\vert {\bf x} - {\bf y} \vert \leq K+2 \leq 2K$ (since $K \geq 2$). Thus, in order that the applications be made,
  it is sufficient that  each of the following conditions is met:  $\e^{3/2} \leq  10^{-11} c_1^2$,
$$
 n    \geq  \max \bigg\{ 
 10^{32} \e^{-75/2} c^{-18} \, \, , \, \, 
 10^{24} 2^{36} c^{-18} \e^{-75/2}  K^{36}  
 \bigg\} \, ,
$$
$$ 
 r  \geq 
 2 \max \bigg\{  10^9 c_1^{-4/5} \, \, , \, \,  15 C^{1/2} \, \, , \, \,  174 \e^{1/2}  K   \bigg\} \, ,
 $$
and
$r \leq 6  \e^{25/6} n^{1/36}$; $n$ is supposed to be a positive integer.

We now make the choice  $r = 10^{44/3} 2^{11/3} c_1^{-4/3}  \big( \log \e^{-1} \big)^{4/3}$ in order that the exponential term in the right-hand side in~(\ref{e.pbound}) equals $\e^4$. We find then that this right-hand side is at most 
$$
 1408 ( K+1)  C \cdot   10^{44/3} 2^{11/3} c_1^{-4/3}  \big( \log \e^{-1} \big)^{4/3}    \e^2 \leq 
 10^{20} K  C  c_1^{-4/3}  \big( \log \e^{-1} \big)^{4/3}    \e^2  
$$
since  $K \geq 1$.
Combining this bound with the upper bound in line~(\ref{e.pupperbound}),  substituting the value of $r$, and using $\e \leq e^{-1}$, $c_3 \leq c_1 \leq 1$ and $C_3 \geq C$,  completes the proof of Proposition~\ref{p.latecoal}. \qed

\section{Polymer coalescence, canopies and intra-canopy weight profiles}

In order to implement the rerooting plan elaborated in Section~\ref{s.reroot}, we need to describe in more precise terms what we mean by the canopy associated to an $f$-rewarded polymer tree, 
or rather the revised notion of an  $(n,*:f,1-\e^{3/2})$-canopy,
and set down notation for coalescence times of polymers. Importantly, we need to make rigorous sense of the notion that the polymer weight profile for polymers emerging from a given root differs by a random constant from this profile where polymers emerge from a nearby mesh point. This task must overcome the challenge presented by the counterexample in Figure~\ref{f.discrete}, which shows that two profiles can fail to differ by a constant even when all the concerned polymers end in a common canopy that does not experience late coalescence. 

This section is devoted to resolving these difficulties. After setting polymer and coalescence notation, we present in Lemma~\ref{l.ec} and its aftermath the concept of an $(n,*:f,s)$-canopy, where here $s$ denotes a time between zero and one. This new notion of canopy refers to an interval of $y$-values such that all the $f$-rewarded line-to-point polymers that end at time one at points~$y$, when viewed in decreasing time, have coalesced by time $s$. This is a stronger notion than that of canopy, which corresponds to $s = 0$. The new canopies are subsets of the old ones, and each is associated to a root. Setting $s = 1 - \e^{3/2}$,  we obtain the revised notion of canopy introduced informally in Section~\ref{s.reroot}. It is our aim to prove that the $f$-rewarded polymer weight profile on any given $(n,*:f,1-\e^{3/2})$ -canopy (with root $r$) differs by a random constant from the profile where we reroot to either the left or right adjacent mesh point $r^-$ or $r^+$. It is at this point that the problem posed by the Figure~\ref{f.discrete} counterexample rears its head. Lemma~\ref{l.middlehighway} is the device permitting a solution of this difficulty. It allows us in Definition~\ref{d.canopyspecialpoint} to introduce a certain {\em special point} into any given   $(n,*:f,1-\e^{3/2})$ -canopy. The special point splits the canopy into two pieces. In Lemma~\ref{l.canopytwopieces}, we establish that, provided that the canopy does not experience late coalescence, the $f$-rewarded weight profile on any such piece differs by the desired random constant from at least one of the two weight profiles of polymers emanating from the neighbouring mesh points $r^-$ and~$r^+$.

\begin{definition}
Let $(x_1,x_2) \in \R^2_\leq$ and $y \in \R$ be such that $(x_1,y)$ and $(x_2,y)$ are elements in  $\polyunique_{n;0}^{1}$.
Define the {\em forward coalescence time}~\hfff{forwardcoal} 
$$
 \tau_{n;\big(  \{ x_1,x_2  \} , 0 \big)}^{\uparrow;( y , 1)} \in [0,1]
$$
from the pair $\{x_1,x_2\}$ at time $0$ to  the element $y$ at time $1$ to be the infimum of the vertical coordinates of points in the intersection $\rho_{n;(x_1,0)}^{(y,1)} \cap \rho_{n;(x_2,0)}^{(y,1)}$.
\end{definition}

Suppose that $f \in \initcond_{\ovbar\coninit}$
for some $\ovbar\coninit \in (0,\infty)^3$.
Let $y \in [-1,1]$. First suppose that $y \in \polyunique_{n;(*:f,0)}^1$, so that the polymer  $\rho_{n;(*:f,0)}^{(y,1)}$ is uniquely defined. 
For the purpose of stating the next result (the notation is not reserved after that), set  $x = \rho_{n;(*:f,0)}^{(y,1)}(0)$ and $i = \lfloor \e^{-1} x \rfloor$.
Second, suppose that
 $(i\e,y), \big((i+1) \e , y \big) \in \polyunique_{n;0}^{1}$. The polymers   $\rho_{n;(i\e,0)}^{(y,1)}$ and  $\rho_{n;((i+1)\e,0)}^{(y,1)}$ are also then well defined. 
 When $y \in [-1,1]$ is such that both of these suppositions is valid, so that all three polymers are well defined, we call $y$ an $(n,\e,f)$-triple uniqueness point; note that Lemma~\ref{l.densepolyunique} ensures that almost every $y \in [-1,1]$ verifies this condition. 

\begin{lemma}\label{l.lcrf}
Let $\ovbar\coninit \in (0,\infty)^3$ and 
suppose that $f \in \initcond_{\ovbar\coninit}$. Let $n \in \N$  $n \in \N$ satisfy
$n > 2^{-3/2} \coninit_1^3 \vee 8 (\coninit_2  + 1)^3$, and let $\e \in (0,1)$.
Suppose that  $\nolatecoal_{n;([-K,K],0)}^{([-1,1],1)}(\e) \cap \regfluc^{(\{-1,1\},1)}_{n;(*,0)}(K-1)$ occurs. Let $y \in [-1,1]$ be an $(n,\e,f)$-triple uniqueness point, and associate $(x,i)$ to $y$ as we just did.
Then    
at least one of the quantities $\tau_{n;\big(\{ i \e,x \} , 0\big)}^{\uparrow;(y,1)}$ and $\tau_{n;\big(\{ x, (i+1) \e \} , 0\big)}^{\uparrow;(y,1)}$ is at most $1 - \e^{3/2}$.
\end{lemma}
\noindent{\bf Proof.} 
The $f$-rewarded line-to-point unit-lifetime polymers that end at $(-1,1)$ and $(1,1)$ are almost surely unique by Lemma~\ref{l.densepolyunique}(2). Thus the polymer sandwich Lemma~\ref{l.sandwich} shows that the occurrence of
$\regfluc^{(\{-1,1\},1)}_{n;(*:f,0)}(K-1)$ entails that  $x \in [-K,K]$.
Each of the polymers $\rho_{n;(i\e,0)}^{(y,1)}$, 
 $\rho_{n;(*,0)}^{(y,1)} = \rho_{n;(x,0)}^{(y,1)}$ and $\rho_{n;\big((i+1)\e,0\big)}^{(y,1)}$ is well defined.
When  $\nolatecoal_{n;([-K,K],0)}^{([-1,1],1)}(\e)$ occurs, either the first and second, or the second and third, of these polymers have a point of intersection whose vertical coordinate is at most $1 - \e^{3/2}$.
Thus, either  $\tau_{n;\big(\{ i \e,x \} , 0\big)}^{\uparrow;(y,1)}$ or  $\tau_{n;\big(\{ x , (i+1) \e \} , 0\big)}^{\uparrow;(y,1)}$  is at most $1 - \e^{3/2}$. \qed

\begin{definition}
Let $x \in \R$ and $(y_1,y_2) \in \R^2_\leq$ be such that $(x,y_1)$ and $(x,y_2)$ are elements in  $\polyunique_{n;0}^{1}$.
Define the {\em backward coalescence time}~\hfff{backwardcoal} 
$$
 \tau_{n;(x,0)}^{\downarrow;(\{y_1,y_2\},1)} \in [0,1]
$$
from the pair $\{y_1,y_2\}$ at time $1$ to  the element $x$ at time $0$ to be the supremum of the vertical coordinates of points in the intersection $\rho_{n;(x,0)}^{(y_1,1)} \cap \rho_{n;(x,0)}^{(y_2,1)}$.

Let $n \in \N$ and let $f \in \initcond_{\ovbar\coninit}$ for some $\ovbar\coninit \in (0,\infty)^3$. Suppose that $(y_1,y_2) \in \R^2_<$  satisfies $y_1,y_2 \in \polyunique_{n;(*:f,0)}^{t}$.
Define the $f$-rewarded line-to-point backward coalescence time
$$
 \tau_{n;(*:f,0)}^{\downarrow;(\{y_1,y_2\},1)} \in [0,t]
$$
from the pair $\{y_1,y_2\}$  at time one (and to time zero) to be the supremum of the vertical coordinates of points in the intersection $\rho_{n;(*:f,0)}^{(y_1,1)} \cap \rho_{n;(*:f,0)}^{(y_2,1)}$.
\end{definition}

Suppose that $f \in \initcond_{\ovbar\coninit}$ for some $\ovbar\coninit \in (0,\infty)^3$. 
Adopt a shorthand that is in force only for this paragraph and the next: set $P = \polyunique_{n;(*:f,0)}^{1} \cap [-1,1]$.
Recall from Lemma~\ref{l.densepolyunique}(2)
that when  $n \in \N$ satisfies
$n > 2^{-3/2} \coninit_1^3 \vee 8 (\coninit_2  + 1)^3$, the set $P$ is a full Lebesgue measure subset of $[-1,1]$.
 
Let $s \in (0,1)$. We identify a subset $R^1_{n;(*:f,0)}(s)$ of $P \times P$
by declaring that an element $(y_1,y_2)$ in this product set belongs to  $R^1_{n;(*:f,0)}(s)$ if and only if 
$\tau_{n;(*:f,0)}^{\downarrow;(\{y_1,y_2\},1)} \geq s$.
\begin{lemma}\label{l.ec}
Suppose that $f \in \initcond_{\ovbar\coninit}$ for some $\ovbar\coninit \in (0,\infty)^3$. Let $n \in \N$. For $y_1,y_2 \in [-1,1]$, write $y_1 R y_2$ for $(y_1,y_2) \in R^1_{n;(*:f,0)}(s)$. Suppose that
$n > 2^{-3/2} \coninit_1^3 \vee 8 (\coninit_2  + 1)^3$.
\begin{enumerate}
\item  $R$ is an equivalence relation on $\polyunique_{n;(*:f,0)}^{1} \cap [-1,1]$. 
\item Let $Q$ denote any equivalence class of  the relation $R$.
There exists an interval $I \subseteq [-1,1]$ such that  
$Q = I \cap \polyunique_{n;(*:f,0)}^{1}$. 
\end{enumerate}
\end{lemma}
\noindent{\em Remark.}
According to the lemma, the interval $I$ could be a singleton, or it might contain at least one of its endpoints. It is in fact not difficult to exclude these circumstances. For example, if $I = \{ a \}$,
then we might consider the limit as $\e \searrow 0$ of the left-neighbouring line-to-point polymers $\rho_{n;(*:f,0)}^{(a - \e,1)}$. 
It can be argued that this limit, which we might denote by $\rho_{n;(*:f,0)}^{(a -,1)}$, is a polymer. Ordinarily, it would equal  $\rho_{n;(*:f,0)}^{(a,1)}$.
Because $I = \{ a \}$, however, the new polymer has a point of intersection with 
$\rho_{n;(*:f,0)}^{(a,1)}$ whose vertical  coordinate exceeds $s$ only at $(a,1)$.
A right limiting polymer  $\rho_{n;(*:f,0)}^{(a +,1)}$
satisfies the same property. Thus, $a$ is not an element of  $\polyunique_{n;(*:f,0)}^{1}$. Since $Q$ is not empty, $I$ cannot then be a singleton.
Similarly, one may find two (but not three) distinct polymers emanating from $(a,1)$ if $I = [a,b)$ with $a < b$.
Although these ideas serve to elucidate the structure of the canopies, we do not need them in this article, so the details of this argument will be omitted.

\noindent{\bf Proof of Lemma~\ref{l.ec}. (1):}
Transitivity of the relation follows because  $\tau_{n;(*:f,0)}^{\downarrow;(\{y_1,y_3\},1)} \geq  \tau_{n;(*:f,0)}^{\downarrow;(\{y_1,y_2\},1)} \wedge  \tau_{n;(*:f,0)}^{\downarrow;(\{y_2,y_3\},1)}$;
the other axioms are trivially verified.

\noindent{\bf (2):} Note first that $\polyunique_{n;(*:f,0)}^{1}$ is dense in $[-1,1]$ by Lemma~\ref{l.densepolyunique}(2). Now let $(y_1,y_2,y_3)$ be an increasing sequence of elements of  $\polyunique_{n;(*:f,0)}^{1} \cap [-1,1]$ such that $y_1 R y_3$. It suffices to argue that $y_1 R y_2$.
Lemma~\ref{l.sandwich} implies that
 $\rho_{n;(*:f,0)}^{(y_1,1)} \preceq \rho_{n;(*:f,0)}^{(y_2,1)}  \preceq \rho_{n;(*:f,0)}^{(y_3,1)}$ 
while the first and third of these polymers are seen to coincide in the strip $\R \times [0,s]$ due to $y_1 R y_3$. Thus, the second polymer also coincides with the other two in this way, so that $y_1 R y_2$.  \qed


\medskip

In the circumstances of Lemma~\ref{l.ec}, the interior in $\R$ of the closure  
of any equivalence class of the relation~$R$ is an open real interval, because  $\polyunique_{n;(*:f,0)}^{1}$ is dense in $[-1,1]$. 
Any such open interval
will be called an $\big(n,*:f,s\big)$-canopy.~\hfff{canopy} 
This definition captures the informal notion presented in  Section~\ref{s.reroot}; more precisely, the formal definition specifies the intersections with $[-1,1]$
of the informally specified objects. 
Whenever $C$ is an  $\big(n,*:f,s\big)$-canopy, the polymer endpoint  $\rho_{n;(*:f),0}^{(x,1)}(0)$ is independent of  $x \in C \cap \polyunique_{n;(*:f,0)}^{1}$.
The shared endpoint will be called the root of the canopy $C$.

\begin{lemma}\label{l.middlehighway}
Suppose that $f \in \initcond_{\ovbar\coninit}$ for some $\ovbar\coninit \in (0,\infty)^3$. 
Let $n \in \N$, $s \in (0,1)$, and let $C$ denote an arbitrary  $\big(n,*:f,s\big)$-canopy. Write $C = (y_1,y_2)$
with $y_1,y_2 \in [-1,1]$ and $y_1 < y_2$, 
and denote by $r$ the root of $C$.
Let $y,z \in (y_1,y_2)$, $y < z$,
be two elements of  $\polyunique_{n;(*:f,0)}^{1}$.
\begin{enumerate}
\item Let $r_- \in \R$ satisfy $r_- < r$,
 $(r_-,y),(r_-,z) \in \polyunique_{n;0}^{1}$ and $\tau_{n;\big(  \{ r_- , r  \} , 0 \big)}^{\uparrow;( y , 1)} \leq s$. Then the polymer $\rho_{n;(r_-,0)}^{(z,1)}$
 is equal to the concatenation of $\rho_{n;(r_-,0)}^{(y,1)} \cap \big( \R \times [0,s]  \big)$
 and  $\rho_{n;(r,0)}^{(z,1)} \cap \big( \R \times  [s,1]  \big)$.
\item Let $r_+ \in \R$ satisfy $r < r_+$,
 $(r_+,y),(r_+,z) \in \polyunique_{n;0}^{1}$ and $\tau_{n;\big(  \{ r, r_+   \} , 0 \big)}^{\uparrow;( z , 1)} \leq s$. Then the polymer $\rho_{n;(r_+,0)}^{(y,1)}$
 is equal to the concatenation of $\rho_{n;(r_+,0)}^{(z,1)} \cap \big( \R \times  [0,s]  \big)$
 and  $\rho_{n;(r,0)}^{(y,1)} \cap \big( \R \times [s,1]  \big)$.
\end{enumerate}
\end{lemma}
\noindent{\bf Proof.} 
The two parts of the lemma are proved similarly and we prove only the first part. To do so, 
note that, since $r_- \leq r$ and $y \leq z$, Lemma~\ref{l.sandwich} implies that 
\begin{equation}\label{e.polysand}
\rho_{n;(r_-,0)}^{(y,1)} \preceq \rho_{n;(r_-,0)}^{(z,1)} \preceq \rho_{n;(r,0)}^{(z,1)} \, .
\end{equation}
Note that 
$$
 \tau_{n; \big( \{ r_- , r \}  , 0 \big)}^{\uparrow;(y,1)}
  \leq  s \leq  
 \tau_{n;(*:f,0)}^{\downarrow; \big( \{ y,z \}  , 1\big)}  = \tau_{n;(r,0)}^{\downarrow; \big( \{ y,z \}  , 1\big)} \, ;
$$
abbreviating  $\tau_- = \tau_{n; \big( \{ r_- , r \}  , 0 \big)}^{\uparrow;(y,1)}$ and  $\tau_+ =  \tau_{n;(r,0)}^{\downarrow; \big( \{ y,z \}  , 1\big)}$, we see that 
$s \in [\tau_-,\tau_+]$.

Note that the polymers $\rho_{n;(r,0)}^{(y,1)}$ and $\rho_{n;(r,0)}^{(z,1)}$ coincide when intersected with the strip $\R \times [0,\tau_+]$,
and that  the polymers
 $\rho_{n;(r_-,0)}^{(y,1)}$ and $\rho_{n;(r,0)}^{(y,1)}$ coincide when intersected with $\R \times [\tau_-,1]$.
Thus, $\rho_{n;(r_-,0)}^{(y,1)}$ coincides with   $\rho_{n;(r,0)}^{(z,1)}$ on $\R \times [\tau_-,\tau_+]$,
which in view of~(\ref{e.polysand}) implies that
\begin{equation}\label{e.middleequality}
\rho_{n;(r_-,0)}^{(y,1)} \cap \big( \R \times [\tau_-,\tau_+]  \big) = \rho_{n;(r_-,0)}^{(z,1)}  \cap \big( \R \times [\tau_-,\tau_+]  \big) =  \rho_{n;(r,0)}^{(z,1)}  \cap \big( \R \times [\tau_-,\tau_+] \big) \, .
\end{equation}
The intersection of $\rho_{n;(r_-,0)}^{(z,1)}$
with $\R \times [\tau_+,1]$ is a polymer whose endpoints 
coincide with the intersection of $\rho_{n;(r,0)}^{(z,1)}$
with this strip. Because there is only one polymer, namely $\rho_{n;(r,0)}^{(z,1)}$, that begins at $(r,0)$ and ends at $(z,1)$, nor can there be more than one polymer whose endpoints are a given pair of points in  $\rho_{n;(r,0)}^{(z,1)}$. Thus, we learn that 
    $\rho_{n;(r_-,0)}^{(z,1)} \cap  \big( \R \times [\tau_+,1] \big) = \rho_{n;(r,0)}^{(z,1)} \cap  \big( \R \times [\tau_+,1]  \big)$. By similar reasoning,  $\rho_{n;(r_-,0)}^{(y,1)} \cap  \big( \R \times [0,\tau_-]  \big) = \rho_{n;(r_-,0)}^{(z,1)} \cap  \big( \R \times [0,\tau_-] \big)$.

We have learnt that $\rho_{n;(r_-,0)}^{(z,1)}$
is the concatenation of the polymers $\rho_{n;(r_-,0)}^{(y,1)} \cap  \big( \R \times [0,\tau_-]  \big)$ and  
$\rho_{n;(r,0)}^{(z,1)}  \cap \big( \R \times [\tau_-,1] \big)$. By~(\ref{e.middleequality}), this assertion also holds when $\tau_-$ is replaced by any element of $[\tau_-,\tau_+]$.
Choosing this element to be $s$ completes the proof of Lemma~\ref{l.middlehighway}(1). \qed

\begin{definition}\label{d.canopyspecialpoint}
Let $\e \in (0,1)$, and let $C$ denote an arbitrary  $\big(n,*:f, 1 - \e^{3/2} \big)$-canopy. Write $C = (x_1,x_2)$
with $x_1,x_2 \in [-1,1]$ and $x_1 < x_2$; denote by $r$ the root of $C$; and set $r_- = \e \lfloor r \e^{-1} \rfloor$.
Define the {\em special point} $\specialpt$ of $C$ to be the infimum of the values of $x \in (x_1,x_2) \cap  \polyunique_{n;(*:f,0)}^{1}$
such that $(r_-,x) \in \polyunique_{n;0}^{1}$ and 
$$
 \tau_{n; \big( \{ r_- , r \}  , 0 \big)}^{\uparrow;(x,1)}
  \leq  1 - \e^{3/2} \, .
$$
We take $\specialpt = x_2$ when no such value of $x$ exists. Thus, $\specialpt$ is a random point that almost surely lies in $[x_1,x_2]$.
\end{definition} 
 \begin{lemma}\label{l.canopytwopieces}
 Suppose that $f \in \initcond_{\ovbar\coninit}$ for some $\ovbar\coninit \in (0,\infty)^3$. Let $K \geq 1$ and let $n \in \N$ satisfy $n >  2^{-3/2} \coninit_1^3 \vee 8 (\coninit_2  + 1)^3 \vee 8 (K+2)^3$. 
Further let $\e \in (0,1)$. Suppose that the event $$
\nolatecoal_{n;([-K,K],0)}^{([-1,1],1)}(\e) \cap \regfluc^{(\{-1,1\},1)}_{n;(*:f,0)}(K-1)
$$ 
occurs.
Let $C$ denote an arbitrary  $\big(n,*:f, 1 - \e^{3/2} \big)$-canopy. Write $C = (x_1,x_2)$
with $x_1,x_2 \in [-1,1]$ and $x_1 < x_2$; denote by $r$ the root of $C$, and set $r_- = \e \lfloor r \e^{-1} \rfloor$ as well as $r_+ = \e \lceil r \e^{-1} \rceil$. Write $\specialpt = \specialpt_C \in [x_1,x_2]$ for the special point of $C$.

The value of
$\weight_{n;(r,0)}^{(x,1)} - \weight_{n;(r_-,0)}^{(x,1)}$ is independent of $x \in (\specialpt,x_2) \cap [-1,1]$, and the value of
$\weight_{n;(r,0)}^{(x,1)} - \weight_{n;(r_+,0)}^{(x,1)}$ is independent of $x \in (x_1,\specialpt) \cap [-1,1]$. 
\end{lemma}
\noindent{\bf Proof.}
We begin by arguing that the  value of
$\weight_{n;(r,0)}^{(z,1)} - \weight_{n;(r_-,0)}^{(z,1)}$ is independent of choices of $z \in (\specialpt,x_2) \subset [-1,1]$
for which $(r^-,z),(r,z) \in \polyunique_{n;0}^{1}$; if it happens that $\specialpt = x_2$, this assertion should be interpreted as being vacuously true.
(Such choices of $z$ will be called $z$-values in the argument that follows. The set of $z$-values is dense 
in~$(\specialpt,x_2)$ by Lemma~\ref{l.densepolyunique}(1) and~(2), where Lemma~\ref{l.densepolyunique}(1) is applicable because $n \geq 8(K+2)^3$ and $\vert r^- - z \vert \leq K+2$.)
To establish the claimed constancy, note that, for any such $z$, there exists, by the definition of $\specialpt$, an element $y \in (\specialpt,z)$ such that  $(r,y), (r^-,y) \in   \polyunique_{n;0}^{1}$
and 
$\tau_{n; \big( \{ r_- , r \}  , 0 \big)}^{\uparrow;(y,1)}
  \leq  1 - \e^{3/2}$. Applying Lemma~\ref{l.middlehighway}(1) with ${\bf s} = 1 - \e^{3/2}$, we find that  the polymer $\rho_{n;(r_-,0)}^{(z,1)}$
 is equal to the concatenation of $\rho_{n;(r_-,0)}^{(y,1)} \cap \big( \R \times [0,1 - \e^{3/2}] \big)$
 and  $\rho_{n;(r,0)}^{(z,1)} \cap \big( \R \times [1 - \e^{3/2},1]  \big)$.
 The polymer  $\rho_{n;(r,0)}^{(z,1)}$ may of course also be viewed as a concatenation of its subpaths in these two strips. The difference in weight 
$\weight_{n;(r,0)}^{(z,1)} - \weight_{n;(r_-,0)}^{(z,1)}$ between these two polymers is thus seen to equal the difference in weight between  $\rho_{n;(r_-,0)}^{(y,1)} \cap \big( \R \times [0,1 - \e^{3/2}]  \big)$ and  $\rho_{n;(r,0)}^{(y,1)} \cap \big( \R \times [0,1 - \e^{3/2}]  \big)$. The latter quantity appears to be independent of $z$, although it does depend on $y$, which was determined by $z$. 
However, we may introduce $\e > 0$, and insist that $y \in (\specialpt,\specialpt + \e)$; since a viable $y$ may be found for any such $\e$, the weight difference in question is seen to be shared among $y$ that exceed $\specialpt$ by an arbitrarily small amount. We see then that the weight difference is indeed independent of  the $z$-value $z$, provided that $z \geq \specialpt + \e$.
Since $\e > 0$ is arbitrary, the quantity is in fact independent of the $z$-value without restriction. 

That the quantity
$\weight_{n;(r,0)}^{(x,1)} - \weight_{n;(r_-,0)}^{(x,1)}$ is independent of $x \in (\specialpt,x_2)$ now follows, because the set of $z$ we have been considering is dense in  $(\specialpt,x_2)$, and this quantity is a continuous function of $x$ in this interval by Lemma~\ref{l.basic}.
  
  It remains to argue that 
$\weight_{n;(r,0)}^{(x,1)} - \weight_{n;(r_+,0)}^{(x,1)}$ is independent of $x \in (x_1,\specialpt)  \subset [-1,1]$. 
To establish this, we first prove that 
\begin{equation}\label{e.secondzvalue}
\weight_{n;(r,0)}^{(z,1)} - \weight_{n;(r_+,0)}^{(z,1)} \, \, \, \textrm{is independent of choices of} \, \, \, z \in (x_1,\specialpt)
\end{equation}
where $z \in [-1,1]$ is such that $(r^-,z)$, $(r,z)$ and $(r^+,z)$ are elements of $\polyunique_{n;0}^{1}$. (In a reuse of terminology, any such $z$ will be called a $z$-value in the ensuing argument. Again, it is Lemma~\ref{l.densepolyunique}(1) and~(2) that assure that the set of $z$-values is dense.) 
To verify this assertion, note that, since each such $z$ is less than $\specialpt$,  
$\tau_{n; \big( \{ r_- , r \}  , 0 \big)}^{\uparrow;(z,1)} >  1 - \e^{3/2}$. Any $z$-value $z$ is an $(n,\e,f)$-triple uniqueness point in the sense specified before 
Lemma~\ref{l.lcrf}; thus, this lemma implies that one or other of 
$\tau_{n; \big( \{ r_- , r \}  , 0 \big)}^{\uparrow;(z,1)}$ and 
$\tau_{n; \big( \{ r , r_+  \}  , 0 \big)}^{\uparrow;(z,1)}$ is at most $1 - \e^{3/2}$. 
We learn then that 
$\tau_{n; \big( \{ r , r_+  \}  , 0 \big)}^{\uparrow;(z,1)} \leq 1 - \e^{3/2}$. This inference permits us consider any pair $(z_1,z_2)$, $z_1 < z_2$, of $z$-values, and to turn to Lemma~\ref{l.middlehighway}(2)
with ${\bf y} = z_1$, ${\bf z} = z_2$ and ${\bf s} = 1 - \e^{3/2}$. From this lemma, we find that the polymer weight $\weight_{n;(r_+,0)}^{(z_1,1)}$
is the sum of the weights of the polymers $\rho_{n;(r_+,0)}^{(z_2,1)} \cap \big(    \R \times [0,1-\e^{3/2}]  \big)$
 and  $\rho_{n;(r,0)}^{(z_1,1)} \cap \big( \R \times [1-\e^{3/2},1]  \big)$.
The weight difference  $\weight_{n;(r_+,0)}^{(z_1,1)} -  \weight_{n;(r,0)}^{(z_1,1)}$
is thus seen to equal the difference of the weights of the polymers  $\rho_{n;(r_+,0)}^{(z_2,1)} \cap \big( \R \times [0,1-\e^{3/2}]  \big)$ and  $\rho_{n;(r,0)}^{(z_2,1)} \cap \big(  \R \times [0,1-\e^{3/2}]  \big)$. This quantity is independent of $z_1$, though it does depend on $z_2$. However, $z_2$ may be chosen to  
 be a $z$-value arbitrarily close to $\specialpt$, by the density of such values, 
 rendering the quantity independent of $z_1$ to the left of any neighbourhood of $\specialpt$.
 Since the neighbourhood is arbitrary, we confirm~(\ref{e.secondzvalue}).

That 
$\weight_{n;(r,0)}^{(x,1)} - \weight_{n;(r_+,0)}^{(x,1)}$ is independent of $x \in (x_1,\specialpt)$ 
follows from the density of $z$-values in $(x_1,\specialpt)$ and the polymer weight continuity Lemma~\ref{l.basic}. This completes the proof of Lemma~\ref{l.canopytwopieces}. \qed

\section{Well-behaved canopy structures are typical}

We are almost ready to rigorously construct the $f$-rewarded polymer weight profiles as Brownian patchwork quilts (that is, to prove Theorem~\ref{t.unifpatchcompare}). This section is devoted to some final preparations for this task.
Recall that it is settled that  the patches that will make up the Brownian quilts will be the intervals obtained by splitting $(n,*:f,1-\e^{3/2})$-canopies at the special points inside these canopies; (the split pieces will be called split canopies when they are formally introduced in the context of the upcoming proof, in Section~\ref{s.quilt}). 
Recall also that we will actually select the value of $\e > 0$ to be random, by searching through consecutively smaller dyadic scales until a certain favourable event is realized. This event will include the absence of late coalescence, but that is not the only desired feature. For example, we will also insist that the number of  $(n,*:f,1-\e^{3/2})$-canopies be not too high, because this cardinality is in essence the size of the random stitch point set in the constructed quilt, which we wish to argue is not typically too big. 

In this section, then, we offer a definition of a {\em normal coalescence} event, which will play the role of this favourable event. Lemma~\ref{l.normalcoal} is the only result in the section. It establishes that the normal coalescence event fails with the same order of probability, $\e^{2 + o(1)}$, that the event of absence of late coalescence does. This is because the other favourable features that we build into the normal coalescence event, such as the   $(n,*:f,1-\e^{3/2})$-canopy cardinality not being abnormally high, are highly typical scenarios.

Let $\ovbar\coninit \in (0,\infty)^3$ and  $f \in \initcond_{\ovbar\coninit}$. For $n \in \N$ and $s \in (0,1)$, 
let $\canopynumber_{n,;(*:f,0)}^{([-1,1],1)}(s)$ denote the cardinality of the set of $(n;*:f,s)$-canopies that intersect $[-1,1]$.~\hfff{canopynumber}

For $n \in \N$, $\e \in (0,1)$, $\chimac > 0$ and $D > 0$, define the {\em normal coalescence} event $\normalcoal_{n;(*:f,0)}^{([-1,1],1)}\big(D,\e,\chimac \big)$ to be~\hfff{normalcoal} 
\begin{eqnarray*}
 & & \nolatecoal_{n;\big([-D (\log \e^{-1})^{1/3},D (\log \e^{-1})^{1/3}],0\big)}^{([-1,1],1)}(\e) \, \cap \, 
\regfluc_{n;(*:f,0)}^{\big( \{-1,1 \}, 1 \big)}\Big( D (\log \e^{-1})^{1/3} -1 \Big) \\ 
  & & \qquad \qquad \qquad \qquad \qquad  \cap \, \, \Big\{ \canopynumber_{n,;(*:f,0)}^{([-1,1],1)}\big( 1  - \e^{3/2} \big) \leq   \e^{-1 - \chimac} \Big\} \, .
\end{eqnarray*}

The unappetizing aspect of explicit hypothesis bounds is perhaps most conspicuous in the next result (and correspondingly, the calculational derivation of the result in Appendix~\ref{s.calcder} is rather long): 
the basic meaning is that, with $D$ and $\chi$ fixed to be large and small, $\normalcoal$ fails with probability at most $\e^{2 + o(1)}$, uniformly in high $n$. The remark following the lemma's statement provides a less forbidding summary of these bounds.

\begin{lemma}\label{l.normalcoal}
Suppose that $n \in \N$, $\e > 0$,  $D > 0$, $\chimac > 0$ and $\ovbar\coninit \in (0,\infty)^3$ satisfy the bounds 
\begin{eqnarray*}
  \e  & \leq & \max \Bigg\{ \,
(\eta_0)^{72} \, , \,   \big( 10^{-41} \chi^{36} D^{-14} \big)^{\chi^{-2}} \, , \,   \exp \big\{ - \beta^{ \chimac^{-1/2}  (3480)^{1/2}  } \big\} \, , \,  \\
 &  & \qquad \qquad (2 \condee)^{- 3480 \chimac^{-1} ( \log \beta)^{2}  }
 \, , \,  \big( 4k_0 \big)^{-\chi^{-1}} \, ,    \,  \exp \big\{ - 2 K_0^2 \big\} \, ,   \, 10^{-2634} c_3^{88} \, ,  \\
  & & \qquad \qquad \qquad
     \exp \bigg\{ - D^{-3} \Big( 78 \Psi_1 \vee 4 \big( (\Psi_2 + 1)^2 + \Psi_3 \big)^{1/2}  \Big)^{3} \bigg\} \, , \,  \exp \Big\{ -  10^7 c_1^{-9} C^{3/8} \Big\} \,  \Bigg\} \, ,
\end{eqnarray*}
$$
    n  \geq  \max \bigg\{ \,   
 c^{-18} \max \Big\{  (\coninit_2 + 1)^9 \, , \,  10^{23} \coninit_1^9 \Big\}  \, , \, 4(K_0)^9 \big( \log \e^{-1} \big)^{K_0} \, ,  \,   10^{740}    D^{36}  c_3^{-84}  
 a_0^{-9}    \e^{-504}  \, \bigg\} \, ,
 $$
$n \e^{3/2} \in \N$,
$D  \geq 10^{16} c_1^{-4/3}$ and  
$\chimac   \leq  2^{-1}   \big( 1  + 500( \log \beta)^{2} \big)^{-1}$. 
Let $f \in \initcond_{\ovbar\coninit}$. Then
 $$
 \PP \Big( \neg \, \normalcoal_{n;(*:f,0)}^{([-1,1],1)}\big(D,\e,\chimac \big) \Big) \, \leq 
 \, \e^2 \cdot \Omega(\e)  \, ,
 $$
 where 
 \begin{equation}\label{e.omegaform}
 \Omega(\e) = 10^{421}    c_3^{-33}      C_3  D  \big( \log \e^{-1} \big)^{43}     \exp \big\{ \beta_3 \big( \log \e^{-1} \big)^{5/6} \big\}  \, .
 \end{equation}
 The positive constants  $k_0$, $\eta_0$, $K_0$, $a_0$, $\beta_3$ and $\beta$ are each specified either in Theorem~\ref{t.maxpoly} or Theorem~\ref{t.disjtpoly}.
  By $\condee$, we denote the finite supremum $\sup_{i \in \N} \conseqmac_i \exp \big\{ - 2 (\log i)^{11/12} \big\}$  associated to the sequence  $\big\{ \conseqmac_i: i \in \N \big\}$ of constants provided by Theorem~\ref{t.maxpoly}.
\end{lemma}
\noindent{\em Remark.}
An inexplicit presentation of the lemma's hypotheses takes the form: that $D$ exceeds a large positive constant, and that $\chimac >0$ is at most a small one; that 
$\e$ is at most a small positive constant determined by 
$D$ and $\chimac$ as well as by the function space parameter 
$\ovbar\coninit$; and that $n$ is at least $\Theta(\chi,D,\ovbar\coninit) \e^{-\Theta}$, where $\Theta$ in the exponent denotes a large universal constant, and the dependence of the other large constant $\Theta$ is indicated. The reader of the ensuing proof who prefers not to examine the separated calculational derivation may readily verify that these conditions ensure the validity of every one of the needed bounds that will be explicitly recorded during the proof.  

\noindent{\bf Proof of Lemma~\ref{l.normalcoal}.} 
In this analysis, we will consider the event 
\begin{eqnarray*}
 & &    \bigcap_{y \in \e\Z \cap [-1,1]} \pdr_{n; \big( -D (\log \e^{-1})^{1/3} ,0 \big)}^{(y,1)} \Big( 1 - \e^{3/2} ,  D \big(\log \e^{-1} \big)^{4/3} \Big) \\
 & & \cap \,  \bigcap_{y \in \e\Z \cap [-1,1]}  \pdr_{n; \big( D (\log \e^{-1})^{1/3} ,0 \big)}^{(y,1)}\Big( 1 - \e^{3/2} ,  D \big(\log \e^{-1}\big)^{4/3}  \Big) \, , 
\end{eqnarray*}
and will use the shorthand $\cap \, \pdr$
to refer to it in this paragraph and the next.
The event  $\normalcoal_{n;(*:f,0)}^{([-1,1],1)}\big(D,\e , \chimac \big)$ is an intersection of a $\nolatecoal$ event, a $\regfluc$ event and a $\canopynumber$ event. Let $\manycan$ denote the intersection of $\cap \, \pdr$, this $\regfluc$ event  and the complement of the $\canopynumber$ event. Although the new event entails several circumstances, each of these is typical, except for the complement of the $\canopynumber$ event: something reflected in the  name $\manycan$.

The lemma will be proved by noting that
\begin{equation}\label{e.fourterms}
 \PP \big( \neg \, \normalcoal \big) \leq \PP \big( \latecoal  \big) +  \PP \big( \manycan  \big) +  \PP \big( \neg \, \regfluc  \big) +  \PP \big( \neg \, \cap \pdr \big) \, ,
\end{equation} 
where we omit all adornments of the events specifying $\normalcoal$.
It is the first term on the right-hand side that dictates the bound in Lemma~\ref{l.normalcoal}: it is of order $\e^{2 + o(1)}$.
The further three terms will be shown to be at least as small. 

We now find upper bounds on the four terms on the right-hand side of~(\ref{e.fourterms}). The order in which we do so is $2134$,  
the analysis of the second term,  $\PP \big( \manycan  \big)$, being slightly more involved.

When $\manycan$ occurs, the bound $\canopynumber_{n,;(*:f,0)}^{([-1,1],1)}\big( 1  - \e^{3/2} \big) >   \e^{-1 - \chimac}$ holds, and
we may select from each $(n;(*:f),1-\e^{3/2})$-canopy that intersects $[-1,1]$ an element of $y \in \polyunique_{n;(*:f,0)}^1 \cap [-1,1]$.
Since such elements~$y$ number at least $\lceil \e^{-1 - \chimac} \rceil$, we may discard certain among them and list those remaining   in the form $\big\{ y_i : i \in \intint{ \lceil \e^{-1 - \chimac} \rceil } \big\}$. The polymers $\rho_{n;(*:f,0)}^{(y_i,1)}$ for   $i \in \intint{ \lceil \e^{-1 - \chimac} \rceil }$ are then pairwise disjoint when intersected with the strip $\R \times [1-\e^{3/2},1]$. By the occurrence of the $\regfluc$ event and the polymer sandwiching Lemma~\ref{l.sandwich},  $\rho_{n;(*:f,0)}^{(y_i,1)}(0) \in \big[ - D (\log \e^{-1})^{1/3},D (\log \e^{-1})^{1/3} \big]$ for all such indices~$i$. 
For this reason, the intersection of $\rho_{n;(*:f,0)}^{(y_i,1)}$
with the strip $\R \times [0,1-\e^{3/2}]$ lies on, or to the right of,
  $\rho_{n;\big(   - D (\log \e^{-1})^{1/3} , 0 \big)}^{(\e \lfloor \e^{-1} y_i \rfloor ,1)}$. (In specifying the latter polymer, note that $\e \lfloor \e^{-1} y_i \rfloor$ is the left-displacment of $y_i$ onto an $\e$-mesh. This polymer is almost surely well defined, by Lemma~\ref{l.densepolyunique}(1) with 
   ${\bf x} =  - D (\log \e^{-1})^{1/3}$ and ${\bf y} = \e \lfloor \e^{-1} y_i \rfloor$.)
We may thus bound
\begin{eqnarray*}
\rho_{n;(*:f,0)}^{(y_i,1)}\big( 1 - \e^{3/2} \big) & \geq & \rho_{n;\big(   - D (\log \e^{-1})^{1/3} , 0 \big)}^{(\e \lfloor \e^{-1} y_i \rfloor ,1)}\big( 1 - \e^{3/2} \big) \\
 & \geq & \big( 1 - \e^{3/2} \big) \e \lfloor \e^{-1} y_i \rfloor  - \e^{3/2}  D (\log \e^{-1})^{1/3} \, - \, \e  D \big(\log \e^{-1} \big)^{4/3} \\
 & \geq & y_i - \e^{3/2} y_i - \e   \, - \, 2\e  D \big(\log \e^{-1} \big)^{4/3} \geq y_i   \, - \, 4\e  D \big(\log \e^{-1} \big)^{4/3} \, ,
\end{eqnarray*}
where the second bound is due to the occurrence of $ \pdr_{n; \big( -D (\log \e^{-1})^{1/3} ,0 \big)}^{(y,1)} \big( 1 - \e^{3/2} ,  D (\log \e^{-1} )^{4/3} \big)$ and $\e^{3/2} \leq 2^{-1}$, the third to $\e \lfloor \e^{-1} y_i \rfloor \geq y_i - \e$ and $\e \leq e^{-1}$, and the fourth to $y_i \in [-1,1]$, $D \geq 1$ and $\e \leq e^{-1}$.
A similar upper bound on 
$\rho_{n;(*:f,0)}^{(y_i,1)}\big( 1 - \e^{3/2} \big)$ holds, and we find that
\begin{equation}\label{e.rhoabs}
\Big\vert \,
\rho_{n;(*:f,0)}^{(y_i,1)}\big( 1 - \e^{3/2} \big) - y_i \, \Big\vert \, \leq \, \e  \cdot 4D \big(\log \e^{-1} \big)^{4/3} \, .
\end{equation}
The polymers $\rho_{n;(*:f,0)}^{(y_i,1)}$ number at least $\e^{-1 - \chimac}$. They may be divided according to which of the intervals $\big[j\e,(j+1)\e \big)$, $j \in \llbracket 0, \lfloor \e^{-1} \rfloor \rrbracket$,
contains $y_i$. Thus, there exists $J  \in \llbracket 0, \lfloor \e^{-1} \rfloor \rrbracket$ for which $\big[J\e,(J+1)\e \big)$ contains at least $\e^{-1 - \chimac} ( \e^{-1} + 1)^{-1} \geq \e^{-\chimac}/2$ (where we used $\e \leq 1$).
Set
$$
I_1 = \big[ J\e - 4D \e \big(\log \e^{-1} \big)^{4/3}, (J+1)\e +  4D \e \big(\log \e^{-1} \big)^{4/3}  \big] \, \, \textrm{ and } \, \, I_2 = \big[J\e,(J+1)\e \big] \, .
$$
We see from~(\ref{e.rhoabs}) that, for indices $i$ such that $y_i \in  \big[J\e,(J+1)\e \big)$,  the polymers $\rho_{n;(*:f,0)}^{(y_i,1)}$ restricted to the strip $\R \times [1-\e^{3/2},1]$
are disjoint polymers that begin in $I_1 \times \{ 1-\e^{3/2} \}$ and end in $I_2 \times \{ 1 \}$.
Thus, we find that
$\maxpoly_{n;\big(  I_1  ,  1 - \e^{3/2} \big)}^{(I_2,1)} \geq \e^{-\chimac}/2$
occurs. That is,  $\manycan$ is a subset of
$$
 \bigcup_{j \in \llbracket 0, \lfloor \e^{-1} \rfloor \rrbracket}   \bigg\{ \, \maxpoly_{n;\big(  \big[ j\e - 4D \e \big(\log \e^{-1} \big)^{4/3}, (j+1)\e +  4D \e \big(\log \e^{-1} \big)^{4/3}  \big]  ,  1 - \e^{3/2} \big)}^{\big( [ j\e,(j+1)\e ],1\big)} \geq 2^{-1} \e^{-\chimac} \, \bigg\} \, .
$$

For given $j \in \llbracket 0, \lfloor \e^{-1} \lfloor \rrbracket$, the probability of the $\maxpoly$ event on display may be bounded above by an application of Theorem~\ref{t.maxpoly}.
The theorem's parameters are set: ${\bf t_1} = 1 - \e^{3/2}$, ${\bf t_2} = 1$, 
${\bf x} =  j\e - 4D \e \big(\log \e^{-1} \big)^{4/3}$, ${\bf y} = j \e$, ${\bf a} = \lceil 8 D \big(\log \e^{-1} \big)^{4/3} \rceil$, ${\bf b} = 1$ and  
${\bf k} = \lfloor 2^{-1} \e^{-\chimac} \rfloor$. With $h =  \lceil 8 D \big(\log \e^{-1} \big)^{4/3} \rceil$, we find that the probability in question is at most
\begin{eqnarray*}
& &
  \big(   4^{-1} \e^{-\chimac} \big)^{-   (145)^{-1} ( \log \beta)^{-2} (0 \vee \log \log \e^{-\chimac})^2}  \cdot h^{(\log \beta)^{-2} (\log \log \e^{-\chimac})^2/{288}  + 3/2} \conseqmac_{\lfloor 2^{-1} \e^{-\chimac} \rfloor} 
  \\
  & \leq  &
  \big(   4^{-1} \e^{-\chimac} \big)^{-   4^{-1} (145)^{-1} ( \log \beta)^{-2} (\log \log \e^{-1})^2}  \cdot h^{(\log \beta)^{-2} (\log \log \e^{-1})^2/{288}  + 3/2}  \condee \exp \big\{  2 (\log  \e^{-1})^{11/12} \big\} \\
  & \leq  &
   \e^{ \chimac  (580)^{-1} ( \log \beta)^{-2} (\log \log \e^{-1})^2}  \cdot   \condee  (2h)^{(\log \beta)^{-2} (\log \log \e^{-1})^2/{288}  + 3/2} \exp \big\{  2 (\log  \e^{-1})^{11/12} \big\} \, ,
\end{eqnarray*}
where recall that  we denote by $\condee$ the finite supremum $\sup_{i \in \N} \conseqmac_i \exp \big\{ - 2 (\log i)^{11/12} \big\}$  associated to the sequence  $\big\{ \conseqmac_i: i \in \N \big\}$ of constants provided by Theorem~\ref{t.maxpoly}.
The form of the first expression is obtained by using $\e \leq 4^{-1/\chimac}$ in the form ${\bf k} \geq 4^{-1} \e^{-\chimac}$.
The first displayed inequality made use of
$0 \vee \log \log \e^{-\chimac} \geq 2^{-1} \log \log \e^{-1}$, which is implied by 
$\log \chi^{-1} \leq  2^{-1} \log \log \e^{-1}$ or equivalently $\e \leq \exp \{ - \chi^{-2} \}$; this inequality is also due to $\chi \leq 1$.

For this application of Theorem~\ref{t.maxpoly}
to be made, since 
$4^{-1} \e^{-\chimac} \leq {\bf k} \leq 2^{-1} \e^{-\chimac}$, it suffices that, writing $\tau =   4D \big(\log \e^{-1} \big)^{4/3}   + 2 \lceil 8 D \big(\log \e^{-1} \big)^{4/3} \rceil$, we have that 
 $$
 4^{-1} \e^{-\chimac} \geq  k_0 \vee \big( 4D  \big(\log \e^{-1} \big)^{4/3} + 2 \lceil 8 D \big(\log \e^{-1} \big)^{4/3} \rceil \big)^3  \, ,
 $$ 
and
\begin{eqnarray*}
 n \e^{3/2}   & \geq  &  \max \bigg\{ \,   2(K_0)^{(12)^{-2} \big( \log \log   (2^{-1} \e^{-\chimac})     \big)^2}   \big( \log ( 2^{-1} \e^{-\chimac}) \big)^{K_0}
    \, , \,  a_0^{-9} \tau^9 \, , \\
    & & \qquad \qquad \qquad     
    10^{325}       \rsc^{-36}   \big( 2^{-1} \e^{-\chimac} \big)^{465}  \max \big\{   1  \, , \,  \tau^{36}   \big\}  
   \,  \bigg\} \, .
\end{eqnarray*}

Taking a union bound over the $\lfloor \e^{-1} \rfloor + 1 \leq 2\e^{-1}$ choices of $j$, we find that $\PP \big( \manycan \big)$ is at most 
\begin{equation}\label{e.ub.manycan}
 \e^{-1 +  \chimac  (580)^{-1} ( \log \beta)^{-2} (\log \log \e^{-1})^2}  \cdot  2 \condee  (2h)^{(\log \beta)^{-2} (\log \log \e^{-1})^2/{288}  + 3/2} \exp \big\{  2 (\log  \e^{-1})^{11/12} \big\} \, .
\end{equation}

Applying Proposition~\ref{p.latecoal}
with ${\bf K} = D (\log \e^{-1})^{1/3}$ and ${\bm \e} = \e$, we find that
\begin{eqnarray}
  & & \PP \Big( \latecoal_{n;\big([-D (\log \e^{-1})^{1/3},D (\log \e^{-1})^{1/3}],0\big)}^{([-1,1],1)}(\e) \Big) \nonumber \\
  & \leq &  \e^2 \cdot 10^{420}    c_3^{-33}      C_3  D  \big( \log \e^{-1} \big)^{127/3}     \exp \big\{ \beta_3 \big( \log \e^{-1} \big)^{5/6} \big\}\, ; \label{e.ub.latecoal}
\end{eqnarray}
since ${\bf K} \geq 2$ if and only if $\e \leq \exp \big\{ - 8D^{-3} \big\}$,
it is sufficient for the application to be made that $n \e^{3/2} \in \N$,
$$
\e  \leq 
\min \Big\{ \exp \big\{ - 8D^{-3} \big\}  \, , \,      (\eta_0)^{72} \, , \, 10^{- 1342}  c_3^{44} \big( D (\log \e^{-1})^{1/3} + 2 \big)^{-8} \,  , \, \exp \big\{ - 2 C^{3/8} \big\} \,
   \Big\} 
$$
and $n \in \N$ verifies
$$
     n  \geq 2 \max \bigg\{ \, 2(K_0)^9 \big( \log \e^{-1} \big)^{K_0}
    \, , \,   
 10^{728}     c_3^{-84}  \e^{-222}   \big( D (\log \e^{-1})^{1/3}  +  2 \big)^{36}  \, , \, 
  a_0^{-9} \big( D (\log \e^{-1})^{1/3}  +   2 \big)^9 2^{6}  \,  \Bigg\} \, .
$$

Applying Lemma~\ref{l.regfluc} with ${\bf R} = D (\log \e^{-1})^{1/3} -1$, and using $\big( 2^{ - 1/2} - 1/2 \big)^{3/2} \geq 2^{-4}$, we learn that
\begin{equation}\label{e.ub.regfluc}
\PP \bigg( \neg \,  \regfluc_{n;(*:f,0)}^{\big( \{-1,1 \}, 1 \big)}\Big( D (\log \e^{-1})^{1/3} -1 \Big) \bigg)
\leq 38 D (\log \e^{-1})^{1/3} \rsC       \e^{   2^{-13} \rsc   D^3} 
\end{equation}
since ${\bf R} \geq 2^{-1} D (\log \e^{-1})^{1/3}$
due to $\e \leq \exp \big\{ - 8D^{-3} \big\}$.
This application is valid when
$$
n \geq 
 c^{-18} \max \Big\{  (\coninit_2 + 1)^9 \, , \,   10^{23} \coninit_1^9     \, , \, 3^{9}  \Big\} 
$$ 
and
$$
 D (\log \e^{-1})^{1/3} - 1  \in \Big[  \,  39 \coninit_1  \,  \vee \, 5  \, \vee \,   3 c^{-3}  \, \vee \, 2 \big( (\coninit_2 + 1 )^2 +  \coninit_3 \big)^{1/2}  \, , \,   6^{-1} c n^{1/9} \,  \Big] \, .
$$

Let $y \in [-1,1]$ be given. Applying Proposition~\ref{p.polyfluc} with ${\bf t_1} = 0$, ${\bf t_2} = 1$, ${\bf x} =  -D (\log \e^{-1})^{1/3}$, ${\bf y} = y$, ${\bf a} = 1 -\e^{3/2}$ and ${\bf r} =  D \big(\log \e^{-1} \big)^{4/3}$, we find that
$$
 \PP \bigg( \neg \, \pdr_{n; \big( -D (\log \e^{-1})^{1/3} ,0 \big)}^{(y,1)} \Big( 1 - \e^{3/2} ,  D \big(\log \e^{-1} \big)^{4/3} \Big) \bigg) \leq 
22 \,  C D \big(\log \e^{-1} \big)^{4/3}    \e^{  10^{-11} c_1  D^{3/4}  }  \, .
$$
Since $y \in [-1,1]$, the application may be made provided that  $\e^{3/2} \leq  10^{-11} c_1^2$; the integer $n$ satisfies
$$
 n    \geq  \max \bigg\{ 
 10^{32} \e^{-75/2} c^{-18} \, \, , \, \, 
 10^{24} c^{-18} \e^{-75/2} \big( D (\log \e^{-1})^{1/3} + 1 \big)^{36}  
 \bigg\} \, ;
$$
and
$$ 
  D \big(\log \e^{-1} \big)^{4/3} \in \bigg[
 \,  10^9 c_1^{-4/5} \, \vee  \,  15 C^{1/2} \, \vee  \,  87 \e^{1/2}   \big( D (\log \e^{-1})^{1/3} + 1 \big) \, \, , \, \, 
 3  \e^{25/6} n^{1/36}  \, \bigg] \, . 
 $$

The same upper bound may be found on the quantity 
$$
\PP \bigg( \neg \, \pdr_{n; \big( D (\log \e^{-1})^{1/3} ,0 \big)}^{(y,1)} \Big( 1 - \e^{3/2} ,  D \big(\log \e^{-1} \big)^{4/3} \Big) \bigg)
$$
by this application, with  ${\bf x}$ instead set equal to $D (\log \e^{-1})^{1/3}$; the conditions on parameters that permit the application are unchanged.
Taking a union bound over the $\lfloor 2 \e^{-1} \rfloor + 1 \leq 4\e^{-1}$ values of $y \in \e\Z \cap [-1,1]$,
we find that the failure probability for the intersection of $\pdr$
events in the definition of the $\normalcoal$ event is at most
\begin{equation}\label{e.ub.pdr}
8 \cdot 22 \,  C D \big(\log \e^{-1} \big)^{4/3}    \e^{  10^{-11} c_1  D^{3/4}  - 1} \, .
\end{equation}

We now assemble the estimates.
The four terms on the right-hand side of~(\ref{e.fourterms})
are bounded above by~(\ref{e.ub.latecoal}),~(\ref{e.ub.manycan}),~(\ref{e.ub.regfluc}) and~(\ref{e.ub.pdr}). That is,
\begin{eqnarray*}
 & &  \PP \Big( \neg \, \normalcoal_{n;(*:f,0)}^{([-1,1],1)}\big(D,\e,\chimac\big) \Big) \\
 & \leq &    \e^2 \cdot 10^{420}    c_3^{-33}      C_3  D  \big( \log \e^{-1} \big)^{127/3}     \exp \big\{ \beta_3 \big( \log \e^{-1} \big)^{5/6} \big\} \\
 & & \, + \,
 \e^{-1 +  \chimac  (580)^{-1} ( \log \beta)^{-2} (\log \log \e^{-1})^2} \\
 & & \qquad \qquad \qquad \cdot \, 
   2 \condee  (2h)^{(\log \beta)^{-2} (\log \log \e^{-1})^2/{288}  + 3/2} \exp \big\{  2 (\log  \e^{-1})^{11/12} \big\}  \\
& & \, + \,  38 D (\log \e^{-1})^{1/3} \rsC       \e^{   2^{-13} \rsc   D^3}   \, + \, 
8 \cdot 22 \,  C D \big(\log \e^{-1} \big)^{4/3}    \e^{  10^{-11} c_1  D^{3/4}  - 1}  \, .
\end{eqnarray*}
The second, third and fourth summands are all at most $\e^2$ provided that
\begin{eqnarray*}
 & & \e^{-3 +  \chimac  (580)^{-1} ( \log \beta)^{-2} (\log \log \e^{-1})^2} \\
  & \leq & 
   2^{-1} \condee^{-1}  \big( 2 \lceil 8 D \big(\log \e^{-1} \big)^{4/3}\rceil   \big)^{-(\log \beta)^{-2} (\log \log \e^{-1})^2/{288}  - 3/2} \exp \big\{-  2 (\log  \e^{-1})^{11/12} \big\} \, ,
\end{eqnarray*}
 $$
      \e^{   2^{-13} \rsc   D^3 - 2} \leq   (38)^{-1} D^{-1} (\log \e^{-1})^{-1/3} \rsC^{-1} 
 $$
 and
$$
  \e^{  10^{-11} c_1  D^{3/4}  - 3} \leq 
8^{-1} 22^{-1} \, c_1 C^{-1} D^{-1} \big(\log \e^{-1} \big)^{-4/3} \, .  
$$ 
Confirming these three bounds is a calculational matter which is undertaken  alongside the calculational derivation of Lemma~\ref{l.normalcoal} in Appendix~\ref{s.calcder}.
Briefly, however,
the first of these conditions is ensured by our insistence that $\e > 0$ be small enough; and
the second and third hold because we insist that $D$ exceed a certain positive lower bound, and then that $\e > 0$ be smaller than a positive constant determined by that lower bound.  

Since $c_3 \leq 1$, $C_3 \wedge D \geq 1$, $\beta_3 \geq 0$ and $\e \leq e^{-1}$, we conclude that
$$
 \PP \Big( \neg \, \normalcoal_{n;(*:f,0)}^{([-1,1],1)}\big(D,\e,\chimac\big) \Big) \leq    \e^2 \cdot 10^{421}    c_3^{-33}      C_3  D  \big( \log \e^{-1} \big)^{43}     \exp \big\{ \beta_3 \big( \log \e^{-1} \big)^{5/6} \big\} \, .
$$
This completes the proof of Lemma~\ref{l.normalcoal}. \qed

\section{Uncovering the Brownian patchwork quilt: the derivation of the main result}\label{s.quilt}

We have built the tools needed to derive Theorem~\ref{t.unifpatchcompare}, and this section is devoted to proving this result. Recall that it is our job 
to demonstrate that
the weight profiles $[-1,1] \to \R: y \to \weight_{n;(*:f,0)}^{(y,1)}$, indexed by 
$(n,f) \in \N \times  \initcond_{\ovbar\coninit}$,
are  uniformly Brownian patchwork $\big(2,3,1/252\big)$-quiltable.
The definition of quiltable entails that a certain cast of characters must be introduced, and proved to enjoy certain properties and relationships. 
It may be helpful for the reader's  bearings that we recall these things now for our particular setting.

Our index set $\indexset$ is $\initcond_{\ovbar\coninit}$; its generic element will be called $f$.
The random function $X_{n,f}$, indexed by $(n,f) \in \N \times \initcond_{\ovbar\coninit}$, which we have been endeavouring to depict as a Brownian patchwork quilt, is, naturally, $[-1,1] \to \R: y \to \weight_{n;(*:f,0)}^{(y,1)}$.
The characters whose presence is called for are 
\begin{itemize}
\item the two sequences $p$ and $q$, which should verify $p_j \leq j^{-2 + o(1)}$ and $q_j \leq j^{-1/252 + o(1)}$;
\item the error event $E_{n,f}$, whose $\PP$-probability must be at most $q_n$;
\item the fabric sequence elements $F_{n,f;i}$, indexed by $(n,f,i) \in \N \times \initcond_{\ovbar\coninit} \times \N$, which will be expected to uniformly withstand $L^{3-}$-comparison to Brownian bridge above scale $\exp \big\{ - g n^{1/12} \big\}$ for some small constant $g > 0$; 
\item and
the stitch points set $S_{n,f}$, whose cardinality will have tail bounded above by the $p$-sequence, and whose elements we will record in the form $\stitch_{n,f;i}$, $i \in \intint{ \# S_{n_0 + n,f}}$. 
\end{itemize}
And, of course, we must have the fundamental relationship that ensures that our weight profiles are indeed being exhibited as patchwork quilts: when the error event $E_{n,f}$ fails to occur, $X_{n,f}$ must equal ${\rm Quilt}[\overline{F}_{n,f},S_{n,f}]$.

Moreover, for given $\ovbar\coninit \in (0,\infty)^3$, 
it suffices to find $n_0 \in \N$ such that 
the above conditions are verified merely if $n \in \N$ satisfies $n \geq n_0$. This weakening of the conditions is permitted because the $q$-sequence may be taken equal to one on an initial finite interval, so that $E_{n,f}$ may be chosen to be the entire sample space for $n < n_0$; thus, no demand is made for such $n$.

We structure our derivation of Theorem~\ref{t.unifpatchcompare}
by first setting up and discussing some structures that are entailed by the occurrence of the  $\normalcoal$ event (with suitable parameter settings) that was analysed in the preceding section. After we have made a little progress in this way, we will be able to specify the characters, and prove the properties, whose need has just been recalled, and thus reach the desired conclusion.

We begin this discussion by recalling the quantities~$D$ and~$\chimac$ that enter as parameters in the definition of the $\normalcoal$ event. 
Henceforth, we set~$D$ equal to $10^{16} c_1^{-4/3}$, a choice that is made in order to permit the application of Lemma~\ref{l.normalcoal}.

Let $\chimac > 0$ given.
We define an $\N$-valued random variable $\Gamma_n$. Its reciprocal will play the role of the parameter~$\e$ in the elaborated rough guide from Section~\ref{s.reroot}.
We have suggested that we will make a decreasing search for $\e$ through dyadic scales. In light of this, it would be natural to  define $\Gamma_n$ to equal~$2^j$
where $j \in \N$ is chosen to be minimal such that   the event 
$\normalcoal_{n;(*:f,0)}^{([-1,1],1)}\big(D, 2^{-j} , \chimac \big)$ occurs.
This definition would seem to entail that the event 
$\normalcoal_{n;(*:f,0)}^{([-1,1],1)}\big(D,  \Gamma_n^{-1} , \chimac \big)$ always occurs. 
However, for this to make sense, it is necessary that the mesh condition $\Gamma_n^{-3/2} \in n^{-1} \Z$ be verified, and so we revisit the definition in order to respect this condition.
For each $j \in \N$, we will write  $\ulcorner \!\! \ulcorner 2^j \urcorner \!\! \urcorner$ for the smallest real value $u \geq 2^j$ such that $u^{-3/2} n \in \Z$. We actually define $\Gamma_n$
to equal   $\ulcorner \!\! \ulcorner 2^j \urcorner \!\! \urcorner$ 
where $j \in \N$ is chosen to be minimal such that   the event 
$\normalcoal_{n;(*:f,0)}^{([-1,1],1)}\big(D, \ulcorner \!\! \ulcorner 2^j \urcorner \!\! \urcorner , \chimac \big)$ occurs.
This certainly entails that 
$\normalcoal_{n;(*:f,0)}^{([-1,1],1)}\big(D,  \Gamma_n^{-1} , \chimac \big)$ always occurs. 
It is due to this occurrence that
\begin{equation}\label{e.cannumbound}
 \canopynumber_{n,;(*:f,0)}^{([-1,1],1)}\big( 1 -  \Gamma_n^{-3/2} \big) \leq   \Gamma_n^{1 + \chimac} \, .
\end{equation}
Before continuing, we mention that
this variation in definition of $\Gamma_n$ to accomodate the mesh condition is a fairly minor detail. For reference shortly, we note that, for given $j \in \N$, the quantity $\ulcorner \!\! \ulcorner 2^j \urcorner \!\! \urcorner$ 
is less than $2^{j+1}$ provided that $n \geq (1 - 2^{-3/2})^{-1} 2^{3j/2}$. We will see that the $j$-value associated to $\Gamma_n$ will very typically verify this condition, so that the concerned dyadic scale does not typically differ between the proposed and actual definitions of $\Gamma_n$. 


Since 
$\normalcoal_{n;(*:f,0)}^{([-1,1],1)}\big(D,  \Gamma_n^{-1} , \chimac \big)$  occurs,  we may find a bound on the tail of $\Gamma_n$ by applying Lemma~\ref{l.normalcoal}. 
We consider $i \in \N$ and apply  this lemma  with ${\bm \e} =    \ulcorner \!\! \ulcorner 2^i \urcorner \!\! \urcorner^{-1}$, and with ${\bf D}$ equal to~$D$, namely $10^{16} c_1^{-4/3}$. 
Provided that $\chi > 0$ is small enough that $\chi \leq 2^{-1} \big( 1 + 500 (\log \beta)^2 \big)^{-1}$, 
we conclude that,   when $n,i \in \N$ satisfy  $i \geq i_0  \big( \chi,\ovbar\coninit \big)$ and $n \geq n_1(\ovbar\coninit) + \Cnew \, 2^{504 i}$,
\begin{equation}\label{e.gammanineq.old}
 \PP \big( \Gamma_n \geq  \ulcorner \!\! \ulcorner 2^i \urcorner \!\! \urcorner \big) \, \leq \,      \ulcorner \!\! \ulcorner 2^i \urcorner \!\! \urcorner^{-2} \cdot \Omega \big(   \ulcorner \!\! \ulcorner 2^i \urcorner \!\! \urcorner^{-1} \big)    \, .
\end{equation}
Here, $2^{-i_0}$ is the upper bound on $\e$ in  Lemma~\ref{l.normalcoal}.
The quantity  $n_1(\ovbar\coninit)$
equals the first of the  three terms of which the triple maximum lower bound on $n$ is composed. The constant $\Cnew$ is specified by increasing the value
$10^{740} c_3^{-84} a_0^{-9} \big( 10^{16} c_1^{4/3} \big)^{36}$
by a $K_0$-determined constant so that the assumed lower bound $n \geq n_1(\chi,\ovbar\coninit) + \Cnew \, 2^{504 i}$ is enough to ensure that $n$ is at least the second, as well as the third, of the three quantities in the triple maximum lower bound in Lemma~\ref{l.normalcoal}. 
Naturally, the function $\Omega$ is specified by the value in the lemma where  ${\bf D}$ equals the above value.

For the values of $(n,i)$ that we are considering, it follows from a comment made a moment ago that   $\ulcorner \!\! \ulcorner 2^{i-1} \urcorner \!\! \urcorner$ 
is less than $2^i$. In view of this, and since $\e \to \Omega(\e)$ is decreasing,
we see from~(\ref{e.gammanineq.old}) that
\begin{equation}\label{e.gammanineq}
 \PP \big( \Gamma_n \geq 2^i \big) \, \leq \,    2^{-2i} \cdot \Omega(2^{-i-1})    \, .
\end{equation}
The form of this bound permits us to leave behind the $\ulcorner \!\! \ulcorner \cdot \urcorner \!\! \urcorner$ notation, because the discrepancy between the `dyadic scale' and actual definitions of $\Gamma_n$
has in practice been taken care of.

For $n \in \N$ that satisfies  $n \geq n_1(\ovbar\coninit)$,
let $\imax(n) \in \N$ be the maximal value of $i \in \N$ such that
$n \geq n_1(\ovbar\coninit) +  \Cnew \, 2^{504i}$.
Define the {\em error} event $E_n = \big\{ \Gamma_n \geq 2^{\imax(n)} \big\}$. Since $2^{\imax(n)}$ is up to a bounded factor equal to $n^{1/{504}}$ provided that $n$ exceeds a certain value $n_2(\ovbar\coninit)$, we find that~(\ref{e.gammanineq}) implies that 
\begin{equation}\label{e.errorub}
\PP (E_n) \leq n^{-2/{504}} \cdot  \Omega \big( n^{-1/503} \big) \, .
\end{equation}
A possibly large constant on this right-hand side has been absorbed by a suitable choice of $n_2(\ovbar\coninit)$, in view of the form~(\ref{e.omegaform}) of $\Omega$.

Suppose now that $n \geq n_2(\ovbar\coninit)$ and that $E_n^c$ occurs. Consider any  $\big(n;*:f,1-\Gamma_n^{-3/2} \big)$-canopy $C$; recall that $C \subseteq [-1,1]$.
We now discuss the behaviour of the $f$-rewarded line-to-point weight profile $C \to \R: y \to \weight_{n;(*:f,0)}^{(y,1)}$.
For any such $C$,
 the  root $\rho_{n;(*:f,0)}^{(y,1)}(0)$  is independent of $y \in C$, and we denote by $\kappa = \kappa(C) \in \N$
the quantity  $\kappa = \lfloor \Gamma_n \, \rho_{n;(*:f,0)}^{(y,1)}(0) \rfloor$ for any $y \in C$. 
Moreover, recalling Definition~\ref{d.canopyspecialpoint}, the closure of the canopy $C$ contains a special point $\specialpt = \specialpt_C$. In 
 the case that the special point lies in the interior of $C$, we  
consider the two random functions of $y \in \big(\inf C , \specialpt \big) \subset [-1,1]$ and  $y \in \big( \specialpt, \sup C  \big) \subset [-1,1]$ respectively given by
$$
 \weight_{n;(*:f,0)}^{(y,1)}  -  \weight_{n;\big((\kappa+1)\Gamma_n^{-1} ,0\big)}^{(y,1)}  \, \, \, \textrm{and} \, \, \, \weight_{n;(*:f,0)}^{(y,1)}  -  \weight_{n;(\kappa\Gamma_n^{-1},0)}^{(y,1)} 
 \, .
$$
When $\specialpt = \sup C$, we consider only the first of these functions; and 
when $\specialpt = \inf C$, only the second.
Each function under consideration is {\em constant}, as we learn by applying Lemma~\ref{l.canopytwopieces}
 with parameter settings ${\bm \e} = \Gamma_n^{-1}$ and  ${\bf K} = D \big( \log \Gamma_n \big)^{1/3}$. 
The lemma's conclusion is valid only when a $\nolatecoal$ and a $\regfluc$ event occur. It is the occurrence of 
$\normalcoal_{n;(*:f,0)}^{([-1,1],1)}\big(D,  \Gamma_n^{-1} , \chimac \big)$ which ensures that these events come to pass.
Meanwhile, the lemma's hypothesis that
 $n >  2^{-3/2} \coninit_1^3 \vee 8 (\coninit_2  + 1)^3 \vee 8 ({\bf K}+2)^3$ is satisfied, after a possible increase in the value of $n_2 = n_2(\chi,\ovbar\coninit)$, because the occurrence of $E_n^c$ entails that $\log \Gamma_n \leq \imax(n) \log 2$, an inequality whose  
 right-hand side is up to a bounded factor at most $\log n$.

In this way, on the event $E_n^c$, we split each $\big(n;*:f,1-\Gamma_n^{-3/2} \big)$-canopy $I$ into two pieces (when the special point lies in the interior of $C$), or leave $I$ untouched (in the other case). We call the intervals so formed  $\big(n;*:f,1 - \Gamma_n^{-3/2} \big)$-{\em split canopies}. Let $C$ be such a split canopy. Note that the two weight difference functions in the preceding paragraph may naturally be indexed by the quantities
$\kappa$ and $\kappa+1$. 
We associate to~$C$ the {\em root neighbour index} $\rootneigh = \rootneigh\big(1-\Gamma_n^{-3/2} , C \big) \in \{\kappa,\kappa+1\}$;  this is the index of the weight difference function that is {\em constant} in the variable $y \in C$.


We write $\cannum$, the split canopy cardinality, for the number of  $\big(n;*:f,1 - \Gamma_n^{-3/2} \big)$-split canopies. Each of these intervals is a subset of  
an $\big(n;*:f,1 - \Gamma_n^{-3/2} \big)$-canopy, and thus is contained in $[-1,1]$; moreover, the union of the closures of split canopies comprises $[-1,1]$.
 We may thus record  the collection of split canopies in the form $\big( \zeta_{(n,f),i},\zeta_{(n,f),i+1} \big)$, $1 \leq i \leq \cannum$, where the increasing real-valued sequence $\big\{ \zeta_{(n,f),i}: i \in \intint{\cannum + 1} \big\}$  satisfies $\zeta_{(n,f),1} = -1$ and $\zeta_{(n,f),\cannum+1} = 1$. This sequence's terms are of two types, namely boundary points of $\big(n;*:f,1-\Gamma_n^{-3/2} \big)$-canopies, and special points interior to such canopies. No two consecutive terms of the latter type are possible, so that $\cannum \leq 2 \, \canopynumber_{n,;(*:f,0)}^{([-1,1],1)}\big(1 - \Gamma_n^{-3/2} \big)$. 

We use a shorthand under which, for each $i \in \intint{\cannum}$, we write $\rootneighuse_i$ for the value of the root neighbour index $\rootneigh\big(1 - \Gamma_n^{-3/2} , (\zeta_{(n,f)i},\zeta_{(n,f),i+1}) \big)$ associated to $i\textsuperscript{th}$ of the above  $\big(n;*:f,1 - \Gamma_n^{-3/2} \big)$-split canopies. 
On the event $E_n^c$, and for each $1 \leq i \leq \cannum$,   we define $Y_{(n,f),i}:[-1,1] \to \R$ by setting $Y_{(n,f),i}(y) =  \weight_{n;\big( \rootneighuse_i \cdot \Gamma^{-1}_n, 0\big)}^{(y,1)}$ for  $y \in [-1,1]$. We also use  a formal device that permits $Y_{(n,f),i}$ to be defined on the whole probability space.
For given $(n,i)$, it remains to specify $Y_{(n,f),i}$
on the event $E_n \cup \big\{   i \geq \cannum + 1 \big\}$.
We specify that the conditional distribution of $Y_{(n,f),i}$, given this event, equals the standard Brownian bridge law $\mc{B}_{1;0,0}^{[-1,1]}$ (which we will denote by $\mu$ later in the proof).

\noindent{\bf Proof of Theorem~\ref{t.unifpatchcompare}.}
We may now begin the formal derivation of this result, because we are ready to specify the sequences and events (whose need has been recalled at the beginning of this section) that will show that our weight profiles are suitably quiltable. Recall also that it is enough to verify the concerned conditions only when $n \in \N$
satisfies $n \geq n_0$, for $n_0$ determined by $\ovbar\coninit$.

We will next specify the definitions and then explain why they enjoy the necessary properties and relationships. It is the verification that fabric sequence elements uniformly withstand $L^{3-}$-comparison to Brownian bridge above a small scale that still requires a little effort. We present this verification after the others.

First the definitions.

\noindent{\em The sequences $p$ and $q$.} 
We may set $\chi > 0$
so that $2(1+\chi)^{-1} = 2 - \e/2$ in order that we may take $p_j =  j^{-2 + \e/2} \cdot 16 \Omega  \big( (j/4)^{-(1+\chi)^{-1}}  \big)$.
We set $q_n = n^{-1/252} \Omega(n^{-1/503})$.

\noindent{\em The error event $E_{n,f}$.}
We set $E_{n,f} = E_n$. 

\noindent{\em The fabric sequence elements.}
We set $F_{(n,f),i} = Y_{(n,f),i}$ for each $(n,f,i) \in \N \times   \initcond_{\ovbar\coninit} \times \N$.

\noindent{\em The stitch points set $S_{n,f}$.} This set's cardinality will be chosen to be $\cannum - 1$. Its elements will be $\stitch_{(n,f),i} = \zeta_{(n,f),i}$ for $2 \leq i \leq \cannum$.  In a convention that is anyway compatible with the endpoint $\zeta$ values, we set $\stitch_{(n,f),1} = -1$  and  $\stitch_{(n,f),\# S_{n,f} + 1} = 1$.

And now the verification of the desired properties.

\noindent{\bf The tail of the two sequences.}
 Since $\Omega(\e) = \e^{o(1)}$ as $\e \searrow 0$,  there exists $j_0 \in \N$ such that $16 \Omega\big( 2^{-1}(j/4)^{-(1+\chi)^{-1}}  \big) \leq j^{\e/2}$ when $j \geq j_0$.
Thus, $p_j \leq j^{-2 + \e}$ whenever $j \geq j_0$.

Since $\Omega(\e) = \e^{o(1)}$ as $\e \searrow 0$, the sequence $q$ verifies $q_n \leq n^{-1/252 + \e}$
for $n$ sufficiently high.

\noindent{\bf The error event's tail.}
In view of the definition of $q_n$ and~(\ref{e.errorub}), we have that
$\PP(E_{n,f}) \leq q_n$.

\noindent{\bf The $p$-sequence dominates the tail of the stitch points' cardinality.}
Recall that $\# S_{n,f} \leq \cannum \leq 2 \, \canopynumber_{n,;(*:f,0)}^{([-1,1],1)}\big(1 - \Gamma_n^{-3/2} \big)$.
By~(\ref{e.cannumbound}) and~(\ref{e.gammanineq}),
$$
\PP \big( \# S_{n,f} \geq 2^{i(1+\chi) + 1} \big) \leq \PP \big( \Gamma_n \geq 2^i \big) \, \leq \,  2^{-2i} \cdot \Omega(2^{-i-1})
$$ 
for $i \in \N$ satisfying $i \geq i_0(\chi,\ovbar\coninit)$. Since $\e \to \Omega(\e)$ is decreasing, 
$$
\PP \big( \# S_{n,f} \geq j \big) \leq 
 j^{-2(1+\chi)^{-1}}  \cdot 16 \, \Omega \big( 2^{-1} (j/4)^{-(1+\chi)^{-1}}  \big) = p_j
$$
for $j$ sufficiently large.

\noindent{\bf The  weight profiles are indeed patchwork quilts.}
We must check that, when the error event $E_{n,f}$ fails to occur, $X_{n,f}$ equals ${\rm Quilt}[\overline{F}_{n,f},S_{n,f}]$.
This amounts to confirming that, on~$E_{n,f}^c$, in each patch $\big[ \stitch_{(n,f),i-1},\stitch_{(n,f),i} \big]$, $i \in \intint{\# S_{n,f} + 1}$, 
the weight profile $y \to \weight_{n;(*:f),0}^{(y,1)}$ differs from the fabric sequence element $F_{n,f;i} = Y_{n,i}$ by a constant (albeit a random one).
Each patch is a split canopy, except the extreme ones, and these are subintervals of split canopies. Thus, this property follows by our construction.



\noindent{\bf Fabric sequence elements uniformly withstand $L^{3-}$-comparison to Brownian bridge above scale $\exp \big\{ - g \, n^{1/12} \big\}$.}
All that remains to complete the proof of our main theorem is to confirm this assertion for a suitable choice of the constant $g > 0$.
To wit, writing  $\mu = \mc{B}_{1;0,0}^{[-1,1]}$,
we have to show that, for any measurable set $A \subseteq \mc{C}_{0,0} \big( [-1,1] , \R \big)$ such that $\mu(A) \geq \exp \big\{ - g n^{1/12} \big\}$, there exists, for every $\eta \in (2/3)$, a constant $C_0$ such that 
\begin{equation}\label{e.presentaim}
 \PP \big( Y_{n,i}^{[-1,1]} \in A  \big) \leq C_0 \big( \mu(A) \big)^{2/3 - \eta} \, .
 \end{equation}
  This constant may depend on $\eta$, but not on $(n,f,j) \in \N \times  \initcond_{\ovbar\coninit} \times \N$ for $n \geq n_0$.


We now present over several paragraphs an argument leading to this conclusion.
Let $\big\{ a_j: j \in \N \big\}$ denote the enumeration of rationals in $[-1,1]$ of the form $p 2^{-i}$ with $p,i \in \N$ coprime in which these values are recorded in increasing order of the dyadic scale~$i$, and increasingly for values of given scale: the sequence begins $a_1 = -1$, $a_2 = 0$, $a_2 = 1$, and continues by enumerating $-1/2$ and $1/2$, and then $-3/4$, $-1/4$, $1/4$ and $3/4$.

For all $n,m \in \N$, we define the  random function $Z_{n,m}: [-1,1] \to \R$, setting  
$$
Z_{n,m}(y) = \weight_{n;( a_m , 0)}^{(y,1)} \, \, \textrm{ for } \, \,  y \in [-1,1] \, .
$$

For $n \in \N$ and  $i \in \intint{\cannum}$,
we introduce the $\N$-valued random variable $M_{n,i}$ by declaring that $a_{M_{n,i}} = \rootneighuse_i \cdot \Gamma^{-1}_n$. Thus, the random function $Y_{n,i}: [-1,1] \to \R$ may be written $Y_{n,i}(y) = Z_{n,M_{n,i}}(y)$ for  $y \in [-1,1]$. 
The root neighbour index  $\rootneighuse_i$
is defined only when $E_n^c \cap \big\{ \cannum \leq i \big\}$ occurs.  To specify $M_{n,i}$
on the whole probability space, we further set $M_{n,i} = 0$ on $E_n \cup \big\{ i \geq \cannum + 1 \big\}$.

For any $n,i \in \N$, $a_{M_{n,i}}$
is a dyadic rational in $[-1,1]$ whose dyadic scale is at most the base-two logarithm of $\Gamma_n$. The number of such rationals in $[-1,1]$ is $2 \Gamma_n + 1$. 
Thus, for each $k \in \N$, $k \geq 1$, 
$$
\PP \big( \sup_{i \in \N} M_{n,i} \geq 2^k + 2  \, , \, E_n^c \big) \leq 
\PP \big( \Gamma_n > 2^{k-1}, \, E_n^c \big) =  
\PP \big( 2^k \leq \Gamma_n \leq 2^{\imax(n) - 1} \big)   \, ,
$$ 
where, in verifying the equality, it may be useful to recall that $E_n$ and $\imax(n)$ have been defined in the paragraph ending at~(\ref{e.errorub}).
When $k \leq \imax(n) - 1$, we may apply~(\ref{e.gammanineq})
with ${\bf i} = k$. The hypothesis that $n \geq n_1(\ovbar\coninit) + \Cnew \, 2^{504k}$
is validated by the definition of $\imax(n)$. We thus find that
\begin{equation}\label{e.whateverk}
\PP \big(  \sup_{i \in \N} M_{n,i} \geq 2^k +  2 \, , \, E_n^c \big)  \leq 2^{-2k} \cdot  \Omega(2^{-k-1})
\end{equation}
for $i_0 \leq k \leq \imax(n) - 1$. Since the event $2^k \leq \Gamma_n \leq 2^{\imax(n) - 1}$ cannot occur if $k \geq \imax(n)$,
the bound~(\ref{e.whateverk}) is seen to be valid provided merely that $k \geq i_0$.

Since $\e \to \Omega(\e)$ is decreasing, we find that
 $$
\PP \big(  \sup_{i \in \N} M_{n,i} \geq j \, , \, E_n^c \big)  \leq 2^2  j^{-2} \cdot \Omega \big( j^{-1}/4 \big)  
$$
whenever $n \in \N$
and $j \geq j_0$, where $j_0 = 2^{i_0 + 1}$.

Let $p_{n,i,j} = \PP \big( M_{n,i} = j \, , \, E_n^c \big)$. 
Let $A \subseteq \mc{C}_{0,0} \big([-1,1],\R \big)$ be measurable. For $n,i \in \N$, note that 
$$
 \PP \big( Y_{n,i}^{[-1,1]} \in A \, , \, E_n^c  \big) = \sum_{j \in \N} \PP \big( Y_{n,i}^{[-1,1]} \in A \, , \, M_{n,i} = j \, , \, E_n^c  \big) \, .
$$
Let $\ell \in \N$. Note that, when $\ell \geq j_0$,
$$ 
\sum_{j = \ell}^{\infty} \PP \Big( Y_{n,i}^{[-1,1]} \in A \, , \, M_{n,i} = j\, , \, E_n^c  \Big)  \leq \PP \big( M_{n,i} \geq \ell \, , \, E_n^c  \big) \leq    \ell^{-2} \cdot  4 \, \Omega  \big( \ell^{-1}/4 \big)  \, ,
$$
and that
$$
 \sum_{j = 1}^{\ell - 1} \PP \Big( Y_{n,i}^{[-1,1]} \in A \, , \, M_{n,i} = j\, , \, E_n^c  \Big) \leq  
 \sum_{j = 1}^{\ell - 1} \PP \Big( Z_{n,j}^{[-1,1]} \in A \, , \, E_n^c  \Big)  \, .
$$

Recall that $Z_{n,j}(y) = \weight_{n;(a_j,0)}^{(y,1)}$. The process $Z_{n,j}: [-1,1] \to \R$
is an instance of the random function $\mc{L}$ seen in Theorem~\ref{t.bridge}. 
Indeed, the process $\mc{L}$ in this theorem is specified by three parameters, and with the settings ${\bf K} = -1$, ${\bf \ipdval} = 2$ and ${\bf x} = a_j$, it coincides with $Z_{n,j}$ on $[-1,1]$. Set $g > 0$
to be the value $G^{-1}$ specified by Theorem~\ref{t.bridge}.
 We 
 find that, provided that the event $A$ is such that $\mu(A) \geq \exp \big\{ - g n^{1/12} \big\}$ and that $n$ is at least a certain constant $n_2$,
$$
 \PP \big( Z_{n,j}^{[-1,1]} \in A  \big)  \, \leq \,  \mu(A) \cdot \conbrac  \exp \Big\{  \conbrac \big( \log  \mu(A)^{-1}   \big)^{5/6} \Big\} \, ,
$$
where the constant $\conbrac$  is finite; note that the final sentence of Theorem~\ref{t.bridge} implies that $\conbrac$ may be chosen independently of $j \in \N$ provided that $4 \leq 2^{-1} c (n+1)^{1/9}$ (and this bound is assured by $n \geq n_2$ 
by suitable choice of the value of $n_2$).
%
%

We find that, for such events $A$ and for $n \in \N$  at least $n_2$,
$$
 \sum_{j = 1}^{\ell - 1} \PP \Big( Y_{n,i}^{[-1,1]} \in A \, , \, M_{n,i} = j \, , \, E_n^c  \Big) \, \leq \, (\ell - 1)  \cdot \mu (A) \cdot \conbrac  \exp \Big\{ \conbrac \big( \log  \mu(A)^{-1}    \big)^{5/6} \Big\} \, .
$$


When $M_{n,i}=0$, $Y_{n,i} = Y_{n,i}^{[-1,1]}$
was artificially declared to have the law $\mu$. Thus,
$$
\PP \big( Y_{n,i}^{[-1,1]} \in A \, , \, M_{n,i} = 0 \big) \, \leq \, \mu (A) \, .
$$

Set $\ell = \lceil \mu(A)^{-1/3} \rceil$. 
We find that, when $n$ is at least $n_2$,
\begin{eqnarray*}
 \PP \big( Y_{n,i}^{[-1,1]} \in A  \big) & \leq & 
   \mu(A)^{2/3}  \Omega [A]   \, + \,  \mu(A)^{2/3}  \conbrac  \exp \big\{ \conbrac \big( \log \mu(A)^{-1}   \big)^{5/6} \big\} \, + \, \mu(A) \\
   & \leq & \mu(A)^{2/3} \big( \Omega[A] +  \conbrac + 1 \big)   \exp \big\{ \conbrac \big( \log \mu(A)^{-1}   \big)^{5/6} \big\}  
    \end{eqnarray*}
  provided that $\mu(A) \leq e^{-1}$. Here, $\Omega[A]$ denotes $4 \, \Omega \big( 4^{-1} \lceil \mu(A)^{-1/3} \rceil^{-1} \big)$. Recalling form~(\ref{e.omegaform}) of $\Omega$, we see that the last right-hand side is bounded above by $\mu(A)^{2/3}  \exp \big\{ \conbrac \big( \log \mu(A)^{-1}   \big)^{5/6} \big\}$ if we make a suitable increase in the value of $\conbrac$. 
 
  Our present aim is to establish~(\ref{e.presentaim})
  with the stated uniformity in the constant $C_0$. This aim is achieved by the bound that we have derived, provided that $C_0$ is raised if need be so that values of $\mu(A)$ above a certain given small value may be accomodated.
Which is to say:
that  fabric sequence elements $F^{[-1,1]}_{(n,f),j}$ uniformly withstand $L^{3-}$-comparison to Brownian bridge above scale $\exp \big\{ - g n^{1/12} \big\}$ has been demonstrated; and
the proof of  Theorem~\ref{t.unifpatchcompare} is complete. \qed

\newpage

\appendix

\section{Glossary of notation}\label{s.glossary}

This article uses quite a lot of notation. Each line of the list that we now present recalls one of the principal pieces of notation; provides a short summary of its meaning; and gives the number of the page at which the notation is either introduced or formally defined. The summaries are, of course, imprecise: phrases in quotation marks  are merely verbal approximations of a precise meaning. 

\vspace{10mm}

\bigskip
\def\qq{&}

\begin{center}
\halign{
\!\!\!\!\!\!\!\!\!\!\!\!#\quad \!\!\! \! \hfill&#\quad\hfill&\!\!\!\!\!\quad\hfill#\cr
 $\weight_{n;(*:f,0)}^{(y,1)}$ \qq the maximum $f$-rewarded line-to-point weight into $(y,1)$ \fff{fweight}
  $\initcond_{\ovbar\coninit}$ \qq class of initial data (specifying choices of $f$), indexed by vector $\ovbar\coninit \in (0,\infty)^3$ \fff{initcond}
  $\mc{B}_{0,0}^{[a,b]}$ \qq the law of standard Brownian 
bridge $B:[a,b] \to \R$, $B(a) = B(b) = 0$ \fff{bridge}
 staircase \qq (a geometric depiction of) an unscaled Brownian LPP path  \fff{staircase}
energy \qq the value assigned to a staircase by Brownian LPP  \fff{energy} 
geodesic \qq a staircase of maximum energy given its endpoints  \fff{geodesic}     
$R_n$ \qq the linear {\em scaling map} \fff{scalingmap}
zigzag \qq the image of a staircase under the scaling map  \fff{zigzag}  
polymer \qq a zigzag that is the image of a geodesic -- and thus of maximum {\em weight}  \fff{polymer} 
weight  \qq the scaled energy, assigned to any zigzag  \fff{weight} 
compatible triple \qq a triple $(n,t_1,t_2)$ so that $[t_1,t_2]$
is the lifetime of a zigzag  \fff{comptriple}    
$\tot$ \qq the difference $t_2 - t_1$; the lifetime of a given polymer in most applications \fff{tot}
$\rho_{n;(x,t_1)}^{(y,t_2)}$ \qq the polymer of journey $(x,t_1) \to (y,t_2)$  \fff{polynot}
$\weight_{n;(x,t_1)}^{(y,t_2)}$ \qq the weight of the polymer $\rho_{n;(x,t_1)}^{(y,t_2)}$  \fff{maxweight} 
$\maxpoly_{n;(I,t_1)}^{(J,t_2)}$ \qq the maximum cardinality of a set of polymers travelling $(I,t_1) \to (J,t_2)$ \fff{maxcard}
 $\rho_{n;(*:f,0)}^{(y,1)}$ \qq the $f$-rewarded line-to-point polymer ending at $(y,1)$  \fff{fpolymer}
 $\polyunique_{n;0}^1$ \qq the (full measure) set  of $(x,y) \in \R^2$ for which $\rho_{n;(x,0)}^{(y,1)}$ is well defined \fff{polyunique}
  $\polyunique_{n;(*:f,0)}^1$ \qq the set, also of full measure, of $y \in \R$ for which  $\rho_{n;(*:f,0)}^{(y,1)}$ is well defined \fff{fpolyunique}
$\pdr_{n;(x,t_1)}^{(y,t_2)}\big(a,r\big)$ \qq $\big\{ \!\!$ $\rho_{n;(x,t_1)}^{(y,t_2)}$ has `normalized' fluctuation $\leq r$ at lifetime fractions $a$ and $1-a$ $\!\!\big\}$  \fff{polydevreg}
$\regfluc_{n;(*:f,0)}^{\big( \{-1,1 \}, 1 \big)}(R)$ \qq $\big\{\!\!$ $f$-rewarded polymers into $[-1,1] \times \{1\}$ begin in $[-R-1,R+1] \times \{0\}$ $\!\!\big\}$ \fff{regfluc}
$\latecoal_{n;([-K,K],0)}^{([-1,1],1)}(\e)$ \qq $\big\{\!\!$  three polymers start $\e$-close; all meet at end; no contact in $\R \times [1 - \e^{3/2}]$  $\!\!\big\}$ \fff{latecoal} 
$\tau_{n;\big(  \{ x_1,x_2  \} , 0 \big)}^{\uparrow;( y , 1)}$ \qq the first time that $\rho_{n;(x_1,0)}^{(y,1)}$ and $\rho_{n;(x_2,0)}^{(y,1)}$ meet \fff{forwardcoal}
$\tau_{n;(x,0)}^{\downarrow;(\{y_1,y_2\},1)}$ \qq the last time that  $\rho_{n;(x,0)}^{(y_1,1)}$ and $\rho_{n;(x,0)}^{(y_2,1)}$ meet \fff{backwardcoal} 
$\big(n,*:f,s\big)$-canopy \qq maximal interval towards which all $f$-rewarded polymers meet after time $s$ \fff{canopy} 
 $\canopynumber_{n,;(*:f,0)}^{([-1,1],1)}(s)$ \qq the cardinality of the set of $(n;*:f,s)$-canopies that intersect $[-1,1]$ \fff{canopynumber}
$\normalcoal_{n;(*:f,0)}^{([-1,1],1)}\big(D,\e,\chimac \big)$ \qq an event specifying that coalescence structure behaves `normally' \fff{normalcoal}
}\end{center}

\newpage

\section{Polymer uniqueness}\label{s.polyunique}

Here we prove the polymer uniqueness Lemma~\ref{l.densepolyunique}.
We first retreat to unscaled coordinates, and prove the counterpart result in a more general form, addressing uniqueness of the energy maximizer among systems of pairwise disjoint staircases with given endpoints.

To state this general result, Lemma~\ref{l.severalpolyunique},
we introduce some notation. The concerned staircases are not in fact entirely disjoint, not least because their endpoints are shared. 
Adopting a weaker notion than disjointness, we say that two staircases are  horizontally separate if there is no planar horizontal interval of positive length that is a subset of a horizontal interval in both staircases.







Now our notation for collections of pairwise horizontally separate staircases. 
For $\ell \in \N$, let $(x_i,s_i)$ and $(y_i,f_i)$, $i \in \intint{\ell}$, be a collection of pairs of elements of $\R \times \N$. (The symbols $s$ and $f$ are used in reference to the staircases' heights at the {\em start} and {\em finish}.)

Let $\staircase^\ell_{(\bar{x},\bar{s}) \to (\bar{y},\bar{f})}$
denote the set of $\ell$-tuples $(\phi_1,\cdots,\phi_\ell)$,
where $\phi_i$ is a staircase from $(x_i,s_i)$ to $(y_i,f_i)$
and each pair $(\phi_i,\phi_j)$, $i \not= j$, is  horizontally separate. 
(This set may be empty; in order that it be non-empty, it is necessary that  $y_i \geq x_i$ and $f_i \geq s_i$
for  $i \in \intint{\ell}$.) 
Note also that 
$\staircase^1_{(x_1,s_1) \to (y_1,f_1)}$ equals 
$\staircase_{(x_1,s_1) \to (y_1,f_1)}$. 

We also associate an energy to each member of  $\staircase^\ell_{(\bar{x},\bar{s}) \to (\bar{y},\bar{f})}$.
Each of the $\ell$ elements of any $\ell$-tuple in $\staircase^\ell_{(\bar{x},\bar{s}) \to (\bar{y},\bar{f})}$ has an energy, as we described in Subsection~\ref{s.staircases}. Define the energy $E\big( \phi \big)$ of any $\phi = \big( \phi_1,\cdots,\phi_\ell \big) \in \staircase^\ell_{(\bar{x},\bar{s}) \to (\bar{y},\bar{f})}$ to be $\sum_{j=1}^\ell E(\phi_j)$.

When  $\upright^\ell_{(\bar{x},\bar{s}) \to (\bar{y},\bar{f})} \not= \emptyset$, we further define the maximum $\ell$-tuple energy
$$
   M^\ell_{(\bar{x},\bar{s}) \to (\bar{y},\bar{f})}  = \sup \Big\{ E(\phi): \phi \in  \upright^\ell_{(\bar{x},\bar{s}) \to (\bar{y},\bar{f})} \Big\} \, .
$$

\begin{lemma}\label{l.severalpolyunique}
Let $\ell \in \N$ and let $(x_i,s_i)$ and $(y_i,f_i)$, $i \in \intint{\ell}$, be a collection of  pairs of points in $\R \times \N$
such that  $\upright^\ell_{(\bar{x},\bar{s}) \to (\bar{y},\bar{f})}$ is non-empty.
Then, 
except on a $\PP$-null set, there is a unique element of
$\upright^\ell_{(\bar{x},\bar{s}) \to (\bar{y},\bar{f})}$
whose energy attains $M^\ell_{(\bar{x},\bar{s}) \to (\bar{y},\bar{f})}$. 
\end{lemma}

\noindent{\bf Proof of Lemma~\ref{l.severalpolyunique}.} 
We model the proof on a resampling argument that shows 
a simple 

\medskip

\noindent{\bf Claim}. Standard Brownian motion $B:[0,1] \to \R$, $B(0) = 0$, has a unique maximizer.

\medskip

We begin by deriving the claim. In this toy case, we may note that the event that the maximizer is non-unique is a subset of the union over pairs of disjoint intervals $I,J \subset [0,1]$ with rational endpoints such that $\sup I < \inf J$ of the event $A_{I,J;B}$ that the maximum value $B[{\rm max}, I]  := \sup_{x \in I} B(x)$ coincides with its counterpart $B[{\rm max}, J]$.
Writing $\PP$ for the probability measure that carries the process $B$, we may augment the probability space associated to this measure by equipping it with a further {\em resampled} process $B_{J}^{{\rm re}}$ that will also be distributed as standard Brownian motion on $[0,1]$ under $\PP$.
There are some variations on the rules for specifying this new process that prove the desired result in this toy case; we select one that is directly adaptable to the Brownian LPP polymer uniqueness problem that we have in mind. 
Let $B':[0,\infty) \to \R$
denote a standard Brownian motion that is independent of $B$.
Label $J = [j_1,j_2]$ for $j_1,j_2 \in [0,1]$.
For $x \in [0,1]$, we define 
\begin{equation*}
 B_{J}^{{\rm re}}(x) = 
\begin{cases} 
 B(x) & x\in [0,j_1]  \, , \\
   B(j_1) + B'(x - j_1)   & x \in J \, , \\ 
   B(j_1) + B'(j_2 - j_1) + B(x) - B(j_2)   & x \in [ j_2 , 1 ] \, .
\end{cases}
\end{equation*}
The process $B_{J}^{{\rm re}}$ is readily seen to have independent Gaussian increments of the necessary variance and thus to have the law of standard Brownian motion.

As such, we may note that $\PP\big(A_{I,J;B}\big) = \PP\big( A_{I,J;B_{J}^{{\rm re}}} \big)$, where naturally the latter notation refers to the event of coincidence of interval maxima for the resampled process. Let $\sigma[B]$ denote the $\sigma$-algebra generated by $B:[0,1] \to \R$. 
Write $\PP_{\sigma[B]}(\cdot) = \EE \big[ {\bf 1}_{\cdot}  \big\vert \sigma[B] \big]$ for conditional probability given knowledge of $B$.
Note then that 
$$
 \PP\big( A_{I,J;B_{J}^{{\rm re}}} \big) = \EE \Big[  \PP_{\sigma[B]}\big( A_{I,J;B_{J}^{{\rm re}}} \big) \Big] \, .
$$

Note that  $A_{I,J;B_{J}^{{\rm re}}}$ occurs precisely when $B(j_1) + \sup_{x \in [0,j_2 - j_1]} B'(x)$ equals the resampled process maximum on $I$, which we may denote by $B_{J}^{{\rm re}}[{\rm max}, I]$.
Because $I$ is assumed to lie to the left of $J$, $B_{J}^{{\rm re}}[{\rm max}, I]$ equals $B[{\rm max}, I]$.
Conditionally on $B$, this equality characterizing the occurrence of $A_{I,J;B_{J}^{{\rm re}}}$   may be expressed as asserting that the conditionally random quantity $\sup_{x \in [0,j_2 - j_1]} B'(x)$ equals the known value $B[{\rm max}, I] - B(j_1)$. 
The conditionally random quantity has the conditional law of the maximum of a standard Brownian motion on the given length interval $[0,j_2-j_1]$; by the reflection principle, this law has a density and thus is non-atomic. Hence,  $\PP_{\sigma[B]}\big( A_{I,J;B_{J}^{{\rm re}}} \big)$ is seen to be zero, $\PP$-almost surely. 

In this way, we prove the claim: 
we
find that 
 $\PP\big( A_{I,J;B_{J}^{{\rm re}}} \big)$, and thus also $\PP(A_{I,J;B})$, equals zero, as we sought to show.

We now adapt this argument to prove Lemma~\ref{l.severalpolyunique}.
A horizontal planar line segment with integer height is called a {\em horizontal rational interval} if its two endpoints have rational $x$-coordinate.
Recall that our data  $\ell \in \N$ and $(x_i,s_i)$, $(y_i,f_i)$, $i \in \intint{\ell}$, is supposed such that 
  $\upright^\ell_{(\bar{x},\bar{s}) \to (\bar{y},\bar{f})}$ is non-empty. 
  
 Say that  an element $\phi$ of the set 
 $\upright^\ell_{(\bar{x},\bar{s}) \to (\bar{y},\bar{f})}$  
 {\em inhabits} a given horizontal rational interval~$I$ if the interval~$I$ is contained in one of the horizontal segments of one of the $\ell$ staircase components of $\phi$. Say that such an element {\em avoids}~$I$ if 
each of these components is disjoint from $I$.
  Let  $A_I$ denote the event that there exist an element of  $\upright^\ell_{(\bar{x},\bar{s}) \to (\bar{y},\bar{f})}$  of maximal weight that inhabits~$I$ and another such element that avoids~$I$.

Note then that the event that there are two elements of  $\upright^\ell_{(\bar{x},\bar{s}) \to (\bar{y},\bar{f})}$ each of whose energies attains $M^\ell_{(\bar{x},\bar{s}) \to (\bar{y},\bar{f})}$ is, up to 
 a $\PP$-null set error, contained in $\bigcup A_I$, where the union is taken over all horizontal rational intervals~$I$.  

We now prove by resampling that $\PP(A_I) = 0$ for any given horizontal rational interval~$I$. 

We augment the probability space carrying the law~$\PP$
with an auxiliary ensemble $B_{I}^{{\rm re}}: \Z \times \R \to \R$ whose law will coincide with that of $B$.
To do so, let  $B':\R \to \R$ denote standard two-sided Brownian motion, independent of the ensemble $B$ under $\PP$.
Write $I = [z_1,z_2] \times \{ k \}$ for $z_1,z_2 \in \R$, $z_1 \leq z_2$, and $k \in \Z$.
For $x \in \R$, we define 
\begin{equation*}
 B_{I}^{{\rm re}}(k,x) = 
\begin{cases} 
 B(k,x) & x \leq z_1  \, , \\
   B(k,z_1) + B'(x - z_1)   & x \in [z_1,z_2] \, , \\ 
   B(k,z_1) + B'(z_2 - z_1) + B(k,x) - B(k,z_2)   & x \geq z_2  \, .
\end{cases}
\end{equation*}

Note that the increments of $B_{I}^{{\rm re}}(k,\cdot)$ coincide with those of $B(k,\cdot)$ away from $I$, while these increments are determined by independent Brownian randomness on this interval. The process $B_{I}^{{\rm re}}(k,\cdot)$ thus shares the law of  $B(k,\cdot)$.

We also set  $B_{I}^{{\rm re}}(j,x) = B(j,x)$
for all $(j,x) \in \Z \times \R$, $j \not= k$, so that the new ensemble shares the law of $B$ under $\PP$.

Any $\ell$-tuple $\phi$ of staircases has an energy $E(\phi)$
specified by increments of the ensemble $B$; it also has a counterpart energy specified in terms of the new ensemble $B_{I}^{{\rm re}}$. We will write $E_{I}^{{\rm re}}(\phi)$ for this new, resampled, energy.
Similarly, the event $A_I$ is specified by the randomness of $B$, and has a counterpart, which we denote by $A^{{\rm re}}_I$,
when the role of this randomness is played by  $B_{I}^{{\rm re}}$. 

The equality in law between $B$ and  $B_{I}^{{\rm re}}$ implies that
$\PP(A_I) = \PP\big(A^{\rm re}_I\big)$.
To prove Lemma~\ref{l.severalpolyunique}, we have seen that it is enough to show that 
$\PP(A_I)  = 0$ for  a given horizontal rational interval $I$. We now complete the proof by fixing such an interval~$I = [z_1,z_2] \times \{ k \}$ and arguing that $ \PP\big(A^{\rm re}_I\big) = 0$.

Any element $\phi \in \upright^\ell_{(\bar{x},\bar{s}) \to (\bar{y},\bar{f})}$ that avoids~$I$ satisfies $E_I^{{\rm re}}(\phi) = E(\phi)$. 
Any element  $\phi \in \upright^\ell_{(\bar{x},\bar{s}) \to (\bar{y},\bar{f})}$ that inhabits~$I$ undergoes a random but $\phi$-independent change of energy under the resampling experiment, satisfying
$E_I^{{\rm re}}(\phi) = E(\phi) + \Theta$, where $\Theta =    B'(z_2-z_1) - \big( B (k,z_2) -  B (k,z_1) \big)$.

Let the avoidance sum energy $\mathsf{ASE}$ denote the supremum of $E(\phi)$ over those  $\phi \in \upright^\ell_{(\bar{x},\bar{s}) \to (\bar{y},\bar{f})}$ that avoid~$I$. Let the inhabitance sum energy $\mathsf{ISE}$ denote the supremum of $E(\phi)$ over those $\phi \in \upright^\ell_{(\bar{x},\bar{s}) \to (\bar{y},\bar{f})}$ that inhabit~$I$. When these quantities are considered after resampling, with $E_I^{{\rm re}}(\phi)$ in place of $E(\phi)$, we prefix the term {\em resampled} to their names.

Since $A^{{\rm re}}_I$ is the event that the resampled avoidance and inhabitance sum energies are equal, this event occurs precisely when $\Theta$ equals the difference  between the avoidance and inhabitance sum energies. 
Let $\sigma[B]$ denote the $\sigma$-algebra
 generated by the ensemble $B:\Z \times \R \to \R$
 and write $\PP_{\sigma[B]}$ for the associated conditional probability. 
We have then that
$$
 \PP \big( A^{{\rm re}}_I \big) = \EE \Big[ \PP_{\sigma[B]} \big( \Theta = \mathsf{ASE} - \mathsf{ISE} \big)  \Big] \, .
$$
To the observer of the ensemble $B$,  $\mathsf{ASE} - \mathsf{ISE}$ is a known quantity, while  
the conditional distribution of $\Theta$ is a normal random variable of mean~$B(k,z_1) -  B(k,z_2)$ and variance $z_2 - z_1$. Thus, $\PP_{\sigma[B]} \big( \Theta = \mathsf{ASE} - \mathsf{ISE} \big)$ equals zero, $\PP$-almost surely, and so 
$\PP \big( A^{{\rm re}}_I \big)$ equals zero, as we sought to show. 
 This completes the proof of Lemma~\ref{l.severalpolyunique}. \qed

\medskip

Let $h: \R \to \R \cup \{-\infty\}$ be measurable.
For $x,y \in \R$, $x \leq y$ and $s,f \in \Z$, $s \leq f$, we associate to any element 
$\phi \in \upright^1_{(x,s) \to (y,f)}$ 
 the $h$-rewarded energy $E(\phi) + h(x)$. 
 The maximum line-to-point $h$-rewarded energy is then defined to be 
 $$
  M^1_{(*:h,s) \to (y,f)}  : = \sup \Big\{ E(\phi) + h(x): \phi \in  \upright^1_{(x,s) \to (y,f)} \, , \, x \leq y \Big\} \, .
$$
\begin{lemma}\label{l.hpolyunique} 
Let $y \in \R$.
Let $h: \R \to \R \cup \{-\infty \}$ be measurable, with $h(x) > -\infty$ for some $x < y$, and  $\limsup_{x \to -\infty}h(x)/\vert x \vert < 0$.
Let $s,f \in \Z$ satisfy $s \leq f$.
Then, 
except on a $\PP$-null set, 
$M^1_{(*:h,s) \to (y,f)}$ is finite, and  
there is a unique choice of  $(x,\phi)$ such that $x \in (-\infty,y]$, $\phi  \in   \upright^1_{(x,s) \to (y,f)}$
and the $h$-rewarded energy of $\phi$ attains $M^1_{(*:h,s) \to (y,f)}$. 
\end{lemma}
\noindent{\bf Proof.}
In this result, we consider the $h$-rewarded energy $E(\phi) + h(x)$ of staircases $\phi$ that end at a given location $(y,f) \in \R \times \N$
and that begin at a location $(x,s) \in \R \times \N$ whose height $s \leq f$ is fixed but whose first coordinate $x$ may vary on $(-\infty,x]$.

We must first argue that the maximum $h$-rewarded energy $M^1_{(*:h,s) \to (y,f)}$ adopted among such staircases $\phi$ is a finite quantity.
Due to our assumption that $\limsup_{x \to -\infty} \vert x \vert^{-1} h(x)$ is negative, it suffices to argue that 
$M^1_{(x,s) \to (y,f)}$,
the energy maximum over such staircases that begin at $(x,s)$, grows sublinearly in $\vert x \vert$ the limit $x \to -\infty$.
(In fact, we also need to know that the supremum of 
$M^1_{(x,s) \to (y,f)}$ over $x$ in any bounded interval bordered on the right by $y$ is almost surely finite.)
A crude estimate on 
$M^1_{(x,s) \to (y,f)}$ is
$$
M^1_{(x,s) \to (y,f)} \, \leq \, \sum_{k = s}^f \Big( \sup_{z \in [x,y]} B(k,z) -  \inf_{z \in [x,y]} B(k,z)  \Big)
$$
The distribution of the summand is unchanged if we suppose that $B(k,x)$ is zero. 
Supposing that this is so, the summand $A - B$ may be viewed in the form $A + (-B)$.
It is then a sum of two terms that share their distribution. By the reflection principle, this distribution assigns mass $2 \cdot (2\pi)^{-1/2} (y-x)^{-1} \exp \{ - 2^{-1} (y-x)^{-1}  r^2   \}$
to any interval $(r,\infty)$ for $r \geq 0$.
We see then that
\begin{equation}\label{e.xybound}
 \PP \Big(  
M^1_{(x,s) \to (y,f)} \geq 2 K r \Big) \leq 2K \cdot 2 \cdot (2\pi)^{-1/2} (y-x)^{-1} \exp \{ - 2^{-1} (y-x)^{-1}  r^2   \}
\end{equation}
for any $r \geq 0$; here and later, we write $K = f - s +1$.
Fixing any $\eta >0$, we may 
choose $r_x = \big( 2^{1/2} + \eta \big) (y-x)^{1/2} \big( \log (y-x) \big)^{1/2}$, we may sum this bound over $x \in y - \N$ to find that 
the event 
$M^1_{(x,s) \to (y,f)} \geq 2K r_x$ may occur for only finitely many such $x$.

Some understanding of the local regularity of 
$M^1_{(x,s) \to (y,f)}$ as $x$ varies is now needed to treat the general case of $x \in (-\infty,y]$.
Let $x_1,x_2 \leq y$.
Note that 
$$
M^1_{(x_1,s) \to (y,f)} - 
M^1_{(x_2,s) \to (y,f)} \geq   B(s,x_2)   -    B(s,x_1)  \, .
$$
Writing this in the form  
$M^1_{(x_2,s) \to (y,f)} - M^1_{(x_1,s) \to (y,f)} \leq   B(s,x_1)   -    B(s,x_2)$, we may fix 
 $x_1 \leq y$, and find that
$$
\sup_{x_2 \in [x_1 - 1, x_1]} \Big( M^1_{(x_2,s) \to (y,f)} - 
M^1_{(x_1,s) \to (y,f)} \Big) \leq   B(s,x_1)  \, -  \,  \inf_{x_2 \in [x_1 - 1,x_1]} B(s,x_2)  \, .
$$

For any $r > 0$,
the right-hand side exceeds $r$ with probability  $2 \cdot (2\pi)^{-1/2}  \exp \{ - 2^{-1}   r^2   \}$ by a second use of the reflection principle.
Combining with~(\ref{e.xybound}), we see that, for any $x \leq y - 1$,
\begin{equation}\label{e.xyk}
 \PP \Big(  
 \sup_{z \in [x-1,x]} M^1_{(z,s) \to (y,f)} \geq 2(K+1) r \Big) \leq 2(2K + 1) (2\pi)^{-1/2} (y-x)^{-1} \exp \{ - 2^{-1} (y-x)^{-1}  r^2   \}.
\end{equation}
Maintaining our choice $r_x = \big( 2^{1/2} + \eta \big) (y-x)^{1/2} \big( \log (y-x) \big)^{1/2}$, we see that, for any given $\eta > 0$,
the set of $x \leq y$ for which 
$M^1_{(x,s) \to (y,f)} \geq 2(K+1)r_{x-1}$ is almost surely bounded below. 
We thus confirm that $M^1_{(x,s) \to (y,f)}$ grows sublinearly in $\vert x \vert$ as $x \to -\infty$.
To prove the finiteness of  $M^1_{(*:h,s) \to (y,f)}$, it remains only to prove the almost sure finiteness of the supremum of $M^1_{(x,s) \to (y,f)}$  as $x$ varies over any bounded interval bordered by $y$ on the right. However, this is again a consequence of~(\ref{e.xybound}) and~(\ref{e.xyk}).

To complete the proof of Lemma~\ref{l.hpolyunique}, it remains to argue that there is a unique choice of $(x,\phi)$
with $x \leq y$ and $\phi \in D^1_{(x,s) \to (y,f)}$
for which $E(\phi) + h(x)$ attains the maximum value $M^1_{(*:h,s) \to (y,f)}$. 
The resampling proof of Lemma~\ref{l.severalpolyunique} works here.
Note the simplification that we consider staircases rather than $\ell$-tuples of staircases: or formally, we take $\ell = 1$.

We revisit the resampling proof, now 
considering any given horizontal rational interval $I$ contained in $(-\infty,y] \times \{ s,f \}$, and 
defining any staircase in  $D^1_{(x,s) \to (y,f)}$ for any $x \leq y$
to inhabit or avoid $I$ as we did before. The event $A_I$
is now specified by the existence of one staircase in the union of   $D^1_{(x,s) \to (y,f)}$ over $x \leq y$
that inhabits $I$ and of another in the same union that avoids $I$. 
As before, it is enough to show that $\PP(A_I) = 0$ for any given horizontal rational interval $I$, this time  contained in $(-\infty,y] \times \{ s,f \}$. 
The proof that this probability is zero is unchanged. 

This completes the proof of Lemma~\ref{l.hpolyunique}. \qed

\medskip

\noindent{\bf Proof of Lemma~\ref{l.densepolyunique}: (1).}  The staircase set $\upright^1_{(2 n^{2/3}x,0) \to (n+2 n^{2/3}y,n)}$, which is in correspondence under the scaling map $R_n:\R^2 \to \R^2$ with the collection of $n$-zigzags in question,
is non-empty if and only if $y \geq x - n^{1/3}/2$.
When this last condition is satisfied, the polymer weight  $\mathsf{W}_{n;(x,0)}^{(y,1)}$ equals
$$
   2^{-1/2} n^{-1/3}   \Big( M^1_{(2 n^{2/3}x , 0) \to (n+2 n^{2/3}y,n)} - n -  2n^{2/3}(y-x)  \Big) \, .
$$
The uniqueness of the $n$-polymer in question then follows from 
 Lemma~\ref{l.severalpolyunique} with ${\bm \ell} = 1$, ${\bf x_1} = 2 n^{2/3}x$, ${\bf y_1} = n+2 n^{2/3}y$, ${\bf s_1} = 0$ and ${\bf f_1} = n$.

\noindent{\bf (2).} Let $y \in [-1,1]$. It is enough to check that $y \in \polyunique_{n;(*:f,0)}^1$ almost surely.
 Recall that the $f$-rewarded line-to-point polymer weight  $\mathsf{W}_{n;(*:f,0)}^{(y,1)}$ equals
\begin{eqnarray*}
 & &   \sup_{x \in (-\infty,n^{1/3}/2 + y]} \bigg(  2^{-1/2} n^{-1/3}   \Big( M^1_{(2 n^{2/3}x , 0) \to (n+2 n^{2/3}y,n)} - n -  2n^{2/3}(y-x)  \Big)  \, + \, f(x) \bigg)  \\
& = &    2^{-1/2} n^{-1/3}  \sup_{u \in (-\infty,n + v]} \Big(    M^1_{(u , 0) \to (n+v,n)} - n -   v+ u    +  h(u) \Big)    \, ,
\end{eqnarray*}
where $h:\R \to \R \cup \{ - \infty \}$ is given by
$h(x) = 2^{1/2} n^{1/3}  f\big(n^{-2/3}x/2 \big)$.  We took $v= 2 n^{2/3}y$ and $u = 2 n^{2/3}x$ to obtain the displayed equality.
We seek to apply Lemma~\ref{l.hpolyunique} with ${\bf h}:  \R \to \R \cup \{-\infty\}$
given by ${\bf h}(u) =  - n -   v+ u    +  h(u)$, ${\bf s} = 0$, ${\bf f} = n$ and ${\bf y} = n + v$. Indeed,
$\mathsf{W}_{n;(*:f,0)}^{(y,1)}$ has been shown to equal  $2^{-1/2} n^{-1/3}  M^1_{(*:{\bf h},{\bf s}) \to ({\bf y},{\bf f})}$, so that the application of the lemma will 
 imply that $y \in \polyunique_{n;(*:f,0)}^1$ almost surely. 
It remains to check that Lemma~\ref{l.hpolyunique}'s hypotheses  are satisfied. In this regard, note that
\begin{eqnarray*}
{\bf h}(u) & \leq & - n -v + u + 2^{1/2} n^{1/3} \coninit_1 \big(1 + n^{-2/3} \vert u \vert/2  \big) \\
   & \leq &  - n -v + u + 2^{1/2} n^{1/3} \coninit_1  +   2^{-1/2} \coninit_1 n^{-1/3} \vert u \vert \, ,
\end{eqnarray*}
so that the hypothesis that $\limsup_{u \to -\infty}{\bf h}(u)/\vert u \vert < 0$
is satisfied when $n > 2^{-3/2} \coninit_1^3$. The other hypothesis, that ${\bf h}(u) > -\infty$ for some $u \leq n+v$, is validated due to  $\sup_{x \in [-\coninit_2,\coninit_2]} f(x) \geq - \coninit_3$
and $2n^{2/3}\coninit_2 \leq n + v$ or equivalently
$n \geq 8(\coninit_2  - y)^3$. \qed

\section{Calculational derivations}\label{s.calcder}

In Subsection~\ref{s.rolehyp}, we explained that one aspect of the proofs of two of our results,  Proposition~\ref{p.latecoal}  and Lemma~\ref{l.normalcoal}, has been separated from the rest. 
The missing piece is  the `calculational derivations' of these results. Here, we present these derivations. In the case of each concerned result, it is now our task to verify that the result's hypotheses are sufficient for the purpose of verifying every condition that is invoked during the course of the result's proof.

A guideline of the format of the derivations' presentation is that first all of the conditions invoked in the proof of the result in question are recorded. These conditions may then be split into lists, according to which parameters are being bounded. Each list is then separately simplified to produce a shorter list that implies the original. Sometimes in carrying out this simplification, new conditions, called {\em further conditions}, will be introduced. The simplified lists may then be compared, with some further simplifications being noted. The derivation ends when a final list is drawn up whose conditions coincide with the hypotheses of the result under review.

Some notational conventions govern the presentation of the derivations.  Conditions are given square bracketed names such as $[1]$ or $[n4]$, shown on the left. 
These names may be recycled from one calculational derivation to the next.
The notation $[1,n4]$ means `conditions $[1]$ and $[n4]$'. Implication is denoted by a right arrow, so that $[1,n4] \to [r3]$ means `conditions $[1]$ and $[n4]$ imply condition $[r3]$'.  

The notation $[n2,4,5]$ is used as a shorthand for $[n2,n4,n5]$. This meaning would be ambiguous, were there named conditions $[4]$ or $[5]$, but we employ this condition only when the meaning is unambiguous.

There are three parts to the section. In the first, a small explicit computation from  Proposition~\ref{p.latecoal}'s proof is provided.
The next two provide the calculational derivations of this proposition and of Lemma~\ref{l.normalcoal}.

\subsection{A short piece of working in the proof of Proposition~\ref{p.latecoal}}

We begin by recording a short piece of working needed at the very end of the proof of this result.

\noindent{\bf Deriving the upper bound in the conclusion of  Proposition~\ref{p.latecoal}.}
At the very end of the proof of  Proposition~\ref{p.latecoal}, it is claimed that 
the sum of  $10^{20} K  C  c_1^{-4/3}  \big( \log \e^{-1} \big)^{4/3}    \e^2$ and  the expression in line~(\ref{e.pupperbound}) is at most $\e^2 \cdot 10^{420}  K  c_3^{-33}      C_3   \big( \log \e^{-1} \big)^{42}     \exp \big\{ \beta_3 \big( \log \e^{-1} \big)^{5/6} \big\}$. Here is the working: 
\begin{eqnarray*}
 & &
\e^2 \cdot 10^{356}  K r^4  c_3^{-27} C_3   \big( \log \e^{-1} \big)^{36}     \exp \big\{ \beta_3 \big( \log \e^{-1} \big)^{5/6} \big\}
+  10^{20} K  C  c_1^{-4/3}  \big( \log \e^{-1} \big)^{4/3}    \e^2  \\
 & =  & 
\e^2 \cdot 10^{356}  K \cdot 10^{176/3} 2^{44/3} c_1^{-16/3}  \big( \log \e^{-1} \big)^{16/3} \cdot   c_3^{-27} C_3   \big( \log \e^{-1} \big)^{36}     \exp \big\{ \beta_3 \big( \log \e^{-1} \big)^{5/6} \big\} \\
 & & \qquad \qquad \qquad
+  \, 10^{20} K  C  c_1^{-4/3}  \big( \log \e^{-1} \big)^{4/3}    \e^2 \\
& \leq  & 
\e^2 \cdot 2^{-1} 10^{420}  K  c_3^{-33}      C_3   \big( \log \e^{-1} \big)^{42}     \exp \big\{ \beta_3 \big( \log \e^{-1} \big)^{5/6} \big\}
+  10^{20} K  C  c_3^{-4/3}  \big( \log \e^{-1} \big)^{4/3}    \e^2 \\  
& \leq & \e^2 \cdot 10^{420}  K  c_3^{-33}      C_3   \big( \log \e^{-1} \big)^{42}     \exp \big\{ \beta_3 \big( \log \e^{-1} \big)^{5/6} \big\} \, .
\end{eqnarray*} 
In the penultimate inequality, we used $10^{356 + 176/3}  2^{44/3} = 10^{414 + 2/3} 26007 \cdots = 10^{414 + 2/3 + 4.4 \cdots} \in [ 10^{419}, 1/2 \cdot 10^{420} ]$, $\e \leq e^{-1}$ and $c_3 \leq c_1 \leq 1$. 
In the final inequality, we used $C_3 \geq C$, $c_3 \leq 1$ and $\e \leq e^{-1}$.

\subsection{Proposition~\ref{p.latecoal}: derivation}

Here we present the calculational derivation of Proposition~\ref{p.latecoal}.

We begin by collecting together all the conditions that are invoked during the proof of this proposition.
The three parameters in the proposition's statement are $n$, $\e$ and $K$. It is in terms of these that the various conditions are expressed, sometimes by means of a quantity $r > 0$.
The value of this quantity is during the proof of the proposition set to be equal to
$$
 r = 10^{44/3} 2^{11/3} c_1^{-4/3}  \big( \log \e^{-1} \big)^{4/3} \, .
$$

\noindent{\em Record of all the conditions invoked during the proof of Proposition~\ref{p.latecoal}.}

In the application of Theorem~\ref{t.disjtpoly}, the bounds   $\e \leq (1 - 2^{-3/2})^{2/3}$,   
$r \geq 1$ and $K+2 \leq \e^{-1/2}$
are used. The conditions ${\bf \tot} \geq 1/2$ and  $({\bf \tot})^{2/3} \geq 1/2$ are also used. Since ${\bf \tot} = 1 - \e^{3/2}$, these are implied by $\e \leq 2^{-2/3}$.
In summary, these bounds are implied by
$$
  \e \leq 2^{-2/3}
$$
$$
r \geq 1
$$ 
and 
$$
K+2 \leq \e^{-1/2} \, .
$$

The conditions needed to satisfy the hypotheses of  Theorem~\ref{t.disjtpoly} are noted to be
$$
 (r+2)\e \leq \min \Big\{ \, (\eta_0)^{36} \, , \, 10^{-616}  c_3^{22}  3^{-115} \,  , \, \exp \big\{ - C^{1/4} \big\} \, \Big\} \, ,
$$
$$
    n  \geq 2 \max \bigg\{ \, 2(K_0)^9 \big( \log \e^{-1} \big)^{K_0}
    \, , \, 10^{606}   c_3^{-48}  3^{240}  \rsc^{-36} \e^{-222}  \max \big\{   1  \, , \,   (K+2)^{36} 2^{24}  \big\} \, , \, 
  a_0^{-9} (K+2)^9 2^{6}  \,  \Bigg\} \, .
$$
$$
  2 \big( K +1 + 2\e \big) \leq (r+2)^{-1/2} \e^{-1/2}  \big( \log (r+2)^{-1}\e^{-1} \big)^{-2/3} \cdot 10^{-8} c_3^{2/3} 3^{-10/3}
$$
and
$$
(r+2) \e \leq e^{-4/3}
$$

Finally, the hypothesis that  $({\bf \tot})^{-2/3} \vert {\bf y} - {\bf x} \vert \leq  {\bm \e}^{-1/2} \big( \log {\bm \e}^{-1} \big)^{-2/3} \cdot 10^{-8} c_3^{-2/3} 3^{-10/3}$ is ensured by  
$$
  2 \big( K +1 + 2\e \big) \leq (r+2)^{-1/2} \e^{-1/2}  \big( \log (r+2)^{-1}\e^{-1} \big)^{-2/3} \cdot 10^{-8} c_3^{-2/3} 3^{-10/3}
$$
since $\vert {\bf x} - {\bf y} \vert \leq K +1 + 2\e$ and  $({\bf \tot})^{2/3} \geq 1/2$, where we also used that the function $u \to u^{-1/2} \big( \log u^{-1} \big)^{-2/3}$ is decreasing on $(0,e^{-4/3})$ alongside~(\ref{e.etwobounds}).

The application of Lemma~\ref{l.maxpolytriple} makes use of $\e \leq 2^{-2/3}$.
After this, the bounds $\e \leq 1/3$,  $K \geq 1$ and $r \geq 1$ are used.

When Theorem~\ref{t.polyfluc} is applied, the bounds $r \geq 1$, 
$$
K \geq 2
$$ and 
$$
\e \leq 1/4
$$ 
are used.
In order that this theorem can be applied, it is noted that it suffices that 
 $$
 \e^{3/2} \leq  10^{-11} c_1^2 \, ,
 $$
$$
 n    \geq  \max \bigg\{ 
 10^{32} \e^{-75/2} c^{-18} \, \, , \, \, 
 10^{24} 2^{36} c^{-18} \e^{-75/2}  K^{36}  
 \bigg\} \, ,
$$
$$ 
 r  \geq 
 2 \max \bigg\{  10^9 c_1^{-4/5} \, \, , \, \,  15 C^{1/2} \, \, , \, \,  174 \e^{1/2}  K   \bigg\} \, ,
 $$
and
$$ 
r \leq 6  \e^{25/6} n^{1/36} \, .
$$ 

\noindent{\bf Dividing the conditions into several lists.}

The listed conditions are now partitioned into lists. A further condition, $[\e4]$, is added at this stage.

\noindent{\em The $r$ list.} In this list, we gather all conditions that concern the parameter $r$. Lower bounds on $r$ appear first.

$$
[r1] \, \, r \geq 1
$$ 
$$ 
[r2] \, \,  r  \geq 
 2 \max \bigg\{  10^9 c_1^{-4/5} \, \, , \, \,  15 C^{1/2} \, \, , \, \,  174 \e^{1/2}  K   \bigg\} \, ,
 $$
$$
[r3] \, \,  (r+2)\e \leq \min \Big\{ \, (\eta_0)^{36} \, , \, 10^{-616}  c_3^{22}  3^{-115} \,  , \, \exp \big\{ - C^{1/4} \big\} \, \Big\} \, ,
$$
$$
[r4] \, \,  2 \big( K +1 + 2\e \big) \leq (r+2)^{-1/2} \e^{-1/2}  \big( \log (r+2)^{-1}\e^{-1} \big)^{-2/3} \cdot 10^{-8} c_3^{2/3} 3^{-10/3}
$$
$$ 
[r5] \, \, r \leq 6  \e^{25/6} n^{1/36} \, .
$$
The needed conditon that $(r+2) \e \leq e^{-4/3}$ is implied by $[r3]$ (since $c_3 \leq 1$) and we do not label it. 

We also have the formula for the value of $r$:

$$
[rV] \, \,  r = 10^{44/3} 2^{11/3} c_1^{-4/3}  \big( \log \e^{-1} \big)^{4/3} \, .
$$

\noindent{\em The $n$ list.} This list consists of the remaining conditions that concern $n$. 

$$
 [n1] \, \,    n  \geq 2 \max \bigg\{ \, 2(K_0)^9 \big( \log \e^{-1} \big)^{K_0}
    \, , \, 10^{606}   c_3^{-48}  3^{240}  \rsc^{-36} \e^{-222}  \max \big\{   1  \, , \,   (K+2)^{36} 2^{24}  \big\} \, , \, 
  a_0^{-9} (K+2)^9 2^{6}  \,  \Bigg\} \, .
$$
$$
 [n2] \, \, n    \geq  \max \bigg\{ 
 10^{32} \e^{-75/2} c^{-18} \, \, , \, \, 
 10^{24} 2^{36} c^{-18} \e^{-75/2}  K^{36}  
 \bigg\} \, ,
$$

\noindent{\em The $\e$ list.} This list consists of the remaining conditions that concern $\e$, as well as the new condition~$[\e4]$.
 
$$
[\e1] \, \, \e \leq 1/4
$$ 
 $$
 [\e2] \, \,  \e^{3/2} \leq  10^{-11} c_1^2 \, ,
 $$
 $$
[\e3] \, \, K+2 \leq \e^{-1/2} \, .
$$
$$
[\e4] \, \, \e \leq 10^{-65} c_1^6
$$

\noindent{\em The $K$ list.} This list consists of the remaining condition that concerns $K$.

$$
[K1] \, \, K \geq 2 \, .
$$ 

\noindent{\bf Beginning the analysis of the conditions.}
First we make a claim.

\noindent{\em Claim.} When all the listed conditions are in force, the bound 
$$
\e \leq 10^{-65} c_1^6
$$
implies that
$$
 (r + 2) \e \leq \e^{1/2} \, .
$$
\noindent{\em Proof.}
By $[rV]$, the condition $(r + 2) \e \leq \e^{1/2}$ is given by 
$$
  10^{44/3} 2^{11/3} c_1^{-4/3}  \big( \log \e^{-1} \big)^{4/3}    + 2 \leq \e^{-1/2}
$$
which is implied by 
$$
 2 \cdot  10^{44/3} 2^{11/3} c_1^{-4/3}  \big( \log \e^{-1} \big)^{4/3}   \leq \e^{-1/2}
$$
and 
$$
4 \leq \e^{-1/2} \, ;
$$
the second of these, namely $\e \leq 4^{-2}$, follows from $[\e2]$ and $c_1 \leq 1$. 
Note that
$$
 2 \cdot  10^{44/3} 2^{11/3} c_1^{-4/3}  \big( \log \e^{-1} \big)^{4/3}   \leq \e^{-1/2}
 $$
holds provided that 
$$
\big( \log \e^{-1} \big)^{4/3}  \leq \e^{-1/4}
$$
and
$$
 2 \cdot  10^{44/3} 2^{11/3} c_1^{-4/3}   \leq \e^{-1/4} \, .
$$
The second of these bounds is equivalent to 
$$
 \e \leq  2^{-4}  10^{-176/3} 2^{-44/3} c_1^{16/3}    \, .
$$
Since $2^{4 +  44/3}  10^{176/3} = 2^{56/3} 10^{176/3} = 1.93 \cdots \times 10^{64}$, and $c_1 \leq 1$,
we see that the last displayed bound is implied by $\e \leq 10^{-65} c_1^6$, which is the hypothesis in the claim.

The first of the two bounds, namely $\big( \log \e^{-1} \big)^{4/3}  \leq \e^{-1/4}$, is equivalent to 
$$
 \e^{-1} \leq \exp \big\{ \e^{-3/16} \big\}
$$
which is implied by 
$$
 \e^{-1} \leq (6!)^{-1} \big( \e^{-3/16} \big)^6 = (6!)^{-1} \e^{-9/8} 
$$
which is implied by $\e \leq (6!)^{-8}$, a condition which follows from the claim's hypothesis $\e \leq 10^{-65} c_1^6$ alongside $c_1 \leq 1$. This proves the claim. \qed

\noindent{\em Simplifying the $r$ list.}
We use the claim to simplify the $r$ list. 

Note that $[\e4]$ is the hypothesis of the claim.

We see from the claim that, when $[\e4]$ holds, $[r3]$ is implied by 
$$
[r6] \, \,  \e \leq \min \Big\{ \, (\eta_0)^{72} \, , \, 10^{- 1232}  c_3^{44}  3^{- 230} \,  , \, \exp \big\{ - 2 C^{1/4} \big\} \, \Big\} \, ,
$$
Note that $[r2] \to [r1]$ due to $c_1 \leq 1$.

Since $\e \leq 1/2$ (by $[\e4]$), the condition $[r4]$ is implied by 
$$
 \, \,  2 \big( K  + 2 \big) \leq (r+2)^{-1/2} \e^{-1/2}  \big( \log \e^{-1} \big)^{-2/3} \cdot 10^{-8} c_3^{2/3} 3^{-10/3} \, .
$$
By the claim, it is thus implied by 
$$
 \, \,  2 \big( K  + 2 \big) \leq  \e^{-1/4}  \big( \log \e^{-1} \big)^{-2/3} \cdot 10^{-8} c_3^{2/3} 3^{-10/3} \, .
$$
We have argued that we have $\big( \log \e^{-1} \big)^{4/3} \leq \e^{-1/4}$, so we see that $[r4]$ is implied by 
$$
 \, \,  2 \big( K  + 2 \big) \leq  \e^{-1/8}  \cdot 10^{-8} c_3^{2/3} 3^{-10/3} \, .
$$
Thus, $[r4]$ is implied by 
$$
 [r7] \, \, \e \leq   \big( K  + 2 \big)^{-8}  c_3^{16/3} 2^{-8} 10^{-64} 3^{-80/3} \, .
$$

The condition $[r5]$ is, in view of $[rV]$ and $\log \e^{-1} \leq \e^{-1}$, implied by 
$$ 
 10^{44/3} 2^{11/3} c_1^{-4/3}   \leq 6  \e^{25/6 + 4/3} n^{1/36} 
$$
or equivalently 
$$ 
 10^{528} 2^{132} 6^{-36}  c_1^{-48}  \e^{-198}   \leq  n \, .
$$
Since  $2^{132} 6^{-36} = 5.27 \cdots \times 10^{11}$, this last condition is implied by
$$
  n \geq  10^{540}  c_1^{-48}  \e^{-198} \, .
$$
Since $c_3 \leq c_1$, $c \leq 1$ and $\e \leq 1$, the last is implied by the second condition in $[n1]$. 
That is, $[n1] \to [r5]$.

Thus 
$[r2,6,7,n1,\e4] \to [r1,2,3,4,5]$. 
We see then that we may simplify the $[r]$ list and present it in the form $[r2,6,7]$.

\noindent{\em Simplifying the $n$ list.}
We now make some straightforward comments about the $n$ conditions.

In the second condition of $[n1]$, we see the term   $\max \big\{   1  \, , \,   (K+2)^{36} 2^{24}  \big\}$. Obviously this expression is attained by the latter quantity.
Thus, the lower bound on $n$ in this second condition equals
$$
 10^{606}  2^{24}    3^{240}    c_3^{-48}  \rsc^{-36} \e^{-222}   (K+2)^{36}
 $$
 Since  $2^{24}    3^{240} = 5.41 \cdots \times 10^{121}$, the displayed value is bounded above by  
$$
 10^{728}     c_3^{-48}  \rsc^{-36} \e^{-222}   (K+2)^{36} \, .
 $$
Thus, writing
$$
 [n3] \, \,    n  \geq 2 \max \bigg\{ \, 2(K_0)^9 \big( \log \e^{-1} \big)^{K_0}
    \, , \,   
 10^{728}     c_3^{-48}  \rsc^{-36} \e^{-222}   (K+2)^{36}  \, , \, 
  a_0^{-9} (K+2)^9 2^{6}  \,  \Bigg\} \, ,
$$
we see that $[n3] \to [n1]$.

The condition $[n2]$ is comprised of two lower bounds on $n$. 
The first of these is implied by the second bound in $[n3]$ due to $\e \leq 1$, $c_3 \leq 1$ and $c \leq 1$. 
The second of these is also implied by the same bound for the same reasons. 
Thus, $[n3] \to [n1]$.

We see then that a simplified form of the $n$ list may take the form of $[n3]$.

\noindent{\em Simplifying the $\e$ list.} 
 
$[\e2]$ may be rewritten
$\e \leq  10^{-22/3} c_1^{4/3}$ and thus is implied by 
$$
[\e5] \, \, 10^{-8} c_1^{4/3} \, .
$$
Note that $[\e5] \to [\e1]$ is due to $c_1 \leq 1$. 

$[\e3]$ is equivalent to 
$$  
 [\e6] \, \,  \e \leq   (K+2)^{-2} \, .
$$
Note that $[\e4] \to [\e5]$ is due to $c_1 \leq 1$.

Thus, $[\e4,6] \to [\e1,2,3,4]$, and thus $[\e4,6]$ is a simplified form for the $\e$-list.

\noindent{\em The $K$ list.} This one-condition list will not be simplified.

\noindent{\bf Gathering the simplified lists.}

We may summarise our progress by recalling that the union of our simplified lists is: 
$[r2,6,7]$, $[n3]$, $[\e4,6]$ and $[K1]$. These conditions collectively imply all the bounds needed in the proof of  Proposition~\ref{p.latecoal}.
Our next step is to record this collection of conditions again, taking the opportunity to reallocate some into other lists where they more suitably belong. 
We rename the conditions that change lists, indicating the old and the new names as we record the conditions now.

\noindent{\em The updated  $r$ list.} 
$$ 
[r2 = R1] \, \,  r  \geq 
 2 \max \bigg\{  10^9 c_1^{-4/5} \, \, , \, \,  15 C^{1/2} \, \, , \, \,  174 \e^{1/2}  K   \bigg\} \, ,
 $$

\noindent{\em The updated  $n$ list.}  
$$
 [n3 = N1] \, \,    n  \geq 2 \max \bigg\{ \, 2(K_0)^9 \big( \log \e^{-1} \big)^{K_0}
    \, , \,   
 10^{728}     c_3^{-48}  \rsc^{-36} \e^{-222}   (K+2)^{36}  \, , \, 
  a_0^{-9} (K+2)^9 2^{6}  \,  \Bigg\} \, ,
$$

\noindent{\em The updated  $\e$ list.} 
$$
[\e4 = E1] \, \, \e \leq 10^{-65} c_1^6
$$
$$  
 [\e6 = E2] \, \,  \e \leq   (K+2)^{-2} 
$$
$$
[r6 = E3] \, \,  \e \leq \min \Big\{ \, (\eta_0)^{72} \, , \, 10^{- 1232}  c_3^{44}  3^{- 230} \,  , \, \exp \big\{ - 2 C^{1/4} \big\} \, \Big\} \, ,
$$
$$
[r7 = E4] \, \,   \e \leq   \big( K  + 2 \big)^{-8}  c_3^{16/3} 2^{-8} 10^{-64} 3^{-80/3} \, .   
$$

\noindent{\em The  $K$ list.} 
$$
[K1] \, \, K \geq 2 \, .
$$

\noindent{\em Further analysis.}
A few further simplifications will be made. The most important of these address the condition $[R1]$, which is the remaining instance of a condition that involves the parameter $r$.
This condition should be eliminated and replaced by an upper bound on $\e$. To do this, recall the value of $r$:
$$
[rV] \, \,  r = 10^{44/3} 2^{11/3} c_1^{-4/3}  \big( \log \e^{-1} \big)^{4/3} \, .
$$ 
Since $\e \leq 1$, due to $[E2]$, and $c_1 \leq 1$, $[rV]$ implies the first of the three conditions in $[R1]$.
The third of these conditions is in light of $[E2]$ implied by $r \geq 174$; this bound is implied by the first condition in $[R1]$, since $c_1 \leq 1$.

Using $[rV]$, the second condition in $[R1]$ is seen to be equivalent to
$$
  10^{44/3} 2^{11/3} c_1^{-4/3}  \big( \log \e^{-1} \big)^{4/3}  \geq  30 C^{1/2}
$$
or equivalently
$$
 \e  \leq  \exp \Big\{ -  \big( 10^{-44/3} 2^{-11/3} c_1^{4/3}  \cdot 30     C^{1/2} \big)^{3/4} \Big\} \, .
$$
Since $10^{-44/3} 2^{-11/3}  \cdot 30$ and $c_1$ are both at most one,  this last is  implied by
$$
[R2] \, \, 
 \e  \leq  \exp \big\{ -  2 C^{3/8}  \big\} \, .
$$  
Note that $[R2]$ implies the third condition in $[E3]$ due to  $C \geq 1$.

We may summarise these inferences by introducing
$$
[E5] \, \,  \e \leq \min \Big\{ \, (\eta_0)^{72} \, , \, 10^{- 1232}  c_3^{44}  3^{- 230} \,  , \, \exp \big\{ - 2 C^{3/8} \big\} \, \Big\} \, ,
$$
and noting that $[E1,2,4,5]$ implies the above $E$ list, $[E1,2,3,4]$, as well as the $R$ list $[R1]$.

Having eliminated the parameter $r$ from our conditions, we now further simplify the $E$ list, which is presently $[E1,2,4,5]$.
Note that $[E4] \to [E2]$. 
Writing
$$
[E6] \, \,  \e \leq \min \Big\{ \, (\eta_0)^{72} \, , \, 10^{- 1342}  c_3^{44}  \,  , \, \exp \big\{ - 2 C^{3/8} \big\} \, \Big\} \, ,
$$
we see that $[E6] \to [E5]$ is due to $3^{230} = 5.46 \cdots  \times 10^{109}$. Moreover, $[E6] \to [E1]$ because the second condition in $[E6]$ implies $[E1]$
due to $c_3 \leq c_1$.
Thus we may update the $E$ list to be $[E4,6]$.

Regarding the condition $[N1]$, we merely note that, since $c_3 \leq c$, we have $[N2] \to [N1]$, where
$$
 [N2] \, \,    n  \geq 2 \max \bigg\{ \, 2(K_0)^9 \big( \log \e^{-1} \big)^{K_0}
    \, , \,   
 10^{728}     c_3^{-84}  \e^{-222}   (K+2)^{36}  \, , \, 
  a_0^{-9} (K+2)^9 2^{6}  \,  \Bigg\} \, .
$$
The status report is that all conditions have been shown to be implied by the conditions: $[E1,4,6]$, $[N2]$ and $[K1]$.

The upper bounds on $\e$, $[E4,6]$, may be written:
$$
  \e \leq \min \Big\{   \,   \big( K  + 2 \big)^{-8}  c_3^{16/3} 2^{-8} 10^{-64} 3^{-80/3}  \, , \,      (\eta_0)^{72} \, , \, 10^{- 1342}  c_3^{44}  \,  , \, \exp \big\{ - 2 C^{3/8} \big\} \, \Big\} \, ,
$$
which since $c_3 \leq 1$ is implied by
$$
  \e \leq \min \Big\{   \,      (\eta_0)^{72} \, , \, 10^{- 1342}  c_3^{44}  \big( K  + 2 \big)^{-8}  \,  , \, \exp \big\{ - 2 C^{3/8} \big\} \, \Big\} \, .
$$ 
The lower bound on $n$, $[N2]$, is 
$$
     n  \geq 2 \max \bigg\{ \, 2(K_0)^9 \big( \log \e^{-1} \big)^{K_0}
    \, , \,   
 10^{728}     c_3^{-84}  \e^{-222}   (K+2)^{36}  \, , \, 
  a_0^{-9} (K+2)^9 2^{6}  \,  \Bigg\} \, .
$$
The lower bound on $K$, $[K1]$, is
$$
  K \geq 2 \, .
$$ 
The last three displayed bounds collectively form the hypotheses of Proposition~\ref{p.latecoal}. This completes the calculational derivation of this proposition.

\subsection{Lemma~\ref{l.normalcoal}: derivation}

We now turn to the calculational derivation of this lemma. The result
 has parameters
 $n \in \N$, $\e > 0$,  $D > 0$, $\chimac > 0$ and $\ovbar\coninit \in (0,\infty)^3$. It demands certain hypotheses on these parameters.
 Naturally, our job is to verify that the hypotheses made on the parameters are sufficient to imply all the conditions used during the proof of the result.
 
 Before we begin, we note that one condition, which is not an estimate, is needed in the derivation, when Proposition~\ref{p.latecoal} is applied. This is $n \e^{3/2} \in \N$.
 This condition is indeed hypothesised in Lemma~\ref{l.normalcoal} and we will not refer to it again.
 
\noindent{\bf Record of all the conditions used during the proof of Lemma~\ref{l.normalcoal}}

The proof consists of derivations of upper bounds on the four terms on the right-hand side of~(\ref{e.fourterms}), and then the assembly of these four estimates. We record the used conditions accordingly.

\noindent{\em Bounding $\PP(\manycan)$.}

The first hypothesis is invoked implicitly, when Lemma~\ref{l.densepolyunique}(1) is used.
The condition needed is 
 $\e \lfloor \e^{-1} y_i \rfloor \geq  - D (\log \e^{-1})^{1/3} - n^{1/3}/2$, where recall that the various points $y_i$ each belong to $[-1,1]$.
 Thus, it is sufficient that
 $1 + \e \leq           D (\log \e^{-1})^{1/3} + n^{1/3}/2$. Imposing $\e \leq 1$,
 we see that it suffices that 
 $$
  1 \leq  D (\log \e^{-1})^{1/3}
 $$  
 and 
 $$
  2 \leq n^{1/3} \, .
 $$
  
Next we use $\e^{3/2} \leq 2^{-1}$ and $\e \leq e^{-1}$.
After this, $\e \leq 4^{-1/\chi}$.
Then we use $\e \leq \exp \big\{ - \chi^{-2} \big\}$ and $\chi \leq 1$.  
  
In the application of Theorem~\ref{t.maxpoly},
it suffices that 
 $$
 4^{-1} \e^{-\chimac} \geq  k_0 \vee \big( 4D  \big(\log \e^{-1} \big)^{4/3} + 2 \lceil 8 D \big(\log \e^{-1} \big)^{4/3} \rceil \big)^3  \, ,
 $$ 
and
$$
 n \e^{3/2}   \geq   \max \bigg\{ \, 2(K_0)^{m^2} \big( \log ( 2^{-1} \e^{-\chimac}) \big)^{K_0}
    \, , \,  a_0^{-9} \tau^9 \, , \,    
    10^{325}       \rsc^{-36}   \big( 2^{-1} \e^{-\chimac} \big)^{465}  \max \big\{   1  \, , \,  \tau^{36}   \big\}  
   \,  \bigg\} \, ,
$$ 
where we write $\tau =   4D \big(\log \e^{-1} \big)^{4/3}   + 2 \lceil 8 D \big(\log \e^{-1} \big)^{4/3} \rceil$.

\noindent{\em Bounding $\PP(\latecoal)$.}

In applying Proposition~\ref{p.latecoal}, we make use of
$$
\e  \leq 
\min \Big\{ \exp \big\{ - 8D^{-3} \big\}  \, , \,      (\eta_0)^{72} \, , \, 10^{- 1342}  c_3^{44} \big( D (\log \e^{-1})^{1/3} + 2 \big)^{-8} \,  , \, \exp \big\{ - 2 C^{3/8} \big\} \,
   \Big\} 
$$ 
and
$$
     n  \geq 2 \max \bigg\{ \, 2(K_0)^9 \big( \log \e^{-1} \big)^{K_0}
    \, , \,   
 10^{728}     c_3^{-84}  \e^{-222}   \big( D (\log \e^{-1})^{1/3}  +  2 \big)^{36}  \, , \, 
  a_0^{-9} \big( D (\log \e^{-1})^{1/3}  +   2 \big)^9 2^{6}  \,  \Bigg\} \, .
$$

\noindent{\em Bounding $\PP\big(\neg \, \regfluc\big)$.}

In applying Lemma~\ref{l.regfluc}, we use the conditions
$$
n \geq 
 c^{-18} \max \Big\{  (\coninit_2 + 1)^9 \, , \,   10^{23} \coninit_1^9     \, , \, 3^{9}  \Big\} \, ,
$$ 
and
$$
 D (\log \e^{-1})^{1/3} - 1  \in \Big[  \,  39 \coninit_1  \,  \vee \, 5  \, \vee \,   3 c^{-3}  \, \vee \, 2 \big( (\coninit_2 + 1 )^2 +  \coninit_3 \big)^{1/2}  \, , \,   6^{-1} c n^{1/9} \,  \Big] \, .
$$

\noindent{\em Bounding $\PP\big(\neg \, \cap \pdr \big)$.}

The application of Proposition~\ref{p.polyfluc} is permitted in view of:
$$
\e^{3/2} \leq  10^{-11} c_1^2 \, ,
$$
$$
 n    \geq  \max \bigg\{ 
 10^{32} \e^{-75/2} c^{-18} \, \, , \, \, 
 10^{24} c^{-18} \e^{-75/2} \big( D (\log \e^{-1})^{1/3} + 1 \big)^{36}  
 \bigg\} 
$$
and
$$ 
  D \big(\log \e^{-1} \big)^{4/3} \in \bigg[
 \,  10^9 c_1^{-4/5} \, \vee  \,  15 C^{1/2} \, \vee  \,  87 \e^{1/2}   \big( D (\log \e^{-1})^{1/3} + 1 \big) \, \, , \, \, 
 3  \e^{25/6} n^{1/36}  \, \bigg] \, . 
 $$

\noindent{\em Assembling the estimates.}
Upper bounds are used at this moment of the proof. These are:
\begin{eqnarray*}
 & & \e^{-3 +  \chimac  (580)^{-1} ( \log \beta)^{-2} (\log \log \e^{-1})^2} \\
  & \leq & 
   2^{-1} \condee^{-1}  \big( 2 \lceil 8 D \big(\log \e^{-1} \big)^{4/3}\rceil   \big)^{-(\log \beta)^{-2} (\log \log \e^{-1})^2/{288}  - 3/2} \exp \big\{-  2 (\log  \e^{-1})^{11/12} \big\} \, ,
\end{eqnarray*}
 $$
      \e^{   2^{-13} \rsc   D^3 - 2} \leq   (38)^{-1} D^{-1} (\log \e^{-1})^{-1/3} \rsC^{-1} 
 $$
 and
$$
  \e^{  10^{-11} c_1  D^{3/4}  - 3} \leq 
8^{-1} 22^{-1} \, c_1 C^{-1} D^{-1} \big(\log \e^{-1} \big)^{-4/3} \, .  
$$ 
  
\noindent{\bf Expressing the gathered conditions in lists.}
We now partition the various conditions into several lists.
A few conditions will be expressed in terms of a parameter $\tau$, whose value is now set: 
 $$
[\tau] \, \,  \tau =   4D \big(\log \e^{-1} \big)^{4/3}   + 2 \lceil 8 D \big(\log \e^{-1} \big)^{4/3} \rceil \, .
 $$

\noindent{\em The [$\e$] list.} Here we record all conditions on $\e$ that do not involve~$n$.

 $$
 [\e1] \, \,  1 \leq  D (\log \e^{-1})^{1/3}
 $$  
 $$
 [\e2] \, \, \e^{3/2} \leq 2^{-1}
 $$ 
 $$
 [\e3] \, \, \e \leq e^{-1}
 $$
 $$
 [\e4] \, \, \e \leq 4^{-1/\chi}
 $$
 $$
 [\e5] \, \, \e \leq \exp \big\{ - \chi^{-2} \big\}
 $$
 $$
 [\e6] \, \, 4^{-1} \e^{-\chimac} \geq  k_0 
 $$
 $$
 [\e7] \, \,  4^{-1} \e^{-\chimac} \geq  \tau^3  \, ,
 $$
$$
 [\e8] \, \,  \e  \leq  \exp \big\{ - 8D^{-3} \big\} 
$$
$$   
  [\e9] \, \,   \e  \leq 
  4^{-2}\big( D (\log \e^{-1})^{1/3} + 2 \big)^{-4} 
$$   
$$   
 [\e10] \, \,   \e  \leq 
(\eta_0)^{72} 
$$
$$ 
  [\e11] \, \,   \e  \leq 
10^{- 1342}  c_3^{44} \big( D (\log \e^{-1})^{1/3} + 2 \big)^{-8}
$$
$$ 
 [\e12] \, \,   \e  \leq 
\exp \big\{ - 2 C^{3/8} \big\} 
$$
$$
  [\e13] \, \,  D (\log \e^{-1})^{1/3} - 1  \geq   5 
$$ 
$$ 
  [\e14] \, \,  D (\log \e^{-1})^{1/3} - 1  \geq        3 c^{-3} 
$$  
$$  
 [\e15] \, \,  D (\log \e^{-1})^{1/3} - 1  \geq  \,   39 \coninit_1  \, \vee \, 2 \big( (\coninit_2 + 1 )^2 +  \coninit_3 \big)^{1/2}
$$
$$
  [\e16] \, \,  \e^{3/2} \leq  10^{-11} c_1^2 
$$
$$ 
  [\e17] \, \,   D \big(\log \e^{-1} \big)^{4/3} \geq
 \,  10^9 c_1^{-4/5}
$$
$$ 
  [\e18] \, \,    D \big(\log \e^{-1} \big)^{4/3} \geq   15 C^{1/2} 
    $$
    $$
  [\e19] \, \,   D \big(\log \e^{-1} \big)^{4/3} \geq     87 \e^{1/2}   \big( D (\log \e^{-1})^{1/3} + 1 \big)
    $$
\begin{eqnarray*}
  [\e20] \, \,  & & \e^{-3 +  \chimac  (580)^{-1} ( \log \beta)^{-2} (\log \log \e^{-1})^2} \\
  & \leq & 
   2^{-1} \condee^{-1}  \big( 2 \lceil 8 D \big(\log \e^{-1} \big)^{4/3}\rceil   \big)^{-(\log \beta)^{-2} (\log \log \e^{-1})^2/{288}  - 3/2} \exp \big\{-  2 (\log  \e^{-1})^{11/12} \big\} 
\end{eqnarray*}
 $$
  [\e21] \, \,       \e^{   2^{-13} \rsc   D^3 - 2} \leq   (38)^{-1} D^{-1} (\log \e^{-1})^{-1/3} \rsC^{-1} 
 $$
$$
  [\e22] \, \,   \e^{  10^{-11} c_1  D^{3/4}  - 3} \leq 
8^{-1} 22^{-1} \, c_1 C^{-1} D^{-1} \big(\log \e^{-1} \big)^{-4/3}  
$$ 
   In fact, condition $[\e9]$ is redundant.

\noindent{\em The [$n$] list.} Here are recorded all conditions concerning $n$. 
 $$
 [n1] \, \,  2 \leq n^{1/3} 
 $$
$$
 [n2] \, \, n \e^{3/2}   \geq  
2(K_0)^{(12)^{-2} \big( \log \log   (2^{-1} \e^{-\chimac})     \big)^2}   \big( \log ( 2^{-1} \e^{-\chimac}) \big)^{K_0}
 $$
 $$
  [n3] \, \, n \e^{3/2}   \geq    a_0^{-9} \tau^9 
 $$
 $$   
  [n4] \, \, n \e^{3/2}   \geq   10^{325}       \rsc^{-36}   \big( 2^{-1} \e^{-\chimac} \big)^{465}  \max \big\{   1  \, , \,  \tau^{36}   \big\}   
$$ 
$$
  [n5] \, \,    n  \geq 4(K_0)^9 \big( \log \e^{-1} \big)^{K_0}
$$
$$     
 [n6] \, \, n \geq 2 \cdot 10^{728}     c_3^{-84}  \e^{-222}   \big( D (\log \e^{-1})^{1/3}  +  2 \big)^{36} 
 $$
 $$ 
 [n7] \, \,  n \geq 2 a_0^{-9} \big( D (\log \e^{-1})^{1/3}  +   2 \big)^9 2^{6}  
$$
$$
 [n8] \, \, n \geq 
 c^{-18} \max \Big\{  (\coninit_2 + 1)^9 \, , \,   10^{23} \coninit_1^9   \Big\}
$$ 
$$ 
 [n9] \, \, n \geq 3^{9} c^{-18}   
$$ 
$$
 [n10] \, \, n    \geq   
 10^{32} \e^{-75/2} c^{-18} 
$$ 
$$
 [n11] \, \, n \geq 10^{24} c^{-18} \e^{-75/2} \big( D (\log \e^{-1})^{1/3} + 1 \big)^{36}   
$$
$$ 
 [n12] \, \,  D (\log \e^{-1})^{1/3} - 1  \leq    6^{-1} c n^{1/9} 
$$
$$ 
  [n13] \, \,   D \big(\log \e^{-1} \big)^{4/3} \leq 
 3  \e^{25/6} n^{1/36}   
 $$

\noindent{\em The [$\chi$] list.}  Here appear remaining conditions that concern $\chi$.
    $$
 [\chi1] \, \,   \chi \leq 1
    $$  

\noindent{\em The [$D$] list.}  We also introduce a further condition:
$$
 [D1] \, \, D \geq 1
$$

Before trying to simplify these conditions, we state and prove a lemma.  
\begin{lemma}\label{l.tau.ub}
The conditions
$$
 [\e23] \, \,  \e \leq  \big( 10^{-41} \chi^{36} D^{-14} \big)^{\chi^{-2}} \, , 
   $$
$D \geq 1$ and 
$\chi \leq 2^{1/2}$ imply that 
$$
 22 D \big(\log \e^{-1} \big)^{4/3}  \leq \e^{-\chi^2/6} \, .
$$
\end{lemma}
\noindent{\bf Proof.} 
The conclusion of the lemma holds precisely when
$$
   \e^{-1} \leq \exp \big\{ (22)^{-3/4} D^{-3/4} \e^{-\chi^2/8} \big\}
$$
and thus is implied by 
$$
   \e^{-1} \leq  (\lceil 16 \chi^{-2} \rceil!)^{-1} \big( (22)^{-3/4} D^{-3/4} \e^{-\chi^2/8} \big)^{\lceil 16\chi^{-2} \rceil} \, .
$$
Our hypothesis that $\chi^2 \leq 2$ implies that $\lceil 16 \chi^{-2} \rceil \leq 18 \chi^{-2}$. Since $\e \leq 1$, the last displayed condition  is thus implied by
$$
   \e \leq  (\lceil 16 \chi^{-2} \rceil!)^{-1}  (22)^{-27 \chi^{-2}/2} D^{-27\chi^{-2}/2} 
   $$
 since $D \geq 1$. Since  $\lceil 16 \chi^{-2} \rceil \leq 18 \chi^{-2}$ and $\ell! \leq \ell^\ell$, the last is implied by
$$
   \e \leq  \big( 18\chi^{-2} \big)^{-18 \chi^{-2}}   (22)^{-27 \chi^{-2}/2} D^{-27\chi^{-2}/2} = \Big( \chi^{36} \cdot (18)^{-18} (22)^{-27/2} \cdot D^{27/2} \Big)^{-\chi^{-2}} 
   $$
which since $D \geq 1$
and $(18)^{18} 22^{27/2} = 5.21 \cdots \times 10^{40}$
 is implied by our hypothesis that
$\e \leq  \big( 10^{-41} \chi^{36} D^{-14} \big)^{\chi^{-2}}$. \qed

\noindent{\bf Analysing the $[\e]$ list.}
We now analyse the $[\e]$ list, aiming to produce a simplified list of conditions that collectively imply all conditions in the $[\e]$ list.

We begin by noting that Lemma~\ref{l.tau.ub} provides an upper bound on the quantity $\tau > 0$ specified in~$[\tau]$.  Indeed, since 
 $\tau =   4D \big(\log \e^{-1} \big)^{4/3}   + 2 \lceil 8 D \big(\log \e^{-1} \big)^{4/3} \rceil$, we see that, when $D \geq 1$ and $\e \leq e^{-1}$,
 the quantity $\tau$ is at most $22 D \big(\log \e^{-1} \big)^{4/3}$.
 When the hypotheses of the lemma hold, then, we see that  
 $$
  \tau \leq \e^{-\chi^2/6} \, .
 $$

Introducing
$$
 [\e24] \, \, \e \leq 2^{-4\chi^{-2}} \, ,
$$ 
we next note that 
$$
[\e23,D1,\chi1,\e24] \to [\e7] \, . 
$$
Indeed,
the conditions $[\e23]$, $[D1]$ and $[\chi1]$ imply the hypotheses of Lemma~\ref{l.tau.ub}.
$[\e7]$ is under these circumstances implied by $4^{-1} \e^{-\chi} \geq \e^{-\chi^2/2}$ or equivalently $\e^{\chi^2/2 - \chi} \geq 4$. 
Thus, it is implied by $[\chi1,\e24]$. Thus, $[\e23,D1,\chi1,\e24] \to [\e7]$, as we claimed.

\noindent{\em Analysing [$\e20$].}
We introduce some further conditions in order to replace $[\e20]$:
$$
 [\e25] \, \,  \e \leq     \exp \big\{ - \beta^{ \chimac^{-1/2}  (3480)^{1/2}  } \big\}
$$
$$
 [\e26] \, \, \e \leq e^{-e}
$$ 
$$
 [\e27] \, \,  \e
  \leq (2 \condee)^{- 3480 \chimac^{-1} ( \log \beta)^{2}  }
  $$
$$
 [\e28] \, \, \chimac   \leq  2^{-1}   \big( 1  + 500( \log \beta)^{2} \big)^{-1}   
$$   
and
$$
 [\e29] \, \,      \e  \leq  \exp \big\{ - 2^{12} \big\} \, .
$$

\begin{lemma}\label{l.e.20}
We have that
$$
  [D1,\chi1,\e23,\e25,\e26,\e27,\e28,\e29]  \to [\e20]
$$
\end{lemma}
\noindent{\bf Proof.}
From $[D1,\e26]$, we have $D \big(\log \e^{-1} \big)^{4/3} \geq 1$ and thus $\lceil 8 D \big(\log \e^{-1} \big)^{4/3}\rceil  \leq 9 D \big(\log \e^{-1} \big)^{4/3}$.  We see then that, in these circumstances, $[\e20]$ is implied by
\begin{eqnarray*}
 & & \e^{-3 +  \chimac  (580)^{-1} ( \log \beta)^{-2} (\log \log \e^{-1})^2} \\
  & \leq & 
   2^{-1} \condee^{-1}  \big( 18 D \big(\log \e^{-1} \big)^{4/3}   \big)^{-(\log \beta)^{-2} (\log \log \e^{-1})^2/{288}  - 3/2} \exp \big\{-  2 (\log  \e^{-1})^{11/12} \big\} \, ,
\end{eqnarray*}
or equivalently
\begin{eqnarray*}
 & &-3 +  \chimac  (580)^{-1} ( \log \beta)^{-2} (\log \log \e^{-1})^2 \\ 
 & \geq & \log (2 \condee) (\log \e^{-1})^{-1} \\
 & &  \, + \, \big( (\log \beta)^{-2} (\log \log \e^{-1})^2/{288}  + 3/2 \big) \log \Big( 18 D \big(\log \e^{-1} \big)^{4/3} \Big)  (\log \e^{-1})^{-1}  \, + \,     2 (\log  \e^{-1})^{-1/12} 
\end{eqnarray*}
This condition is satisfied provided that each of the following conditions is met:
$$
   \chimac  (580)^{-1} ( \log \beta)^{-2} (\log \log \e^{-1})^2 \geq 6
$$
$$
 \chimac  (580)^{-1} ( \log \beta)^{-2} (\log \log \e^{-1})^2  \geq 6 \log (2 \condee) (\log \e^{-1})^{-1}
$$
$$
 \chimac  (580)^{-1} ( \log \beta)^{-2} (\log \log \e^{-1})^2  \geq 6\big( (\log \beta)^{-2} (\log \log \e^{-1})^2/{288}  + 3/2 \big) \log \Big( 18 D \big(\log \e^{-1} \big)^{4/3} \Big)  (\log \e^{-1})^{-1}  
$$   
and
$$
 \chimac  (580)^{-1} ( \log \beta)^{-2} (\log \log \e^{-1})^2  \geq 12 (\log  \e^{-1})^{-1/12} 
 $$  
During the proof of Lemma~\ref{l.e.20}, we will call these last four bounds $[1]$, $[2]$, $[3]$ and $[4]$.
 
Recall that the quantity $\beta$ is at least $4$. Thus, $\beta \geq 1$, so that  $[1]$ is implied by 
$$
 \log \log \e^{-1} \geq    \chimac^{-1/2}  (580)^{1/2} ( \log \beta) 6^{1/2}
$$
and thus by $[\e25]$. 
$$
  \e \leq     \exp \big\{ - \beta^{ \chimac^{-1/2}  (3480)^{1/2}  } \big\}
$$

$[2]$ is implied by 
$$
   \chimac  (580)^{-1} ( \log \beta)^{-2} \geq 6 \log (2 \condee) (\log \e^{-1})^{-1}
$$
and $[\e26]$ in the form
$(\log \log \e^{-1})^2   \geq 1$.
The last display is 
$$
  \log \e^{-1}
  \geq 3480 \chimac^{-1} ( \log \beta)^{2}   \log (2 \condee)
  $$
and thus is $[\e27]$.

When $[\e26]$, i.e.,
$(\log \log \e^{-1})^2   \geq 1$, holds, $[3]$ is implied by
$$
 \chimac  (580)^{-1} ( \log \beta)^{-2} \geq 6\big( (\log \beta)^{-2}/{288}  + 3/2 \big) \log \Big( 18 D \big(\log \e^{-1} \big)^{4/3} \Big)  (\log \e^{-1})^{-1}  
$$   
and thus by 
$$
 \chimac  (580)^{-1}  \geq 6\big( 1/{288}  + 3( \log \beta)^{2}/2 \big) \log \Big( 18 D \big(\log \e^{-1} \big)^{4/3} \Big)  (\log \e^{-1})^{-1}  
$$   

By Lemma~\ref{l.tau.ub}, $22D (\log \e^{-1})^{4/3} \leq \e^{-\chimac^2/6}$ is implied by $[\e23,D1,\chi1]$. The last display is thus implied by 
$$
 \chimac  (580)^{-1}  \geq 6\big( 1/{288}  + 3( \log \beta)^{2}/2 \big) \log \Big(  \e^{-\chimac^2/6} \Big)  (\log \e^{-1})^{-1}  
$$   
which is 
$$
  (580)^{-1}  \geq  \chimac \big( 1/{288}  + 3( \log \beta)^{2}/2 \big) 
$$   
which is implied by
$$
  288(580)^{-1}  \geq   \chimac  \big( 1  + 500( \log \beta)^{2} \big) 
$$   
which is implied by $[\e28]$.

When $[1]$ holds, $[4]$ is implied by 
$$
 6 \geq 12 (\log  \e^{-1})^{-1/12} 
$$
which is $[\e29]$.

Thus, the hypotheses of Lemma~\ref{l.e.20} imply $[1,2,3,4]$.
This completes the proof of the lemma. \qed

\noindent{\em Analysing $[\e21]$.}
We now add a condition to the $[D]$ list
$$
 [D2] \, \,  D \geq \rsc^{-1/3} 2^{5}
$$
and another to the $[\e]$ list:
 $$
 [\e30]  \, \,   \e \leq   2^{-1}   \rsC^{-1} 
 $$

\begin{lemma}\label{l.e.21}
We have that
$$
 [D1,D2,\chi1,\e23,\e30]  \to [\e21]  \, .
$$
\end{lemma}
\noindent{\bf Proof.}
By Lemma~\ref{l.tau.ub}, $22D (\log \e^{-1})^{4/3} \leq \e^{-\chimac^2/6}$ is implied by $[\e23,D1,\chi1]$.
Since $[\e23]$ implies that  $\e \leq e^{-1}$, we see that, under these circumstances,
the condition $[\e21]$
is 
implied by $[D2]$ in the guise
$2^{-13} \rsc   D^3 - 2 \geq 2$ alongside   
 $$
      \e^2 \leq   2^{-1}  \e^{\chimac^2/6} \rsC^{-1} \, .
 $$
 Since $\chi \leq 1$ by $[\chi1]$,
 the last display  is implied by $[\e30]$. \qed

\noindent{\em Analysing $[\e21]$.}
We now add conditions 
$$
 [D3] \, \,  D  \geq 5^{4/3} \cdot 10^{44/3} c_1^{-4/3}
$$
and
$$
 [\e31] \, \, \e \leq 
 8^{-1} \, c_1 C^{-1} 
$$ 

\begin{lemma}\label{l.e.22}
We have that
$$
 [D1,D3,\chi1,\e23,\e31]  \to [\e22]  \, .
$$
\end{lemma}
\noindent{\bf Proof.}
By Lemma~\ref{l.tau.ub}, $22D (\log \e^{-1})^{4/3} \leq \e^{-\chimac^2/6}$ is implied by $[\e23,D1,\chi1]$.
The condition $[\e22]$
is thus under these circumstances 
implied by
$$
  \e^{  10^{-11} c_1  D^{3/4}  -  3} \leq 
 8^{-1} \, c_1 C^{-1} \e^{\chimac^2/6}  
$$ 
which in turn is implied by $[D3]$ in the guise
$$
10^{-11} c_1  D^{3/4}  -  3 \geq 2
$$
and $[\e31]$. \qed

Since $c_1 \leq 1$, $[\e31] \to [\e30]$. Note that $[\e29] \to [\e26]$.

\noindent{\em Summary of the analysis of $[\e20,21,22]$.}
To obtain $[\e20,21,22]$ using the preceding three lemmas, we invoke the hypotheses  $[D1,D2,D3,\chi1,\e23,\e25,\e26,\e27,\e28,\e29,\e30,\e31]$. As we have just noted, however, $[\e26,30]$ are redundant. Thus, 
$$
  [D1,D2,D3,\chi1,\e23,\e25,\e27,\e28,\e29,\e31]   \to [\e20,21,22] \, .
$$

We now introduce some further conditions:
$$
 [\e32] \, \,  \e \leq \big( 4k_0 \big)^{-\chi^{-1}}
$$
$$
 [\e33] \, \,  \e \leq \exp \big\{ - 4^3 c^{-9} \big\}
$$
$$
 [\e34] \, \,  \e \leq \exp \bigg\{ - D^{-3} \Big( 78 \Psi_1 \vee 4 \big( (\Psi_2 + 1)^2 + \Psi_3 \big)^{1/2}  \Big)^{3} \bigg\}
$$
$$
 [\e35] \, \,   D \big( \log \e^{-1} \big)^{4/3} \geq 10^9 c_1^{-4/5} C^{1/2}
$$
$$
 [\e36] \, \,  \e \leq \exp \Big\{ - D^{-3/4} 10^{27/4} c_1^{-3/5} C^{3/8} \Big\} 
$$
$$
 [\e37] \, \, \e \leq 10^{-1} D^{-4}
$$
\begin{lemma}\label{l.further.e}
We have the following inferences.
\begin{enumerate}
\item $[\e8] \to [\e1]$
\item $[\e29] \to [\e2,3]$
\item $[\e5,28] \to [\e4]$
\item $[\e6] = [\e32]$
\item $[\e29,D1] \to [\e8]$
\item $[\e13,\e23,D1,\chi1,\e37] \to [\e9]$
\item $[\e29,D1] \to [\e13]$
\item $[\e33,D1] \to [\e14]$
\item $[\e34,D1,\e29] \to [\e15]$
\item $[\e11] \to [\e16]$
\item $[\e35] \to [\e17,18]$
\item $[\e36] = [\e35]$
\item $[\e29] \to [\e19]$
\end{enumerate}
\end{lemma}
\noindent{\bf Proof: (1).}
Due to $[\e1]$ being equivalent to $\e \leq \exp \big\{ - D^{-3} \big\}$.

\noindent{\bf (2).} Trivial.

\noindent{\bf (3).}  $[\e28]$ implies that $\chi \leq 1/2$. We also use $\beta \geq e$, which follows from the definition of $\beta$.

\noindent{\bf (4,5).} Trivial.
 
\noindent{\bf (6).} $[\e9]$ is implied by $[\e13]$ and the bound 
$$
 \e \leq 4^{-2} \Big( 4/3 \cdot D \big( \log \e^{-1} \big)^{1/3} \Big)^{-4} \, .
$$  
This display is equivalent to
$$
  \e \leq 4^{-6} 3^4 D^{-4}  \big( \log \e^{-1} \big)^{4/3} \, .
$$
By Lemma~\ref{l.tau.ub}, $[\e23,D1,\chi1]$ imply that the last display is implied by 
$$
  \e \leq 4^{-6} 3^4 D^{-4}  \cdot 22 D \e^{\chi^2/6} \, .
$$
Invoking $[\chi1]$ (which is $\chi \leq 1$), we see that the last is implied by 
$$
 \e \leq \big( 22 \cdot 4^{-6} 3^4 \big)^{6/5} D^{-18/5} \, .
$$
Since $\big( 22 \cdot 4^{-6} 3^4 \big)^{6/5} = 0.368 \cdots$, we see that, using $D \geq 1$ (i.e., $[D1]$), the last condition is implied by $[\e37]$.
This completes the derivation of $({\bf 6})$.

\noindent{\bf (7).}  $[\e13]$ is equivalent to $\e \leq \exp \big\{ - 6^3 D^{-3} \big\}$ which is implied by $[\e29,D1]$.

\noindent{\bf (8).}  $[\e14]$ is equivalent to $\e \leq \exp \big\{ - ( 3 c^{-3}  + 1  )^3 D^{-3} \big\}$. Using $D \geq 1$ (i.e., $[D1]$) and $c \leq 1$, this is implied by $[\e33]$.

\noindent{\bf (9).}  Note first that the bound $1 \leq 2^{-1} D \big( \log \e^{-1} \big)^{1/3}$ is implied by $[D1,\e29]$.

Thus, when $[D1,\e29]$ holds, $[\e15]$ is implied by 
$$
  D \big( \log \e^{-1} \big)^{1/3} \geq        78 \Psi_1 \vee 4 \big( (\Psi_2 + 1)^2 + \Psi_3 \big)^{1/2}  \, .
$$
This display is equivalent to $[\e34]$.

\noindent{\bf (10).} $[\e11] \to [\e16]$ since $c_3 \leq c_1$.

\noindent{\bf (11).} $[\e35] \to [\e17]$ due to $C \geq 1$.  $[\e35] \to [\e18]$ due to $c_1 \leq 1$.

\noindent{\bf (12).} Trivial.

\noindent{\bf (13).} Since $D \big( \log \e^{-1} \big)^{1/3} \geq 1$, $[\e19]$ is implied by $\log \e^{-1} \geq 2 \cdot 87 \e^{1/2}$. Since $\e \leq 1$ by $[\e29]$, the last inequality is implied by $\e \leq e^{-174}$, which itself is implied by $[\e29]$. \qed

\noindent{\em A summary of progress so far.} We are working to simplify the $[\e]$-list, which consists of the conditions $[\e1,\e2,\cdots,\e20]$. Drawing together our inferences made thus far, we see that this collection of conditions is implied by 
$$
 \big[  \e5,\e10,\e11,\e12,\e23,\e24,\e25,\e27,\e28,\e29,\e31,\e32,\e33,\e34,\e36,\e37,D1,D2,D3,\chi1  \big] \, .
$$ 
A little further simplification will now be made of these conditions, and so it is convenient to restate them all now.

 $$
 [\e5] \, \, \e \leq \exp \big\{ - \chi^{-2} \big\}
 $$
 $$   
 [\e10] \, \,   \e  \leq 
(\eta_0)^{72} 
$$
$$ 
  [\e11] \, \,   \e  \leq 
10^{- 1342}  c_3^{44} \big( D (\log \e^{-1})^{1/3} + 2 \big)^{-8}
$$
$$ 
 [\e12] \, \,   \e  \leq 
\exp \big\{ - 2 C^{3/8} \big\} 
$$
$$
 [\e23] \, \,  \e \leq  \big( 10^{-41} \chi^{36} D^{-14} \big)^{\chi^{-2}} \, , 
   $$
   $$
 [\e24] \, \, \e \leq 2^{-4\chi^{-2}} \, ,
$$ 

$$
 [\e25] \, \,  \e \leq     \exp \big\{ - \beta^{ \chimac^{-1/2}  (3480)^{1/2}  } \big\}
$$
$$
 [\e27] \, \,  \e
  \leq (2 \condee)^{- 3480 \chimac^{-1} ( \log \beta)^{2}  }
  $$
$$
 [\e28] \, \, \chimac   \leq  2^{-1}   \big( 1  + 500( \log \beta)^{2} \big)^{-1}   
$$   
and
$$
 [\e29] \, \,      \e  \leq  \exp \big\{ - 2^{12} \big\} \, .
$$
$$
 [\e31] \, \, \e \leq 
 8^{-1} \, c_1 C^{-1} 
$$ 
$$
 [\e32] \, \,  \e \leq \big( 4k_0 \big)^{-\chi^{-1}}
$$
$$
 [\e33] \, \,  \e \leq \exp \big\{ - 4^3 c^{-9} \big\}
$$
$$
 [\e34] \, \,  \e \leq \exp \bigg\{ - D^{-3} \Big( 78 \Psi_1 \vee 4 \big( (\Psi_2 + 1)^2 + \Psi_3 \big)^{1/2}  \Big)^{3} \bigg\}
$$
$$
 [\e36] \, \,  \e \leq \exp \Big\{ - D^{-3/4} 10^{27/4} c_1^{-3/5} C^{3/8} \Big\} 
$$
$$
 [\e37] \, \, \e \leq 10^{-1} D^{-4}
$$
$$
 [D1] \, \, D \geq 1
$$
$$
 [D2] \, \,  D \geq \rsc^{-1/3} 2^{5}
$$
$$
 [D3] \, \,  D  \geq 5^{4/3} \cdot 10^{44/3} c_1^{-4/3}
$$
    $$
 [\chi1] \, \,   \chi \leq 1
    $$  

\noindent{\em Some further inferences.}

We now introduce a further condition:
$$
 [\e38] \, \,  \e \leq \exp \Big\{ -  10^7 c_1^{-9} C^{3/8} \Big\} 
$$

Note that $[\e38] \to [\e12]$ since $c_1 \leq 1$.  Since $c_1 \leq 1$, it is easily checked that $[\e38] \to [\e31]$. Since $c_1 \leq c$ and $C \geq 1$, $[\e38] \to [\e33]$.
Note also that $[\e38,D1] \to [\e36]$ since $c_1 \leq 1$.

Note that $[\e23,D1,\chi1] \to [\e5,\e24,\e37]$.

Note also that $[\e28] \to [\chi1]$.

Introducing 
$$
 [D4] \, \,  D  \geq 10^{16} c_1^{-4/3}
$$
we see that,
since $5^{4/3} 10^{44/3} = 3.96 \cdots \times 10^{15}$,
$[D4] \to [D3]$. We also have $[D4] \to [D1,2]$ since $c_1 \leq 1$ and $c_1 \leq c$.

Introducing
$$
 [\e39] \, \, \e \leq 10^{-2634} c_3^{88} \, ,
$$
we see, since $\exp \big\{ 2^{12} \big\} = \big[  10^{1778}  ,  10^{1779}  \big]$ and $c_3 \leq 1$, $[\e39] \to [\e29]$.
We also make a 

Claim: $[\e38,39,D4] \to [\e11]$.

\noindent{\bf Proof.} In view of the form of $[\e11]$, this result is valid when 
$$
 \Big(D \big( \log \e^{-1} \big)^{1/3} + 2  \Big) \geq \e^{1/2}
$$
and this bound we now derive. It is equivalent to 
$$
 \e^{-1} \leq \exp \Big\{ D^{-3} (\e^{-4} - 2)^3 \Big\}
$$
and thus implied by $\e^{-1} \leq D^{-3} (\e^{-4} - 2)^3$ which since $\e^{-4} \geq 4$ (valid since $\e \leq 2^{-1/2}$ by $[\e39]$) follows from $\e^{-1} \leq D^{-3} \big( 2^{-1} \e^{-4} \big)^3$ or equivalently from $\e \leq D^{-3/11} 2^{-3/11}$. The latter holds by $[D4,\e38]$. \qed

\noindent{\em Concluding the analysis of the $[\e]$ list.}

In summary, then, the original $[\e]$-list, consisting of $[\e1,\e2,\cdots,\e20]$, is implied by 
\begin{equation}\label{e.esimple}
 \big[  \e10,\e23,\e25,\e27,\e28,\e32,\e34,\e38,\e39,D4  \big] \, .
\end{equation}
These conditions collectively may be expressed as:
\begin{eqnarray*}
  \e  & \leq & \max \Bigg\{ \,
(\eta_0)^{72} \, , \,   \big( 10^{-41} \chi^{36} D^{-14} \big)^{\chi^{-2}} \, , \,   \exp \big\{ - \beta^{ \chimac^{-1/2}  (3480)^{1/2}  } \big\} \, , \,  (2 \condee)^{- 3480 \chimac^{-1} ( \log \beta)^{2}  }
 \, , \,  \big( 4k_0 \big)^{-\chi^{-1}} \, ,    \\
 & & \qquad \quad   \, \, \, \exp \bigg\{ - D^{-3} \Big( 78 \Psi_1 \vee 4 \big( (\Psi_2 + 1)^2 + \Psi_3 \big)^{1/2}  \Big)^{3} \bigg\} \, , \,  \exp \Big\{ -  10^7 c_1^{-9} C^{3/8} \Big\}  \, , \, 10^{-2634} c_3^{88} \,  \Bigg\} \, ,
\end{eqnarray*}
as well as 
$$
\chimac   \leq  2^{-1}   \big( 1  + 500( \log \beta)^{2} \big)^{-1} 
$$ 
and
$$
D  \geq 10^{16} c_1^{-4/3} \, .
$$

\noindent{\bf Analysing the $[n]$ list.}

We begin this analysis by presenting some further conditions. Since the conditions concern $\e$, we label them as such.
$$
 [\e40] \, \,   \e \leq \exp \big\{  - (12)^{-4} (  \log K_0  )^4   \big\}
$$
$$
 [\e41] \, \,    \e \leq \exp \{ - 8! \}
$$
$$
 [\e42] \, \,   \e \leq \exp \big\{ - 2 K_0^2 \big\}
$$
$$
  [\e43] \, \,    K_0^{(12)^{-2} \big( \log \log  (2^{-1} \e^{-\chi} ) \big)^2  } \leq \e^{-1}
$$
$$
 [\e44] \, \,     \big(   \log  (  2^{-1} \e^{-\chimac}  )  \big)^{K_0}  \leq \e^{-1}
$$
$$
  [\e45] \, \,    \big(  \log \e^{-1} \big)^{K_0} \leq \e^{-1}
$$

The next lemma will be used to simplify the $[n]$ list.
\begin{lemma}\label{l.n.simplify}
\begin{enumerate}
\item $[\e40,\e41,\chi1] \to [\e43]$
\item  $[\e42]  \to [\e45]$
\item $[\e45,\chi1] \to [\e44]$
\end{enumerate}
\end{lemma}
\noindent{\bf Proof: (1).} Note that $[\e43]$ is equivalent to
$$
  \big( \log \log  (2^{-1} \e^{-\chi} ) \big)^2 (12)^{-2} \log K_0 \leq \log \e^{-1}
$$
and thus is
implied by 
$$
(12)^{-2} \log K_0 \leq \big( \log \e^{-1} \big)^{1/2}
$$
and
$$
\big( \log \log  (2^{-1} \e^{-\chi} ) \big)^2 \leq \big( \log \e^{-1} \big)^{1/2} \, .
$$
These two displayed equations will here be called $[1]$ and $[2]$.

Note that $[1]$ is equivalent to $[\e40]$.

Set $h = \log \e^{-1}$. Note that $[2]$ is equivalent to
$$
  \big(  \log ( \chi h   - \log 2  )   \big)^4 \leq h
$$
When $[\chi1]$ holds, the last is implied by $(\log h)^4 \leq h$
or $h \leq \exp \{ h^{1/4} \}$
which is implied by $h \geq 8!$ or equivalently $[\e41]$.

\noindent{\bf (2).} Again we set $h = \log \e^{-1}$. Note that $[\e45]$
is equivalent to $h \leq \exp \big\{  h K_0^{-1} \big\}$ 
which is implied by $h \leq h^2 K_0^{-2}/2$
or $h \geq 2K_0^2$ which is $[\e42]$.

\noindent{\bf (3).} Trivial. \qed

We now introduce some further conditions expressed in terms of $n$.
$$
 [n14] \, \, n  \geq 2 \e^{-7/2}
$$
$$
 [n15] \, \,  n \geq a_0^{-9} \e^{-3}
$$
 $$   
 [n16] \, \,   n    \geq   10^{186}       \rsc^{-36}     \e^{-3/2 - 466\chimac}        
$$ 
$$
 [n17] \, \,   n \geq 10^{740}    D^{36}  c_3^{-84}  \e^{-234}  
 $$
$$
 [n18] \, \,  n \geq 2^{16} a_0^{-9}  D^9    \e^{-3}    
$$
$$
 [n19] \, \, n    \geq   
 10^{186} c^{-36}  \e^{-75/2 - 466\chimac} 
$$ 
$$
 [n20]  \, \,  n \geq 10^{35}    c^{-18} D^{36}  \e^{-50} 
$$
$$
 [n21]  \, \,  n \geq 10^{35}    c^{-18} D^{36}  \e^{-198} 
$$
$$
 [n22]  \, \,  n \geq 10^{740}    D^{36}  c_3^{-84}  
 a_0^{-9}    \e^{-504}  
 $$

\begin{lemma}\label{l.n.analysis}
\begin{enumerate}
\item $[n16] \to [n1]$
\item $[n14,\e43,\e44] \to [n2]$
\item  $[n15,\e23,D1,\chi1] \to [n3]$
\item  $[n16,\e23,D1,\chi1] \to [n4]$
\item  $[n17,\e8]  \to [n6]$
\item $[n18,\e8] \to [n7]$
\item $[n16] \to [n9]$
\item $[n19] \to [n10,16]$
\item $[n20,\e8] \to [n11]$
\item $[n20,D1] \to [n12]$
\item $[n21] \to [n13,n20]$
\item $[n22] \to [n14,15,16,17,18,19,20,21]$
\end{enumerate}
\end{lemma}
\noindent{\bf Proof: (1).} $[n1]$ is the condition that $n \geq 8$ which is implied by $[n16]$ since $c \leq 1$ and $\e \leq 1$.

\noindent{\bf (2).} When $[\e43,44]$ hold, $[n2]$ is implied by $n \e^{3/2} \geq 2 \e^{-2}$ and thus by $[n14]$.

\noindent{\bf (3).} By Lemma~\ref{l.tau.ub} and the comments that follow it, $[\e23,D1,\chi1]$ imply that $\tau \leq \e^{-\chi^2/6}$.
Under these circumstances, $[n3]$ is implied by $n \e^{3/2}   \geq    a_0^{-9} \e^{-3\chi^2/2}$.
Since $\chi \leq 1$, it is also implied by $n \geq a_0^9 \e^{-3}$.

\noindent{\bf (4).} For the same reason, $[\e23,D1,\chi1]$ imply that $[n4]$ is implied by 
 $$   
   n \e^{3/2}   \geq   10^{325}       \rsc^{-36}   \big( 2^{-1} \e^{-\chimac} \big)^{465}  \max \big\{   1  \, , \,  \e^{-\chi^{12}}   \big\}   
$$ 
Since $\e \leq 1$ and $\chi \leq 1$, this condition is implied by 
 $$   
   n \e^{3/2}   \geq   10^{325}       \rsc^{-36}   2^{-465}  \e^{-466\chimac}        
$$ 
which since $2^{465} = 9.52 \cdots \times 10^{139}$ is implied by $[n16]$.

\noindent{\bf (5).} Note that $2 \leq  D \big( \log \e^{-1} \big)^{1/3}$ holds when $\e \leq \exp \big\{ - 8 D^{-3} \big\}$ occurs and that the latter bound in $[\e8]$.
Thus, when $[\e8]$ occurs, $[n6]$ is implied by 
$$
 n \geq 2 \cdot 10^{728}     c_3^{-84}  \e^{-222}   \big( 2 D (\log \e^{-1})^{1/3}  \big)^{36} 
$$
Since $\log \e^{-1} \leq \e^{-1}$, the last bound is implied by 
$$
 n \geq 2^{37} 10^{728}    D^{36}   c_3^{-84}  \e^{-234}  
 $$
 Since $2^{37} = 1.37 \times 10^{11}$, the last is implied by $[n17]$.

\noindent{\bf (6).}
For the same reason, $[n7]$ is when $[\e8]$ occurs implied by 
$$
 n \geq 2 a_0^{-9} \big( 2 D (\log \e^{-1})^{1/3}  \big)^9 2^{6}  \, .
$$
This is implied by  $[n18]$.

\noindent{\bf (7,8).} This is due to $\e \leq 1$.

\noindent{\bf (9).} When $[\e8]$ holds, $[n11]$ is implied by 
$$
n \geq 10^{24} c^{-18} \e^{-75/2} \big( 2 D (\log \e^{-1})^{1/3}  \big)^{36}   
$$
and thus by 
$$
n \geq 10^{24} 2^{36}   c^{-18} \e^{-75/2 - 12} D^{36} \, .    
$$
Since $2^{36} =  6.87 \cdots 10^{10}$ and $\e \leq 1$, this is implied by 
$$
n \geq 10^{24} 2^{36}   c^{-18} \e^{-50} D^{36} 
$$
and thus by $[n20]$.

\noindent{\bf (10).} $[n12]$ is implied by 
$$
  D (\log \e^{-1})^{1/3}   \leq    6^{-1} c n^{1/9} 
$$
or equivalently
$$
  n \geq  6^9 c^{-9} D^9 (\log \e^{-1})^3    
$$
When $[D1]$ holds, this is implied by $[n20]$ since $\e \leq 1$.

\noindent{\bf (11).} $[n13]$ is implied by 
$$ 
     D^{36} \big(\log \e^{-1} \big)^{48}  
 3^{-36}  \e^{-150} \leq n   
 $$
 and thus by 
$$ 
  n \geq     D^{36}   3^{-36}  \e^{-198} \, . 
 $$
 This last is implied by $[n21]$ since $c \leq 1$. Note that $[n21] \to [n20]$ since $\e \leq 1$.

\noindent{\bf (12).} This is due to  
 $D \geq 1$, $c_3 \leq 1$, $a_0 \leq 1$, 
$\e \leq 1$, $c_3 \leq c$ and $\chi \leq 1$. (Recall that the condition $a_0 \leq 1$ was imposed when this parameter was introduced in Theorem~\ref{t.maxpoly}.) \qed

\noindent{\em Simplifying the $[n]$ list.}
We learn from Lemma~\ref{l.n.analysis} that the original $[n]$ list, consisting of conditions $[n1,2,\cdots,13]$ is implied by 
$$ 
  \big[ n5, n8,  n22, \e8, \e23,  \e43, \e44, D1, \chi1 \big] \, .
$$
Recall that our simplified $[\e]$ list~(\ref{e.esimple}) consists of $\big[  \e10,\e23,\e25,\e27,\e28,\e32,\e34,\e38,\e39,D4  \big]$. 
Since we are assuming these conditions, we may remove from the displayed list those conditions that are implied by them: 
these include $[D1]$, $[\chi1]$, $[\e8]$ and $[\e23]$. (Note that $[\e8]$ is implied because it is on the original $[\e]$ list, $[\e1,2,\cdots,22]$, and the simplified $[\e]$ list.)
Moreover, by Lemma~\ref{l.n.simplify}, we may replace $[\e43,44]$ by $[\e40,41,42]$.
Thus, our original $[n]$ list is, when the simplified $[\e]$ list is assumed, implied by  
$$ 
  \big[ n5, n8,  n22,    \e40, \e41, \e42  \big] \, .
$$
We may drop $[\e41]$ because it is implied by $[\e11]$ since $c_3 \leq 1$.

Finally, we will drop $[\e40]$ by arguing that it is implied by $[\e42]$. To argue this, recall that the parameter $K_0$, which originates in Theorem~\ref{t.maxpoly}, is supposed to satisfy $K_0 \geq 1$.
If $[\e42]$ fails to imply $[\e40]$, then $2 K_0^2 < (12)^{-4} \big( \log K_0 \big)^4$
or equivalently $\exp \big\{  2^{1/4} 12 K_0^{1/2} \big\} \leq K_0$. This last is not valid when $(4!)^{-1} \big( 2^{1/4} 12 \big)^4 K_0 \geq 1$.
Since $K_0 \geq 1$, this condition is not met, so $[\e42] \to [\e40]$.

Thus, we choose our simplified $[n]$ list to be 
$$ 
  \big[ n5, n8,  n22,     \e42  \big] \, .
$$

\noindent{\bf Conclusion.}
The original collection of conditions, on the $[\e]$, $[n]$, $[D]$ and $[\chi]$ lists, are implied by the simplified $[\e]$ and $[n]$ lists, which consist of the conditions:
$$
 \big[  \e10,\e23,\e25,\e27,\e28,\e32,\e34,\e38,\e39,D4, n5,n8,n22,\e42  \big] \, .
$$

These conditions collectively may be expressed as:
\begin{eqnarray*}
  \e  & \leq & \max \Bigg\{ \,
(\eta_0)^{72} \, , \,   \big( 10^{-41} \chi^{36} D^{-14} \big)^{\chi^{-2}} \, , \,   \exp \big\{ - \beta^{ \chimac^{-1/2}  (3480)^{1/2}  } \big\} \, , \,  \\
 &  & \qquad \qquad (2 \condee)^{- 3480 \chimac^{-1} ( \log \beta)^{2}  }
 \, , \,  \big( 4k_0 \big)^{-\chi^{-1}} \, ,    \,  \exp \big\{ - 2 K_0^2 \big\} \, ,   \, 10^{-2634} c_3^{88} \, ,  \\
  & & \qquad \qquad \qquad
     \exp \bigg\{ - D^{-3} \Big( 78 \Psi_1 \vee 4 \big( (\Psi_2 + 1)^2 + \Psi_3 \big)^{1/2}  \Big)^{3} \bigg\} \, , \,  \exp \Big\{ -  10^7 c_1^{-9} C^{3/8} \Big\} \,  \Bigg\} \, ,
\end{eqnarray*}
$$
    n  \geq  \max \bigg\{ \,  4(K_0)^9 \big( \log \e^{-1} \big)^{K_0} \, , \, 
 c^{-18} \max \Big\{  (\coninit_2 + 1)^9 \, , \,   10^{23} \coninit_1^9 \Big\} \, ,  \,   10^{740}    D^{36}  c_3^{-84}  
 a_0^{-9}    \e^{-504}  \, \bigg\} \, ,
 $$
as well as 
$$
\chimac   \leq  2^{-1}   \big( 1  + 500( \log \beta)^{2} \big)^{-1} 
$$ 
and
$$
D  \geq 10^{16} c_1^{-4/3} \, .
$$
Since these are the hypotheses of Lemma~\ref{l.normalcoal}, we have completed the calculational derivation of this result.

\bibliographystyle{plain}

\bibliography{airy}

\end{document}